%% file: kh_lambda.tex
\documentclass[a4,10pt]{article}
\usepackage{amsmath,amscd,amssymb}

\input{notation}
\input{new_theorem}

\begin{document}

\title{Kobayashi-Hitchin correspondence
 for tame harmonic bundles II}
\author{Takuro Mochizuki}
\date{}
\maketitle

\begin{abstract}
\input{abstract}
\end{abstract}

\tableofcontents

\section{Introduction}
\input{1}
\section{Preliminary}
\input{2}

\section{Ordinary metric and some consequences}
\label{section;06.1.23.20}
\input{3}

\section{Continuity of some families of harmonic metrics}
\label{section;06.2.6.1}
\input{4}

\section{The existence of a pluri-harmonic metric}

\input{5}

\section{Filtered local system}
 \label{section;06.2.3.150}
\input{6}

\input{kh_lambdaref}
\noindent
{\it Address\\
Department of Mathematics,
Kyoto University,
Kyoto 606-8502, Japan.\\
takuro@math.kyoto-u.ac.jp
}\\

\end{document}

%% file: notation.tex
\newcommand{\nbigb}{\mathcal{B}}
\newcommand{\nbigc}{\mathcal{C}}
\newcommand{\nbigd}{\mathcal{D}}
\newcommand{\nbige}{\mathcal{E}}
\newcommand{\nbigf}{\mathcal{F}}

\newcommand{\nbigi}{\mathcal{I}}

\newcommand{\nbigk}{\mathcal{K}}
\newcommand{\nbigl}{\mathcal{L}}
\newcommand{\nbigm}{\mathcal{M}}
\newcommand{\nbign}{\mathcal{N}}
\newcommand{\nbigo}{\mathcal{O}}
\newcommand{\nbigp}{\mathcal{P}}

\newcommand{\nbigs}{\mathcal{S}}
\newcommand{\nbigt}{\mathcal{T}}
\newcommand{\nbigu}{\mathcal{U}}
\newcommand{\nbigv}{\mathcal{V}}

\newcommand{\nbigz}{\mathcal{Z}}

\newcommand{\proj}{\mathbb{P}}
\newcommand{\seisuu}{\mathbb{Z}}
\newcommand{\rnum}{{\boldsymbol Q}}

\newcommand{\cnum}{{\boldsymbol C}}
\newcommand{\real}{{\boldsymbol R}}

\newcommand{\DD}{\mathbb{D}}
\newcommand{\EE}{\mathbb{E}}

\newcommand{\vece}{{\boldsymbol e}}

\newcommand{\vecv}{{\boldsymbol v}}
\newcommand{\vecu}{{\boldsymbol u}}

\newcommand{\vecalpha}{{\boldsymbol \alpha}}
\newcommand{\veca}{{\boldsymbol a}}
\newcommand{\vecb}{{\boldsymbol b}}

\newcommand{\vecc}{{\boldsymbol c}}
\newcommand{\vecd}{{\boldsymbol d}}

\newcommand{\vecomega}{{\boldsymbol \omega}}

\newcommand{\vecF}{{\boldsymbol F}}
\newcommand{\vecE}{{\boldsymbol E}}

\newcommand{\vecn}{{\boldsymbol n}}

\newcommand{\lrarr}{\longrightarrow}

\newcommand{\pf}{{\bf Proof}\hspace{.1in}}
\newcommand{\qed}{\mbox{\rule{1.2mm}{3mm}}}

\def\End{\mathop{\rm End}\nolimits}

\def\Image{\mathop{\rm Im}\nolimits}

\def\Re{\mathop{\rm Re}\nolimits}

\def\Gr{\mathop{\rm Gr}\nolimits}

\def\rank{\mathop{\rm rank}\nolimits}

\def\Gr{\mathop{\rm Gr}\nolimits}
\def\Sym{\mathop{\rm Sym}\nolimits}

\def\Res{\mathop{\rm Res}\nolimits}

\def\degpar{\mathop{\rm par\textrm{-}deg}\nolimits}
\def\parch{\mathop{\rm par\textrm{-}ch}\nolimits}

\def\parchern{\mathop{\rm par\textrm{-}c}\nolimits}

\def\tr{\mathop{\rm tr}\nolimits}

\def\vol{\mathop{\rm dvol}\nolimits}
\def\dvol{\mathop{\rm dvol}\nolimits}

\def\id{\mathop{\rm id}\nolimits}

\def\codim{\mathop{\rm codim}\nolimits}
\def\gap{\mathop{\rm gap}\nolimits}

\def\Irr{\mathop{\rm Irr}\nolimits}

\def\wt{\mathop{\rm wt}\nolimits}

\def\diag{\mathop{\rm diag}\nolimits}

\def\Obv{\mathop{\rm Obv}\nolimits}

\newcommand{\Cpoly}[1]{\nbigc^{poly}_{#1}}

\newcommand{\del}{\partial}
\newcommand{\delbar}{\overline{\del}}

\newcommand{\pardeg}{\degpar}

\newcommand{\nibar}{\underline{2}}

\newcommand{\barz}{\overline{z}}
\newcommand{\zbar}{\barz}

\newcommand{\barlambda}{\overline{\lambda}}
\newcommand{\lambdabar}{\barlambda}

\newcommand{\Hbar}{\overline{H}}
\newcommand{\Hbarepsilon}{\overline{H}_{\epsilon}}
\newcommand{\Abar}{\overline{A}}

\newcommand{\prolong}[1]{{}^{\diamond}{#1}}

\newcommand{\prolongg}[2]{{}_{#1}{#2}}

\newcommand{\DDlambda}{\DD^{\lambda}}
\newcommand{\DDlambdastar}{\DD^{\lambda\,\star}}

\newcommand{\laplacianlambda}

\newcommand{\KMS}{{\mathcal{KMS}}}

\newcommand{\Par}{{\mathcal Par}}

\newcommand{\lefttop}[1]{{}^{#1}}

\def\Gal{\mathop{\rm Gal}\nolimits}

\newcommand{\closedopen}[2]{[#1,#2[}
\newcommand{\openclosed}[2]{]#1,#2]}
\newcommand{\openopen}[2]{]#1,#2[}

\newcommand{\Etilde}{\widetilde{E}}
\newcommand{\vecEtilde}{\widetilde{\vecE}}
\newcommand{\thetatilde}{\widetilde{\theta}}

\newcommand{\htilde}{\widetilde{h}}

\newcommand{\gtilde}{\widetilde{g}}
\newcommand{\nbigitilde}{\widetilde{\nbigi}}
\newcommand{\Xtilde}{\widetilde{X}}
\newcommand{\Dtilde}{\widetilde{D}}
\newcommand{\DDtilde}{\widetilde{\DD}}
\newcommand{\DDtildelambda}{\DDtilde^{\lambda}}
\newcommand{\vtilde}{\widetilde{v}}
\newcommand{\vecvtilde}{\widetilde{\vecv}}
\newcommand{\ftilde}{\widetilde{f}}
\newcommand{\vecV}{{\boldsymbol V}}
\newcommand{\Ftilde}{\widetilde{F}}
\newcommand{\ctilde}{\widetilde{c}}
\newcommand{\vecctilde}{\widetilde{\vecc}}

%% file: new_theorem.tex
\newtheorem{thm}{Theorem}[section]
\newtheorem{cor}[thm]{Corollary}

\newtheorem{rem}[thm]{Remark}
\newtheorem{lem}[thm]{Lemma}
\newtheorem{prop}[thm]{Proposition}
\newtheorem{df}[thm]{Definition}

\newtheorem{condition}[thm]{Condition}

%% file: abstract.tex
Let $X$ be a smooth projective complex variety
with an ample line bundle $L$,
and let $D$ be a simple normal crossing divisor.
We establish the Kobayashi-Hitchin correspondence
between tame harmonic bundles on $X-D$
and $\mu_L$-stable parabolic $\lambda$-flat bundles
with trivial characteristic numbers
on $(X,D)$.
Especially, we obtain the quasiprojective
version of the Corlette-Simpson correspondence
between flat bundles and Higgs bundles.

\vspace{.1in}
\noindent
Keywords:
harmonic bundle,
$\lambda$-connection,
Kobayashi-Hitchin correspondence
 \\
MSC: 14J60, 53C07

%% file: 1.tex
\subsection{Main results}

We explain the main results in this paper.
We do not recall history or background about
the study of Kobayashi-Hitchin correspondence and harmonic bundles,
for which we refer to the introductions of
\cite{s5}, \cite{lubke-teleman} or \cite{mochi4},
for example.
The notion of regular filtered $\lambda$-flat bundles
and parabolic $\lambda$-flat bundles
are explained in Subsection \ref{subsection;06.1.23.10}.
(See also Subsections 3.1--3.2 of \cite{mochi4}.
But, we also use a slightly different notation and terminology,
as is explained in Subsection \ref{subsubsection;06.1.24.10}.)
They are equivalent,
and we will not care about the distinction of them.
The notion of filtered local systems is explained
in Section \ref{section;06.2.3.150}.

\subsubsection{Kobayashi-Hitchin Correspondence}

Let $X$ be a smooth complex projective variety
with an ample line bundle $L$.
Let $D$ be a normal crossing divisor of $X$.
Our main purpose is to show the following theorem.
\begin{thm}
 [Theorem \ref{thm;06.1.23.100},
 Proposition \ref{prop;06.1.23.60},
 Proposition \ref{prop;06.1.23.50}]
 \label{thm;06.1.22.10}
Let $(\vecE_{\ast},\DDlambda)$ be
a regular filtered $\lambda$-flat bundle on $(X,D)$.
We put $E:=\vecE_{|X-D}$.
Then, the following conditions are equivalent.
\begin{itemize}
\item
It is $\mu_L$-polystable
with the trivial characteristic numbers
$\pardeg_{L}(\vecE_{\ast})=
 \int_X\parch_{2,L}(\vecE_{\ast})=0$.
\item
There exists a pluri-harmonic metric $h$
 of $(E,\DDlambda)$
 adapted to the parabolic structure.
\end{itemize}
 Such a metric is unique up to obvious ambiguity.
\hfill\qed
\end{thm}

\begin{rem}
The claims of Theorem {\rm\ref{thm;06.1.22.10}}
in the case $\lambda=0$
has already been proved
in our previous paper {\rm\cite{mochi4}}.
Hence, we restrict ourselves
to the case $\lambda\neq 0$
in this paper.
\hfill\qed
\end{rem}

\begin{cor}
 [Corollary \ref{cor;06.1.23.101}]
Let $\Cpoly{\lambda}$ denote the category
of $\mu_L$-polystable
regular filtered $\lambda$-flat bundles on $(X,D)$
with trivial characteristic numbers.
Then, we have the natural equivalence of 
the categories
$\Cpoly{\lambda_1}\simeq
 \Cpoly{\lambda_2}$
for any $\lambda_i\in\cnum$ $(i=1,2)$.
The equivalence preserves the tensor products,
direct sums and duals.
\hfill\qed
\end{cor}

\begin{rem}
 \label{rem;06.2.5.1}
Let $\lambda_i\in\cnum^{\ast}$ $(i=1,2)$.
A $\lambda_2$-connection
$\DD^{\lambda_2}=d''+(\lambda_2/\lambda_1)\cdot d'$
is induced from a $\lambda_1$-connection $\DD^{\lambda_1}=d''+d'$.
Hence we have the obvious functor
$\Obv:\Cpoly{\lambda_1}\lrarr \Cpoly{\lambda_2}$.
But this is not the same 
as the above functor $\Xi_{\lambda_1,\lambda_2}$.
\hfill\qed
\end{rem}

Especially, we obtain a generalization
of the Corlette-Simpson correspondence
between flat bundles and Higgs bundles
in the so-called non-abelian Hodge theory.
\begin{cor}
We have the equivalences of the following two categories:
\begin{itemize}
\item
The category of $\mu_L$-polystable regular filtered Higgs bundles
on $(X,D)$ with trivial characteristic numbers.
\item
The category of $\mu_L$-polystable regular filtered flat bundles
on $(X,D)$ with trivial characteristic numbers.
\hfill\qed
\end{itemize}
\end{cor}

\subsubsection{Bogomolov-Gieseker inequality 
and some formula for the characteristic numbers}

Let $X$, $L$ and $D$ be as above.

\begin{thm}
 [Corollary \ref{cor;06.1.13.150}]
 \label{thm;06.1.24.1}
Let $(\vecE_{\ast},\DDlambda)$ be a $\mu_L$-stable 
regular filtered $\lambda$-flat bundle on $(X,D)$
in codimension two.
Then, we have the following inequality
holds for the parabolic characteristic numbers
for $\vecE_{\ast}$:
\begin{equation}
 \label{eq;06.1.22.1}
 \int_X\parch_{2,L}(\vecE_{\ast})
\leq
 \frac{\int_X\parchern_{1,L}^2(\vecE_{\ast})}{2\rank E}.
\end{equation}
It is a generalization of the so-called Bogomolov-Gieseker inequality.
\hfill\qed
\end{thm}

In the case $\lambda\neq 0$,
we also have some formulas
about the parabolic Chern characteristic numbers,
which are valid for any parabolic $\lambda$-flat bundles
in codimension two.
One of the formulas can be stated simply,
after we see the correspondence of 
regular filtered $\lambda$-flat sheaves and
filtered local systems.
Let$(\vecE_{\ast},\DDlambda)$ be a regular filtered
$\lambda$-flat sheaf on $(X,D)$.
As is explained in Remark \ref{rem;06.2.5.1},
we have the obvious correspondence
of flat $\lambda$-connection 
$\DDlambda=d''+d'$ $(\lambda\neq 0)$
and flat connection $\DD^{\lambda\,f}=d''+\lambda^{-1}d'$.
In particular, we obtain the local system $\nbigl$ on $X-D$
from the flat bundle $(\vecE_{\ast},\DD^{\lambda,f})_{|X-D}$.
Moreover,
the parabolic structure of $(\vecE_{\ast},\DDlambda)$
induces the filtered structure of $\nbigl$,
and we have the more refined claims
as in the following proposition.

\begin{prop}
[Corollary \ref{cor;06.2.3.100} and Corollary \ref{cor;06.2.3.200}]
 \label{prop;06.2.5.3}
Let $\widetilde{\nbigc}(X,D)$ denote the category
of filtered local system on $(X,D)$,
and let $\nbigc_{\lambda}^{sat}(X,D)$
denote the category of
saturated regular filtered $\lambda$-flat sheaves 
on $(X,D)$ for $\lambda\neq 0$.
Then, we have the equivalent functor 
$\Phi_{\lambda}:
 \widetilde{\nbigc}(X,D)\lrarr \nbigc_{\lambda}^{sat}(X,D)$
such that
$\parchern_1(\nbigl_{\ast})
=\parchern_1\bigl(\Phi_{\lambda}(\nbigl_{\ast})\bigr)$
and
 $\int_X\parch_{2,L}(\nbigl_{\ast})
 =\int_X\parch_{2,L}\bigl(\Phi_{\lambda}(\nbigl_{\ast})\bigr)$.
The functor $\Phi_{\lambda}$ preserves the $\mu_L$-stability.
\hfill\qed
\end{prop}

\begin{rem}
From Theorem {\rm\ref{thm;06.1.24.1}}
and Proposition {\rm\ref{prop;06.2.5.3}},
we obtain the Bogomolov-Gieseker inequality
for $\mu_L$-stable filtered local systems
(Corollary {\rm\ref{cor;06.2.5.2}}).
Such a kind of the inequality is discussed
in {\rm{\cite{s8}}}.
\hfill\qed
\end{rem}

\begin{rem}
Let us describe the formula
$\int_X\parch_{2,L}(\nbigl_{\ast})
=\int_X\parch_{2,L}(\Phi(\nbigl_{\ast}))$
in terms of the $\vecc$-truncation
$(\prolongg{\vecc}{E}_{\ast},\DDlambda)$
of saturated regular filtered $\lambda$-flat bundle
$\Phi_{\lambda}(\nbigl_{\ast})$.
For simplicity, we assume $\dim X=2$.
\begin{multline}
\label{eq;06.1.22.3}
 \int_X\parch_{2}(\prolongg{\vecc}{E}_{\ast})=
\frac{1}{2}
 \sum_{i\in S}\sum_{u\in\KMS(\prolongg{\vecc}{E}_{\ast},i)}
\bigl(
 \Re(\lambda^{-1}\alpha)+a
 \bigr)^2\cdot
 r(i,u)\cdot (D_i,D_i)\\
+\frac{1}{2}\sum_{i\in S}
 \sum_{\substack{j\neq i\\ P\in D_i\cap D_j}}
 \sum_{(u_i,u_j)\in\KMS(\prolongg{\vecc}{E}_{\ast},P)}
 \bigl(\Re (\lambda^{-1}\alpha_i)+a_i\bigr)
 \bigl(\Re(\lambda^{-1}\alpha_j)+a_j\bigr)
\cdot r(P,u_i,u_j).
\end{multline}
Here, 
$u=(a,\alpha)$, $u_i=(a_i,\alpha_i)$ and $u_j=(a_j,\alpha_j)$
denote the KMS-spectra of 
$(\prolongg{\vecc}{E},\DDlambda)$,
which are elements of
$\real\times\cnum$.
We put 
$r(i,u):= \rank
 \lefttop{i}\Gr^{F,\EE}_{u}(\prolongg{\vecc}{E})$
for $u\in\KMS(\prolongg{\vecc}{E}_{\ast},i)$,
and $r(P,u_i,u_j):=\rank
 \lefttop{P}\Gr^{F,\EE}_{(u_i,u_j)}(\prolongg{\vecc}{E}_{|P})$
for $(u_i,u_j)\in \KMS(\prolongg{\vecc}{E},P)$ and $P\in D_i\cap D_j$.
And $(D_i,D_j)$ and $\bigl(D_i,c_1(L)\bigr)$
denote the intersection numbers.

We also have some other formulas 
for $\int_X\parch_{2}\bigl(\prolongg{\vecc}{E}_{\ast}\bigr)$
(Proposition {\rm\ref{prop;06.1.13.70}})
or some vanishings for the data of
$(\prolongg{\vecc}{E}_{\ast},\DDlambda)$ at $D$
(Corollary {\rm\ref{cor;06.1.13.150}}
and Proposition {\rm\ref{prop;06.1.13.70}}).
\hfill\qed
\end{rem}

\subsubsection{Vanishing of the characteristic numbers
 and existence of the Corlette-Jost-Zuo metric}

Due to Proposition \ref{prop;06.2.5.3},
we obtain the vanishings
$\pardeg_L(\vecE_{\ast})=\int_X\parch_{2,L}(\vecE_{\ast})=0$,
when $(\vecE_{\ast},\nabla)$ corresponds
to the filtered local system whose parabolic structure is trivial,
in other words,
$\Re(\alpha)+a=0$ is satisfied
for any KMS-spectrum $u=(a,\alpha)\in\KMS(i)$
and for any $i\in S$.
We can apply such a consideration
to the canonical prolongation
of a flat bundle due to P. Deligne \cite{d}.
Let $(E,\nabla)$ be a flat bundle on $X-D$.
Then, it is shown that there exists
the holomorphic vector bundle
$\widetilde{E}$ on $X$ satisfying
(i) $\widetilde{E}_{|X-D}=E$,
(ii) $\nabla \widetilde{E}\subset 
 \widetilde{E}\otimes\Omega_X^{1,0}(\log D)$,
(iii) the real parts of the eigenvalues of 
$\Res_i(\nabla)$ are contained in $\closedopen{0}{1}$.
In that case, we have the naturally defined parabolic 
structure $\vecF$ for which $\Re(\alpha)+a=0$ are satisfied
for any KMS-spectrum $(a,\alpha)$.
Hence, we obtain the vanishing
$\pardeg_{L}(\widetilde{E},\vecF)
=\int_X\parch_{2,L}(\widetilde{E},\vecF)=0$.

This vanishing is significant to understand
the existence theorem of the Corlette-Jost-Zuo metric
from the view point of Kobayashi-Hitchin correspondence.
When $(E,\nabla)$ is semisimple,
we know the existence of
a tame pure imaginary pluri-harmonic metric,
which we call the Corlette-Jost-Zuo metric.
(See \cite{corlette} for the case $D=\emptyset$
and \cite{JZ2} for the general case.
See also \cite{mochi2}.)
Since semisimplicity obviously implies the $\mu_L$-polystability
of $(\widetilde{E},\vecF,\nabla)$
(\cite{sabbah2}, for example),
we can derive the existence of the Corlette-Jost-Zuo metric
from Theorem \ref{thm;06.1.22.10}
due to the vanishing of the characteristic numbers.

\subsection{Methods and difficulty}

\subsubsection{Perturbation of parabolic structure}
\label{subsubsection;06.1.25.1}

Let $X$ be a smooth projective surface,
and let $D$ be a simple normal crossing divisor of $X$.
Let $(E,\vecF,\DDlambda)$ be a parabolic
$\lambda$-flat bundle on $(X,D)$.
For any small $\epsilon>0$,
we take an $\epsilon$-perturbation $\vecF^{(\epsilon)}$
of the parabolic structure,
and then $(E,\vecF^{(\epsilon)},\DDlambda)$
is {\em graded semisimple}
(Subsection \ref{subsubsection;06.1.19.10}).
It can be shown that
the pseudo curvature of ordinary metrics
for $(E,\vecF^{(\epsilon)},\DDlambda)$ $(\epsilon>0)$
satisfy the appropriate finiteness
(Section \ref{section;06.1.23.20}).
By using the theorem of Simpson,
we can take a Hermitian-Einstein metric
$h_{HE}^{(\epsilon)}$
of $(E_{|X-D},\DDlambda)$ which is adapted
to $\vecF^{(\epsilon)}$ $(\epsilon>0)$.
Then, we can easily derive the Bogomolov-Gieseker inequality
(Theorem \ref{thm;06.1.24.1}).
We also obtain the formulas 
by calculating the integrals of the characteristic numbers
for pseudo curvatures,
for example (\ref{eq;06.1.22.3}).

Let us consider the existence of a pluri-harmonic metric
(Theorem \ref{thm;06.1.22.10}).
Ideally,
the limit $\lim_{\epsilon\to 0}h^{(\epsilon)}_{HE}$
should give the desired pluri-harmonic metric
for the given flat parabolic bundle $(E,\vecF,\DDlambda)$.
However, it is not easy to show such a convergence.
It is the main problem
which we have to overcome in this paper.

\subsubsection{Difficulty}
\label{subsubsection;06.1.24.5}

In \cite{mochi4},
we gave an argument to deal with such a convergence
problem for the case $\lambda=0$.
The argument doesn't work in the case $\lambda\neq 0$.
Let us explain what is the difference heuristically
and imprecisely
in the case $\lambda=1$.
Since we have $\pardeg_L(E,\vecF^{(\epsilon)})=0$,
the metrics $h_{HE}^{(\epsilon)}$ give
the harmonic metrics in this case.
Recall that a harmonic metric 
can be regarded as a harmonic map,
at least locally,
and that we know a well established argument
for the convergence of a sequence of harmonic maps
when the energies are dominated (\cite{ES}).
In our case,
the energies of $h_{HE}^{(\epsilon)}$
over $X-D$ are not finite, in general.
Even if we consider the energies over a compact subset
$Z\subset X-D$,
it is not clear how to derive a uniform estimate
which is independent of $\epsilon$.
On the other hand,
the Higgs field is fixed 
for such a convergence problem in the case $\lambda=0$.
In particular, the eigenvalues of the Higgs fields are fixed.
Then, we can derive the estimate
of the local $L^2$-norm of the Higgs fields
independently from $\epsilon$.
Since such $L^2$-norms play the role of the energies,
the local convergence can be easily shown
in the Higgs case,
although we need some technical argument
for global convergence.
On the contrary,
even the local convergence is not
easy to show in the case $\lambda\neq 0$.

\subsubsection{Convergences}

To attack the problem,
we discuss similar convergence problems in the curve case
where the Kobayashi-Hitchin correspondence
was established and well understood by 
the work of C. Simpson \cite{s2}.
Let $C$ be a smooth projective curve,
and let $D$ be a divisor of $C$.
Let $(E,\vecF,\DDlambda)$ be a $\lambda$-flat
stable parabolic bundle on $(C,D)$,
and let $\vecF^{(\epsilon)}$ be $\epsilon$-perturbations.
Note we have
$\det(E,\vecF,\DDlambda)=\det(E,\vecF^{(\epsilon)},\DDlambda)$.
We can take a sequence of harmonic metrics
$h^{(\epsilon)}$ for $(E,\vecF^{(\epsilon)},\DDlambda)$
$(\epsilon\geq 0)$
such that $\det h^{(\epsilon)}=\det h^{(0)}$,
due to the result of Simpson.

First, we will show that
the sequence $\{h^{(\epsilon)}\,|\,\epsilon>0\}$
converges to $h^{(0)}$.
Namely,
let $h_{in}^{(\epsilon)}$ $(\epsilon>0)$
be initial metrics
for $(E,\vecF^{(\epsilon)},\DDlambda)$,
and let $s^{(\epsilon)}$ be the endomorphism
determined by
$h^{(\epsilon)}=h_{in}^{(\epsilon)}\cdot s^{(\epsilon)}$.
Then, we can show the following relations:
\begin{equation}
 M(h^{(\epsilon)}_{in},h^{(\epsilon)})\leq 0,
\quad
 \bigl| \log s^{(\epsilon)}
 \bigr|_{h_{in}^{(\epsilon)}}
\leq C_{1,\epsilon}+C_{2,\epsilon}
\cdot M(h^{(\epsilon)}_{in},h^{(\epsilon)}),
\quad
\bigl\|\DDlambda s^{(\epsilon)}\bigr\|^2
 _{L^2,h^{(\epsilon)}_{in},\omega_{\epsilon}}
 \leq
 \int
 \bigl|\tr \bigl(
s^{(\epsilon)}\cdot G(h^{(\epsilon)}_{in})
 \bigr)\bigr|\dvol_{\omega_{\epsilon}}.
\end{equation}
Here, $M(h^{(\epsilon)}_{in},h^{(\epsilon)})$
denote the Donaldson functionals,
and $\omega_{\epsilon}$ denote appropriate
metrics of $C-D$.
Hence, if we show that
$C_{i,\epsilon}$ can be taken independently from $\epsilon$
for some $\omega_{\epsilon}$,
and if we can construct appropriate family of initial metrics
$h_{in}^{(\epsilon)}$ such that
$G(h^{(\epsilon)}_{in})$ are uniformly bounded
with respect to $\omega_{\epsilon}$
and $h^{(\epsilon)}_{in}$,
then we obtain the $L_1^2$-boundedness
of the family $\{s^{(\epsilon)}\}$.
Then, by using a standard bootstrapping argument,
we can show that the sequence 
$\{s^{(\epsilon)}\}$ is convergent
to the identity in the $C^{\infty}$-sense,
i.e.,
$\{h^{(\epsilon)}\}$ is convergent to $h^{(0)}$
(Section \ref{section;06.2.6.1}).

\vspace{.1in}

Next, suppose that we are given hermitian metrics
$\widetilde{h}^{(\epsilon)}:=
 h^{(\epsilon)}\cdot\widetilde{s}^{(\epsilon)}$
for $\epsilon>0$,
with the following properties:
\begin{itemize}
\item
 $\det\widetilde{h}^{(\epsilon)}=\det h^{(\epsilon)}$.
\item
 $\int |G(\widetilde{h}^{(\epsilon)})|^2\lrarr 0$.
\item
 $\bigl\|\DDlambda s^{(\epsilon)}\bigr\|^2<\infty$.
 (We do not need uniform bound.)
\end{itemize}
Then, we can show that
$\{\widetilde{h}^{(\epsilon)}\}$ is convergent
to $h^{(0)}$.
(See Subsection \ref{subsection;06.1.21.31}
 for more precise claims.)

\vspace{.1in}

We apply the above results
to our convergence problem explained
in Subsection \ref{subsubsection;06.1.25.1}.
Due to the standard Mehta-Ramanathan type theorem
(Proposition \ref{prop;06.1.18.6}),
the restriction
$(E,\vecF,\DDlambda)_{|C}$ is also
stable for almost every very ample $C\subset X$.
Let $h_C$ be a harmonic bundle
of $(E,\vecF,\DDlambda)_{|C}$.
Then, we can show that
$\bigl\{h_{HE\,|\,C}^{(\epsilon)}\bigr\}$ 
is convergent to $h_C$ almost everywhere on $C$
for almost every very ample $C\subset X$,
by using the above result.
Therefore, we obtain a metric $h_{\nbigv}$
defined almost everywhere on $X-D$ such that
 $h_{\nbigv\,|\,C}=h_C$
 almost everywhere on $C$
 for almost every curve $C\subset X$.
With some more additional argument,
we can show that $h_{\nbigv}$ gives
the desired pluri-harmonic metric, indeed
(Subsection \ref{subsection;06.2.12.1}).

\begin{rem}
Perhaps,
the argument of this paper can be used
in the Higgs case,
to show the existence of a pluri-harmonic metric.
However, we remark that
the argument for a convergence given in
{\rm\cite{mochi4}}
can be applied in a wider range.
In fact, we used it to discuss the convergence
of a family of harmonic bundles 
induced by the constant multiplication of Higgs fields.
\hfill\qed
\end{rem}

\subsection{Acknowledgement}

This paper is a result of an effort to understand
the works of C. Simpson, in particular, \cite{s1} and \cite{s5}.
The author is grateful to
N. Budur, H. Konno,
D. Panov and C. Sevenheck for some discussions.
The author thanks A. Ishii and Y. Tsuchimoto 
for their constant encouragement.
He is grateful to the colleagues of Department of Mathematics
at Kyoto University for their cooperation.
The author wrote the first version of this paper
during his stay at Max-Planck Institute for Mathematics.
He acknowledges their excellent hospitality and support.

%% file: 2.tex
\subsection{Generality of regular filtered $\lambda$-flat
 sheaf in complex geometry}

\label{subsection;06.1.23.10}

The notion of a parabolic bundle, filtered bundle
and their characteristic numbers
are explained in Sections 3.1--3.2 of \cite{mochi4}.
We use the notation there.

\subsubsection{$\lambda$-connection}

Let $Y$ be a complex manifold,
and let $\nbige$ be an $\nbigo_Y$-module.
Recall that a $\lambda$-connection of $\nbige$
is defined to be a map
$\DDlambda:\nbige\lrarr\nbige\otimes\Omega^{1,0}_Y$
satisfying the twisted Leibniz rule
$\DDlambda(f\cdot s)=
 f\cdot \DDlambda(s)+\lambda\cdot d_Y(f)\cdot s$,
where $f$ and $s$ denote holomorphic sections
of $\nbigo_Y$ and $\nbige$ respectively.
The maps
$\DDlambda:
 \nbige\otimes\Omega^{p,0}\lrarr 
 \nbige\otimes\Omega^{p+1,0}$
are induced.
When $\DDlambda\circ\DDlambda=0$ is satisfied,
it is called flat.

Let $X$ be a complex manifold,
and let $D$ be a simple normal crossing divisor 
with the irreducible decomposition $D=\bigcup_{i\in S}D_i$.
Let $\nbige_{\ast}=
\bigl(\nbige,\{\lefttop{i}\nbigf\,\big|\,i\in S\}\bigr)$
be a $\vecc$-parabolic sheaf on $(X,D)$
for some $\vecc\in\real^S$.
A flat logarithmic $\lambda$-connection $\DDlambda$
of $\nbige_{\ast}$ is defined to be a map
$\DDlambda:\nbige\lrarr\nbige\otimes\Omega^{1,0}(\log D)$
satisfying the same twisted Leibniz rule as above,
the flatness
$\DDlambda\circ\DDlambda=0$
and $\DDlambda(\lefttop{i}\nbigf_a)\subset
 \lefttop{i}\nbigf_a\otimes\Omega^{1,0}(\log D)$.
Such a tuple $(\nbige_{\ast},\DDlambda)$ will be called
a regular parabolic $\lambda$-flat sheaf.
When the underlying $\vecc$-parabolic sheaf
$\nbige_{\ast}$ is a $\vecc$-parabolic bundle
in codimension $k$,
it is called a regular $\lambda$-flat 
$\vecc$-parabolic bundle in codimension $k$.

\begin{rem}
We often omit to state ``regular'' in this paper,
because we always assume regularity.
Non-regular case is discussed in
{\rm\cite{mochi7}}.
\hfill\qed
\end{rem}

Let $\vecE_{\ast}=
 \bigl(\vecE,\{\prolongg{\vecc}{E}\}\,
 \big|\,\vecc\in\real^S\bigr)$
be a filtered sheaf on $(X,D)$.
A regular $\lambda$-connection of
$\vecE_{\ast}$ is a $\lambda$-connection
$\DDlambda$ of $\vecE$ satisfying
$\DDlambda\bigl(\prolongg{\vecc}{E}\bigr)
 \subset
 \prolongg{\vecc}{E}\otimes\Omega_X^{1,0}(\log D)$.
A tuple $(\vecE_{\ast},\DDlambda)$ is called
a regular filtered $\lambda$-flat sheaf.
When the underlying filtered sheaf
is a filtered bundle in codimension $k$,
it is called a regular filtered $\lambda$-flat bundle
in codimension $k$.

\begin{lem}
 \label{lem;06.2.4.1}
A regular filtered $\lambda$-flat sheaf on $(X,D)$ is
a regular filtered $\lambda$-flat bundle
in codimension one.
\end{lem}
\pf
We have only to check that
there exists a subset $W\subset D$
with $\codim_X(W)\geq 2$,
such that $\prolongg{\vecc}{E}_{\ast\,|\,X\setminus W}$
is a $\vecc$-parabolic bundle
on $(X\setminus W,D\setminus W)$
for some $\vecc$.
We can take $W$ as
$\bigcup_{i\neq j}D_i\cap D_j\subset W$,
and hence we may assume $D$ is smooth.
Since $E=\vecE_{|X-D}$ is locally free
and $\prolongg{\vecc}{E}$ is torsion-free,
we can take $W'\subset D$
with $\codim_X(W')\geq 2$
such that $\prolongg{\vecc}{E}_{|X-W'}$ is locally free.
We may also take a subset
$W''\subset D\setminus W'$ 
with $\codim_X(W'')\geq 2$
such that the parabolic filtration
of $\prolongg{\vecc}{E}_{|D\setminus (W'\cup W'')}$
is filtration in the category of vector bundles.
Then, $W=W'\cup W''$ gives the desired subset.
\hfill\qed

\vspace{.1in}

When $X$ is an $n$-dimensional projective variety
with an ample line bundle $L$,
we can define the $\mu$-stability, $\mu$-semistability,
and $\mu$-polystability of 
regular filtered $\lambda$-flat sheaves
with respect to $L$,
in the standard manner.
``$\mu$-stability with respect to $L$''
will be called $\mu_L$-stability, in this paper.

\subsubsection{KMS-structure}
\label{subsubsection;06.1.26.1}

Let $X$ be a complex manifold,
and let $D$ be a simple normal crossing divisor
with the irreducible decomposition
$D=\bigcup_{i\in S}D_i$.
Let $(\vecE_{\ast},\DDlambda)$ be a regular filtered
$\lambda$-flat bundle in codimension one over $(X,D)$.
For simplicity,
we consider only the case $\lambda\neq 0$.
Let us take any element $\vecc\in\real^S$,
and the $\vecc$-truncation
$\prolongg{\vecc}{E}_{\ast}$ of $\vecE_{\ast}$.
We would like to recall the KMS-structure at $D_i$,
or more precisely, at the generic point of $D_i$.
We may assume that $(\prolongg{\vecc}{E}_{\ast},\DDlambda)$
is a $\vecc$-parabolic bundle.
We have the induced filtration $\lefttop{i}F$ on
$\prolongg{\vecc}E_{|D_i}$.
We put 
$ \lefttop{i}\Gr^F_a(\prolongg{\vecc}{E})
:=\lefttop{i}F_a\bigl(\prolongg{\vecc}{E}\bigr)
 \big/\lefttop{i}F_{<a}\bigl(\prolongg{\vecc}{E}\bigr)$.
Recall that we use the notation:
\[
  \Par(\prolongg{\vecc}{E}_{\ast},i):=
 \bigl\{
 a\,\big|\,c_i-1<a\leq c_i,\,\,\,
 \lefttop{i}\Gr^{F}_a(\prolongg{\vecc}{E})\neq 0
 \bigr\},
\quad
 \Par(\vecE_{\ast},i):=\bigcup_{\vecc\in\real^S}
 \Par(\prolongg{\vecc}{E}_{\ast},i)
\]

Due to the regularity,
we have the residue endomorphism
$\Res_i(\DDlambda)$ 
on $\prolongg{\vecc}{E}_{|D_i}$,
which preserves the filtration $\lefttop{i}F$,
and hence we have the induced endomorphism
$\Gr^F\Res_i(\DDlambda)$
of $\lefttop{i}\Gr^F\bigl(\prolongg{\vecc}{E}\bigr)$.
We remark that the eigenvalues of $\Res_i(\DDlambda)$
are constant on $D_i$.
In particular,
we obtain the generalized eigen decomposition:
\[
 \lefttop{i}\Gr^{F}_a(\prolongg{\vecc}{E})
=\bigoplus _{\alpha\in\cnum}
 \lefttop{i}\Gr^{F,\EE}_{a,\alpha}(\prolongg{\vecc}{E}).
\]
We put 
$\KMS\bigl(
 \prolongg{\vecc}{E}_{\ast},i
 \bigr)
:=\bigl\{
 (a,\alpha)\in\openclosed{c_i-1}{c_i}\times\cnum\,\big|\, 
 \lefttop{i}\Gr^{F,\EE}_{a,\alpha}
 (\prolongg{\vecc}{E}_{|D_i})\neq 0
 \bigr\}$.
Each element of
$\KMS\bigl(
 \prolongg{\vecc}{E}_{\ast},i \bigr)$ or 
$\KMS\bigl(\vecE_{\ast},i\bigr):=
 \bigcup_{\vecc\in\real^S} 
 \KMS\bigl(\prolongg{\vecc}{E}_{\ast},i\bigr)$
is called a KMS-spectrum.

\subsubsection{Prolongment of flat subbundle
and Mehta-Ramanathan type theorem}
\label{subsubsection;06.2.4.30}

To begin with,
we recall a well known fact
about regular singularity of a connection.

\begin{lem}
Let $E$ be a holomorphic bundle on a disc $\Delta$,
and let $\nabla$ be a logarithmic connection of $E$
on $(\Delta,O)$, i.e.,
$\nabla(E)\subset E\otimes\Omega_{\Delta}^{1,0}(\log O)$.
Let $f$ be a flat section of $E_{|\Delta^{\ast}}$.
Then, $f$ naturally gives a meromorphic section of $E$.
\hfill\qed
\end{lem}

\begin{cor}
We put 
$X=\Delta_z\times\Delta_w^{n}$
and $D=\{0\}\times\Delta_w^n$.
Let $E$ be a holomorphic vector bundle
on $X$
and $\nabla$ be the logarithmic connection of $E$
on $(X,D)$.
Let $e$ be a flat section of $E_{|X-D}$.
\begin{itemize}
\item $e$ gives a meromorphic section of $E$.
\item Assume that $e$ is holomorphic on $E$
 and that $e_{|Q}\neq 0$ for some $Q\in D$.
 Then, $e_{|Q'}\neq 0$ for any $Q'\in D$.
\end{itemize}
\end{cor}
\pf
We may assume that we have a holomorphic frame
$\vecv$ of $E$.
We have the expression $e=\sum f_i(z,w)\cdot v_i$.
When we fix $w$,
then $f_{i}(z,w)$ are meromorphic with respect to $z$.
Thus, we have the least integer
$j(w)$ such that the orders of the poles 
of $f_{i}(z,w)$ are less than $j(w)$.
We put $\nbigs_j:=\{w\,|\,j(w)\leq j\}$.
We have $D=\bigcup_{j}\nbigs_j$.
If $\nbigs_j\neq D$, the measure of $\nbigs_j$ is $0$.
Hence, we obtain $\nbigs_j=D$ for some $j$,
which means $e$ is meromorphic.
Thus, we obtain the first claim.

Assume that $e$ is holomorphic
and that $e_{|Q}\neq 0$ for some $Q\in D$.
Recall that we have
the induced connection $\lefttop{D}\nabla$ of $E_{|D}$.
Namely, for any holomorphic section $f\in E_{|D}$,
take a holomorphic $F\in E$ such that $F_{|D}=f$,
and then $\lefttop{D}\nabla(f):=\nabla(F)_{|D}$ is
well defined.
Since we have $\lefttop{D}\nabla (e_{|D})=0$,
we obtain the second claim.
\hfill\qed

\begin{cor}
 \label{cor;06.1.18.1}
We put $X=\Delta^n$,
$D_i=\{z_i=0\}$ and $D=\bigcup_{i=1}^n D_i$.
Let $(E,\nabla)$ be a logarithmic connection
on $(X,D)$,
and let $e$ be a flat section on $X-D$.
\begin{itemize}
\item
$e$ gives a meromorphic section of $E$.
\item
Assume that $e$ is holomorphic.
We put $D_i^{\circ}:=D_i\setminus \bigcup_{j\neq i}D_j$.
If $e_{|Q}\neq 0$ for some $Q\in D_i^{\circ}$,
we have $e_{|Q'}\neq 0$ for any $Q'\in D_i^{\circ}$.
\hfill\qed
\end{itemize}
\end{cor}

Let $X$ be a complex manifold,
and let $D$ be a normal crossing divisor of $X$.
Let $(E,\nabla)$ be a flat bundle on $X-D$.
Recall that P. Deligne gave the extension
$\widetilde{E}$ of $E$
in \cite{d},
such that
(i) $\widetilde{E}_{|X-D}=E$,
(ii) $\nabla (\widetilde{E})
 \subset \widetilde{E}\otimes\Omega^{1,0}(\log D)$,
(iii) the real parts of the eigenvalues
of $\Res_i(\nabla)$ are contained in $\{0\leq t<1\}$.
Such an extension is unique,
or in other words, it is unique as the subsheaf of 
$\iota_{\ast}E$, where $\iota$ denotes the inclusion
$X-D\lrarr X$.
The prolongment can also be done for 
$\lambda$-flat bundle $(E,\DDlambda)$ on $X-D$,
or more precisely, for the associated
flat bundle $(E,\DD^{\lambda\,f})$.

\begin{lem}
 \label{lem;06.2.4.10}
Let $(\vecE_{\ast},\DDlambda)$ be a regular filtered 
$\lambda$-flat bundle on $(X,D)$,
and we put $(E,\DDlambda):=(\vecE_{\ast},\DDlambda)_{|X-D}$.
Let $(\widetilde{E},\DDlambda)$ be the Deligne extension
of $(E,\DDlambda)$.
Then, we have 
$\vecE=\widetilde{E}\otimes\nbigo_X(\ast D)$,
where $\nbigo_X(\ast D)$ denotes the sheaf of
meromorphic functions on $X$ whose poles
are contained in $D$.
\end{lem}
\pf
We have the naturally defined flat section $s$
on $Hom(\prolongg{\vecc}{E},\widetilde{E})_{|X-D}$.
Due to Corollary \ref{cor;06.1.18.1},
$s$ is a meromorphic section,
and hence we obtain the flat inclusion
$\prolongg{\vecc}{E}\lrarr
 \widetilde{E}\otimes\nbigo(N\cdot D)$
for some large integer $N$,
which induce the morphism
$\vecE=\bigcup \prolongg{\vecc}{E}
=\prolongg{\vecc}{E}\otimes\nbigo(\ast D)
\lrarr \widetilde{E}\otimes\nbigo(\ast D)$.
Similarly, we obtain the inclusion
$\widetilde{E}\lrarr 
 \prolongg{\vecc}{E}\otimes\nbigo(N\cdot D)$,
and $\widetilde{E}\otimes\nbigo(\ast D) \lrarr \vecE$.
They are clearly mutually inverse.
\hfill\qed

\begin{lem}
Let $(\vecE_{\ast},\DDlambda)$ be a regular filtered
$\lambda$-flat sheaf on $(X,D)$,
and let $(\widetilde{E},\DDlambda)$ be 
as in the previous lemma.
Then, we have 
$\vecE\simeq\widetilde{E}\otimes\nbigo(\ast D)$
naturally.
\end{lem}
\pf
Due to Lemma \ref{lem;06.2.4.1} and
Lemma \ref{lem;06.2.4.10},
there exists a subset $W\subset D$
with $\codim_X(W)\geq 2$
such that
$\vecE_{|X-W}\simeq 
 \widetilde{E}\otimes\nbigo(\ast D)_{|X-W}$.
Let us fix $\vecc$.
There exists a large integer $N$
such that we have
$\prolongg{\vecc}{E}_{|X-W}\subset 
 \widetilde{E}\otimes\nbigo(N\cdot D)_{|X-W}$.
Since $\widetilde{E}$ is locally free,
we obtain
$\prolongg{\vecc}{E}\subset 
 \widetilde{E}\otimes\nbigo(N\cdot D)$,
and thus
$\vecE\subset \widetilde{E}\otimes\nbigo(\ast D)$.
On the other hand,
there exists a large integer $N'$
such that $\widetilde{E}_{|X-W}\subset
 \prolongg{\vecc}{E}\otimes\nbigo(N'\cdot D)_{|X-W}$.
Hence,
$\widetilde{E}\subset
 \prolongg{\vecc}{E}^{\lor\lor}
\otimes\nbigo(N'\cdot D)$,
where $\prolongg{\vecc}{E}^{\lor\lor}$
denotes the double dual of $\prolongg{\vecc}{E}$.
Hence, we obtain
$\widetilde{E}\otimes\nbigo(\ast D)
\subset \prolongg{\vecc}{E}^{\lor\lor}\otimes\nbigo(\ast D)$.
It is easy to see
$\prolongg{\vecc}{E}^{\lor\lor}
\otimes\nbigo(\ast D)
\simeq
 \prolongg{\vecc}{E}\otimes\nbigo(\ast D)$.
Thus we are done.
\hfill\qed

\begin{lem}
 \label{lem;06.2.4.20}
Let $(\vecE_{\ast},\DDlambda)$ be 
a regular filtered $\lambda$-flat sheaf on $(X,D)$,
and we put $(E,\DDlambda):=(\vecE_{\ast},\DDlambda)_{|X-D}$.
Let $E'$ be a $\lambda$-flat subbundle of $E$.
Then, we have the corresponding 
regular filtered $\lambda$-flat subsheaf 
$\vecE'_{\ast}\subset\vecE_{\ast}$
such that $\prolongg{\vecc}{E'}$ are saturated
in $\prolongg{\vecc}{E}$.
\end{lem}
\pf
Let $\widetilde{E}$ denote the Deligne extension
of $(E,\DDlambda)$.
We have the corresponding subbundle
$\widetilde{E}'\subset\widetilde{E}$.
Therefore, we obtain
$\widetilde{\vecE}':=\widetilde{E}'\otimes\nbigo(\ast D)
\subset \widetilde{E}\otimes\nbigo(\ast D)
=\vecE$.
For each $\vecc$,
the $\vecc$-truncation
$\prolongg{\vecc}{E'}$ is given by the intersection
of $\prolongg{\vecc}{E}$ and $\vecE'$
in $\vecE$.
Or equivalently,
$\prolongg{\vecc}{E'}$ can be given by
the intersection of $\prolongg{\vecc}{E}$
and $\widetilde{E'}(N\cdot D)$
in $\widetilde{E}(N\cdot D)$ for sufficiently large $N$.
Thus, we obtain
$\vecE'_{\ast}\subset\vecE_{\ast}$.
\hfill\qed

\vspace{.1in}

Let us show the Mehta-Ramanathan type theorem
for regular filtered $\lambda$-flat sheaves.
Let $X$ be a smooth projective variety
with an ample line bundle $L$
and a simple normal crossing divisor $D$.
Let $(\vecE_{\ast},\DDlambda)$ be
a regular filtered $\lambda$-flat sheaf on $(X,D)$.
Let $N$ be a sufficiently large number.
We can take a generic hyper-plane section $Y$ of $L^{\otimes\,N}$
satisfying the properties:
(i) $D_Y:=Y\cap D$ is simply normal crossing in $Y$,
(ii) $\pi_1(Y\setminus D)\lrarr \pi_1(X\setminus D)$
is surjective.

\begin{prop}
 \label{prop;06.1.18.6}
Assume $\dim X\geq 2$.
Then,
$(\vecE_{\ast},\DDlambda)$ is $\mu_L$-stable,
if and only if
$(\vecE_{\ast},\DDlambda)_{|Y}$ is $\mu_L$-stable.
\end{prop}
\pf
Let us fix $\vecc$.
If $W\subset \prolongg{\vecc}{E}$ destabilizes,
the restriction $W_{|Y}$ clearly destabilizes.
Hence, the $\mu_L$-stability of
$(\prolongg{\vecc}{E}_{\ast},\DDlambda)_{|Y}$
implies the $\mu_L$-stability of
$(\prolongg{\vecc}{E}_{\ast},\DDlambda)$.
Assume that 
$(\prolongg{\vecc}{E}_{\ast},\DDlambda)_{|Y}$
is not $\mu_L$-stable,
and let $W$
be a subsheaf of $\prolongg{\vecc}{E}_{|Y}$
satisfying
$\DDlambda(W)\subset 
 W\otimes\Omega^{1,0}_Y(\log D_Y)$
and $\pardeg(W_{\ast})/\rank (W)\geq
 \pardeg(\prolongg{\vecc}{E}_{\ast})/\rank E$.
Let $Q$ be any point of $X-D$.
Take a path $\gamma$ connecting
$Q$ and a point $P$ of $Y\setminus D$.
By the parallel transport along the path,
we obtain the vector subspace
$W'_Q\subset E_{|Q}$.
It is independent of choices of $P$ and $\gamma$,
and we obtain the flat subbundle
$W'\subset \prolongg{\vecc}{E}_{|X-D}$.
Due to Lemma  \ref{lem;06.2.4.20},
we obtain the saturated subsheaf
$\widetilde{W}'\subset \prolongg{\vecc}{E}$.
By a general argument,
it can be shown that
there exists a subset $Z\subset D$
with $\codim_X(Z)\geq 2$
such that $\widetilde{W}'_{\ast|X-Z}$
is a parabolic subbundle of 
$\prolongg{\vecc}{E}_{|X-Z}$.
Then, it is easy to check that
$\widetilde{W}'$ destabilizes.
\hfill\qed

\subsubsection{Saturated regular filtered $\lambda$-flat sheaf}
\label{subsubsection;06.2.5.10}

Let $X$ and $D$ be as above.
Let $(\vecE_{\ast},\DDlambda)$ be
a regular filtered $\lambda$-flat sheaf $(\lambda\neq 0)$.
\begin{df}
 $(\vecE_{\ast},\DDlambda)$ is called
saturated, if the following conditions are satisfied:
\begin{itemize}
\item
 There exists a subset $Z\subset D$ 
 with $\codim_X(Z)\geq 2$,
 and each $\prolongg{\veca}{E}$ are determined 
 on $\prolongg{\veca}{E}_{|X-Z}$.
 Namely, for any open subset $U\subset X$,
 we have the following:
\begin{equation}
 \label{eq;06.2.3.2}
\prolongg{\veca}{E}(U)=
 \prolongg{\veca}{E}(U\setminus Z)
 \cap \vecE(U).
\end{equation}
\end{itemize}
\hfill\qed
\end{df}

It is easy to see that a regular filtered $\lambda$-flat {\em bundle}
is saturated.

\begin{lem}
 \label{lem;06.2.4.2}
Let $(\vecE_{\ast},\DDlambda)$ be a saturated
regular filtered $\lambda$-sheaf on $(X,D)$.
Then, each $\vecc$-truncation
$\prolongg{\vecc}{E}$ is reflexive.
\end{lem}
\pf
Recall we have already known that
$\prolongg{\vecc}{\vecE}_{\ast}$ is a filtered bundle
in codimension one (Lemma \ref{lem;06.2.4.1}).
Let $\prolongg{\vecc}{E}^{\lor\lor}$ denote 
the double dual of $\prolongg{\vecc}{E}$.
We have the naturally defined injective map
$\prolongg{\vecc}{E}\lrarr \prolongg{\vecc}{E}^{\lor\lor}$.
Due to the saturatedness,
any sections of $\prolongg{\vecc}{E}^{\lor\lor}$
naturally gives sections of $\prolongg{\vecc}{E}$,
i.e.,
$\prolongg{\vecc}{E}$ is isomorphic to
$\prolongg{\vecc}{E}^{\lor\lor}$.
\hfill\qed

\begin{lem}
 \label{lem;06.2.4.3}
A saturated regular filtered $\lambda$-flat sheaf
$(\vecE_{\ast},\DDlambda)$ on $(X,D)$
is a regular filtered $\lambda$-flat bundle
in codimension two.
\end{lem}
\pf
We have only to show that
there exists a subset $Z\subset D$
with $\codim_X(Z)\geq 3$
such that $\prolongg{\vecc}{E}_{\ast\,|\,X-Z}$
is a $\vecc$-parabolic bundle on $(X-Z,D-Z)$
for any $\vecc$.
Due to 
$\prolongg{\vecc+\vecb}{E}
=\prolongg{\vecc}{E}\otimes\nbigo(\vecb\cdot D)$,
where $\vecb\cdot D=\sum_{i\in S} b_i\cdot D_i$,
we have only to show such a claim for finite
number of tuples $\vecc$.
Due to Lemma \ref{lem;06.2.4.2},
there exists a subset $Z'\subset D$
with $\codim_X(Z')\geq 3$
such that $\prolongg{\vecc}{E}_{|X-Z'}$ is locally free.
Hence, we can assume 
that $\prolongg{\vecc}{E}$ is locally free
from the beginning.

We have the parabolic filtration
$\lefttop{i}F=\{\lefttop{i}F_a\,|\,c_i-1<a\leq c_i\}$
of $\prolongg{\vecc}{E}_{|D_i}$.
We can take the saturation $\lefttop{i}\widetilde{F}_a$
of $\lefttop{i}{F}_a$.
Namely, we put
$G_a:=\prolongg{\vecc}{E}_{|D_i}\big/\lefttop{i}F_a$,
and let $G_{a\,tor}$ denote the torsion-part
of $G_a$.
Let $\pi_a:\prolongg{\vecc}{E}_{|D_i}\lrarr G_{a}$
denote the projection,
and we put 
$\lefttop{i}\widetilde{F}_a:=\pi_a^{-1}
 \bigl( G_{a\,tor} \bigr)$.
\begin{lem}
 \label{lem;06.2.4.4}
$\lefttop{i}\widetilde{F}_a=
 \lefttop{i}F_a$.
\end{lem}
\pf
By our construction,
we have $\lefttop{i}F_a\subset \lefttop{i}\widetilde{F}_a$,
and we also know that
there exists a subset $W\subset D_i$
with $\codim_{D_i}(W)\geq 1$
such that
$\lefttop{i}F_{a\,|\,D_i-W}
=\lefttop{i}\widetilde{F}_{a\,|\,D_i-W}$.

Let $P$ be any point of $D_i$.
Let $g$ be a germ of a section of $\lefttop{i}\widetilde{F}_a$
at $P$,
and let $G$ be a local section of $\prolongg{\vecc}{E}$
on an open subset $U$of $P$ in $X$
such that the germ of the restriction of $G$ to $D_i$ gives $g$.
Then, $G_{|U\setminus W}$ gives a section of
$\prolongg{\vecc'}{E}$ on $U\setminus W$,
where $\vecc'=(c_j')$ is determined by
$c_j'=c_j$ $(j\neq i)$ and $c_i=a$.
Due to the saturatedness,
$G$ is a section of 
$\prolongg{\vecc'}{E}$ on $U$.
Thus, $g$ is the germ of a section of
$\lefttop{i}F_a$,
and $\lefttop{i}F_a=\lefttop{i}\widetilde{F}_a$.
Hence, we obtain Lemma \ref{lem;06.2.4.4}.
\hfill\qed

\vspace{.1in}
Let us return to the proof of Lemma \ref{lem;06.2.4.3}.
Due to Lemma \ref{lem;06.2.4.4},
the associated graded vector bundle 
$\lefttop{i}\Gr^F(\prolongg{\vecc}{E}_{|D_i})$
is torsion free.
Hence, there exists a subset
$Z''_i\subset D_i$
with $\codim_{D_i}Z''_i\geq 2$ such that
$\lefttop{i}F_{|D_i\setminus Z''_i}$ is a filtration
in the category of vector bundles on $D_i''\setminus Z''_i$.
Then, $\prolongg{\vecc}{E}_{\ast\,|\,X-Z''}$
is a $\vecc$-parabolic locally free sheaf on $(X-Z'',D-Z'')$.
Thus we are done.
\hfill\qed

\begin{rem}
By the correspondence of
saturated regular filtered flat bundles
and filtered local systems,
we can obtain more concrete picture
of the saturated regular filtered flat sheaves.
We will see it in Section {\rm\ref{section;06.2.3.150}}.
\hfill\qed
\end{rem}

\subsubsection{Canonical decomposition}

Let $\bigl(\nbige^{(i)}_{\ast},\DD^{\lambda(i)}\bigr)$ $(i=1,2)$
be $\mu_L$-semistable 
regular $\vecc$-parabolic $\lambda$-flat sheaves
such that 
$\mu_L(\nbige^{(1)}_{\ast})
=\mu_L(\nbige^{(2)}_{\ast})$.
Let $f:(\nbige^{(1)}_{\ast},\DD^{\lambda\,(1)})
 \lrarr (\nbige^{(2)}_{\ast},\DD^{\lambda\,(2)})$ 
be a non-trivial morphism.
Let $(\nbigk_{\ast},\DDlambda_{\nbigk})$ 
denote the kernel of $f$,
which is naturally equipped 
with the parabolic structure and the flat $\lambda$-connection.
Let $\nbigi$ denote the image of $f$,
and $\nbigitilde$ denote 
the saturated subsheaf of $\nbige^{(2)}$
generated by $\nbigi$.
The parabolic structures of $\nbige^{(1)}_{\ast}$ 
and $\nbige^{(2)}_{\ast}$
induce the parabolic structures of
$\nbigi$ and $\nbigitilde$, respectively.
We denote the induced parabolic flat sheaves by
$(\nbigi_{\ast},\DDlambda_{\nbigi})$
and $(\nbigitilde_{\ast},\DDlambda_{\nbigitilde})$.
The following lemma can be shown
by the same argument
as the proof of Lemma 3.9 of \cite{mochi4}.
\begin{lem}
\label{lem;06.8.4.3}
$(\nbigk_{\ast},\DDlambda_{\nbigk})$,
$(\nbigi_{\ast},\DDlambda_{\nbigi})$
and $\bigl(\nbigitilde_{\ast},\DDlambda_{\nbigitilde}\bigr)$
are also $\mu_L$-semistable such that
$\mu_L(\nbigk_{\ast})=\mu_L(\nbigi_{\ast})
=\mu_L(\nbigitilde_{\ast})
=\mu_L(\nbige^{(i)}_{\ast})$.
Moreover,
 $\nbigi_{\ast}$ and $\nbigitilde_{\ast}$
are isomorphic in codimension one.
\hfill\qed
\end{lem}

\begin{lem}
\label{lem;06.8.22.20}
Let $\bigl(\nbige^{(i)}_{\ast},\,
 \DD^{\lambda\,(i)}\bigr)$ $(i=1,2)$
be $\mu_L$-semistable reflexive saturated
regular parabolic $\lambda$-flat sheaves
such that 
$\mu_L(\nbige^{(1)}_{\ast})
=\mu_L(\nbige^{(2)}_{\ast})$.
Assume either one of the following:
\begin{enumerate}
\item
One of $(\nbige^{(i)}_{\ast},\DD^{\lambda\,(i)})$
 is $\mu_L$-stable,
and $\rank(\nbige^{(1)})=\rank(\nbige^{(2)})$ holds.
\item
Both of $(\nbige^{(i)}_{\ast},\DD^{\lambda\,(i)})$
 are $\mu_L$-stable.
\end{enumerate}
If there is a non-trivial map 
$f:(\nbige^{(1)}_{\ast},\DD^{\lambda\,(1)})
 \lrarr
 (\nbige^{(2)}_{\ast},\DD^{\lambda\,(2)})$,
then $f$ is isomorphic.
\end{lem}
\pf
If $(\nbige^{(1)}_{\ast},\DD^{\lambda\,(1)})$ is $\mu_L$-stable,
the kernel of $f$ is trivial due to Lemma \ref{lem;06.8.4.3}.
If $(\nbige^{(2)}_{\ast},\DD^{\lambda\,(2)})$ is $\mu_L$-stable,
the image of $f$ and $\nbige^{(2)}$ are 
the same at the generic point of $X$.
Thus, we obtain that $f$ is generically isomorphic
in any case.
Then, we obtain that $f$ is isomorphic
in codimension one, due to Lemma 3.7 of \cite{mochi4}.
Since both of $\nbige^{(i)}_{\ast}$
are reflexive and saturated,
we obtain that $f$ is isomorphic.
\hfill\qed

\begin{cor}
Let $(\nbige_{\ast},\DDlambda)$ be a $\mu_L$-polystable 
reflexive saturated regular parabolic $\lambda$-flat sheaf.
Then, we have the unique decomposition:
\[
 (\nbige_{\ast},\DDlambda)
=\bigoplus_{j} 
 \bigl(\nbige^{(j)}_{\ast},\DD^{\lambda\,(j)}\bigr)
 \otimes\cnum^{m(j)}.
\]
Here, $(\nbige^{(j)}_{\ast},\DD^{\lambda\,(j)})$
 are $\mu_L$-stable
with $\mu_L(\nbige^{(j)}_{\ast})=\mu(\nbige_{\ast})$,
and they are mutually non-isomorphic.
It is called the canonical decomposition in the rest of the paper.
\hfill\qed
\end{cor}

\subsubsection{Perturbation of parabolic structure}
\label{subsubsection;06.1.19.10}

Let $X$ be a smooth projective {\em surface}
with an ample line bundle $L$,
and $D$ be a simple normal crossing divisor
with the irreducible decomposition $D=\bigcup_{i\in S} D_i$.
Let $(\prolongg{\vecc}{E},\vecF,\DDlambda)$ be
a regular $\vecc$-parabolic 
$\lambda$-flat bundle  over $(X,D)$
for some $\vecc\in\real^S$.
Assume $\lambda\neq 0$.
We also assume $c_i\not\in\Par(\prolongg{\vecc}{E},\vecF,i)$
for each $i\in S$,
for simplicity.
Let $\nbign_i$ denote the nilpotent part
of the induced endomorphism $\Gr^F\Res_i(\DDlambda)$ on
$\lefttop{i}\Gr^F_{a}(\prolongg{\vecc}{E})$.
Before proceeding,
we give a definition of graded semisimplicity,
as in the Higgs case.
\begin{df}
The $\lambda$-flat $\vecc$-parabolic bundle
$(\prolongg{\vecc}{E},\vecF,\DDlambda)$ is called 
graded semisimple,
if the nilpotent parts $\nbign_i$ are $0$ for any $i\in S$.
\hfill\qed
\end{df}

We would like to consider perturbation
of parabolic structure,
as in Subsection 3.4 of \cite{mochi4}.
First, we will recall general construction.
Then, we will give two kinds of perturbations.

\vspace{.1in}

Let $\eta$ be a generic point of $D_i$.
We have the weight filtration
$W_{\eta}$ of the nilpotent map
$\nbign_{i,\eta}$ on
$\lefttop{i}\Gr^F\bigl(\prolongg{\vecc}{E}\bigr)_{\eta}$,
which is indexed by $\seisuu$.
Then, we can extend it to the filtration $W$ of
$\lefttop{i}\Gr^F\bigl(\prolongg{\vecc}{E}\bigr)$
in the category of vector bundles on $D_i$
due to $\dim D_i=1$.
By our construction,
$\nbign_i(W_k)\subset W_{k-2}$.
The endomorphism
$\Res_i(\DDlambda)$ preserves the filtration $W$
on $\lefttop{i}\Gr^{F}(\prolongg{\vecc}{E})$,
and the nilpotent part of the induced endomorphisms
on $\Gr^W\bigl(\lefttop{i}\Gr^F(\prolongg{\vecc}{E})\bigr)$
are trivial.
Recall that the flat $\lambda$-connection 
$\DDlambda$ {\em locally} induces
the $\lambda$-connection $\lefttop{i}\DDlambda$ of
the vector bundle $\prolongg{\vecc}{E}_{|D_i}$ on $D_i$.
Since $\lefttop{i}\Gr^F(\lefttop{i}\DDlambda)$
commutes with  $\Res_i\DDlambda$, it preserves the filtration $W$.

\vspace{.1in}

Let us take the refinement of the filtration $\lefttop{i}F$.
For any $a\in\openclosed{c_i-1}{c_i}$,
we have the surjection
$ \pi_a:\lefttop{i}F_a(\prolongg{\vecc}{E}_{|D_i})
\lrarr
 \lefttop{i}\Gr^F_a(\prolongg{\vecc}{E})$.
We put 
$\lefttop{i}\widetilde{F}_{a,k}:=\pi_a^{-1}(W_k)$.
We use the lexicographic order on
$\openclosed{c_i-1}{c_i}\times\seisuu$.
Thus, we obtain the increasing filtration
$\lefttop{i}\widetilde{F}$ indexed by
$\openclosed{c_i-1}{c_i}\times\seisuu$.
Obviously,
the set 
$\widetilde{S}_i:=
\bigl\{(a,k)\in \openclosed{c_i-1}{c_i}\times\seisuu
 \,\big|\,
 \lefttop{i}\Gr^{\widetilde{F}}_{(a,k)}\neq 0\bigr\}$
is finite.

\vspace{.1in}

We explain the perturbation of the weight
for the parabolic structure.
Let $\varphi_i:\widetilde{S}_i\lrarr \openclosed{c_i-1}{c_i}$
be the increasing map
such that $|\varphi_i(a,k)-a|\leq C\cdot \epsilon$
for some $C>0$.
(Since we are interested in the family 
of the filtrations 
 $\vecF^{(\epsilon)}$ $(\epsilon>0)$,
this condition makes sense.)
Then, $\lefttop{i}\widetilde{F}$ and $\varphi_i$ give
the $\vecc$-parabolic filtration 
$\vecF^{(\epsilon)}=\bigl(\lefttop{i}F^{(\epsilon)}\,\big|\,
 i\in S\bigr)$.
Thus, we obtain 
the regular $\vecc$-parabolic $\lambda$-flat bundle
$\bigl(\prolongg{\vecc}{E},\vecF^{(\epsilon)},\DDlambda\bigr)$,
which are called the $\epsilon$-perturbation
of $(\prolongg{\vecc}{E},\vecF,\DDlambda)$.
By construction,
we have the following convergence in $H^{\ast}(X,\real)$.
\[
 \lim_{\epsilon\to 0}
 \parchern_1(\prolongg{\vecc}{E},\vecF^{(\epsilon)})
=\parchern_1(\prolongg{\vecc}{E},\vecF),
\quad\quad
  \lim_{\epsilon\to 0}
 \parch_2(\prolongg{\vecc}{E},\vecF^{(\epsilon)})
=\parch_2(\prolongg{\vecc}{E},\vecF)
\]

The following proposition is standard.
(See Proposition 3.28 of \cite{mochi4}, for example.)
\begin{prop}
Assume that $\bigl(\prolongg{\vecc}{E},\vecF,\DDlambda\bigr)$
is $\mu_L$-stable.
If $\epsilon$ is sufficiently small,
then the $\epsilon$-perturbation
$\bigl(\prolongg{\vecc}{E},
 \vecF^{(\epsilon)},\DDlambda\bigr)$ is also
$\mu_L$-stable.
\hfill\qed
\end{prop}

We will use two kinds of perturbations 
$\varphi_i$ of parabolic weights.

\begin{description}
\item[(I)]
The image of $\varphi_i$ is contained in $\rnum$
for each $i\in S$.
This perturbation will be used
to obtain the formula
for the parabolic characteristic numbers.
\item[(II)]
 For simplicity, we assume 
$\epsilon=m^{-1}$
 and $0<10\rank E\cdot 
 \epsilon<\gap(\prolongg{\vecc}{E},\vecF)$.
(See Subsection 3.1 of \cite{mochi4}
 for $\gap$.)
Let $i\in S$.
For each $a\in\Par(\prolongg{\vecc}{E},\vecF)$,
we take $a'(\epsilon,i)\in m^{-1}\cdot\seisuu$
such that 
$|a'(\epsilon,i)-a|<m^{-1}$.
Let $L(\epsilon,i)\in\real$ be determined by the following:
\[
 L(\epsilon,i)\cdot\rank(E)
:=\sum (a(\epsilon,i)-a)\cdot
 \rank\lefttop{i}\Gr^{F}_a(\prolongg{\vecc}{E})
\]
Then, we put
$a(\epsilon,i):=a'(\epsilon,i)-L(\epsilon,i)$
and 
$\varphi(a,k):=
 a(\epsilon,i)
+k\cdot \epsilon$.
By construction,
we have the following equality:
\[
 \sum_{a,k} \varphi(a,k)\cdot
 \rank\bigl(\lefttop{i}\Gr^{F}_a(\prolongg{\vecc}{E})\bigr)
=
 \sum_a a\cdot
 \rank\bigl(\lefttop{i}\Gr^{\Ftilde}_{a,k}
 (\prolongg{\vecc}{E})\bigr)
\]
Hence, we have
$\parchern_1(\prolongg{\vecc}{E},\vecF)
=\parchern_1(\prolongg{\vecc}{E},\vecF^{(\epsilon)})$.
For each $i$,
we also have some $-1/m<\gamma_i\leq 0$
such that 
$\Par(\prolongg{\vecc}{E},\vecF^{(\epsilon)},i)$
is contained in
$\bigl\{
 c_i+\gamma_i+p/m\,\big|\,
 p\in\seisuu_{\leq 0},\,
 -1<\gamma_i+p/m\leq 0
 \bigr\}$.
\end{description}

\begin{rem}
The construction given in this subsection
is valid, when the base manifold $X$ is a curve.
\hfill\qed
\end{rem}

\subsubsection{Remarks about the terminology and the notation}
\label{subsubsection;06.1.24.10}

We give some remarks about the terminology ``parabolic structure''.
Let $X$ be a complex manifold,
and let $D$ be a simple normal crossing divisor of $X$
with the irreducible decomposition $D=\bigcup_{i\in S}D_i$.
We often discuss 
a regular $\vecc$-parabolic $\lambda$-flat bundle on $(X,D)$
for some $\vecc\in\real^S$.
In our most arguments,
a choice of $\vecc$ are not relevant.
In fact,
$\vecc$ is fixed to be $(0,\ldots,0)$ in many references
where the parabolic structure is discussed.
But, it is sometimes convenient
to avoid the case $c_i\in\Par(\prolongg{\vecc}{E}_{\ast},i)$,
for example, when we consider a perturbation of 
the parabolic structure.
That is the main reason why we consider
general $\vecc$-parabolic structure.

In the following argument,
we often assume 
$c_i\not\in\Par(\prolongg{\vecc}{E}_{\ast},i)$ implicitly,
and we often omit to distinguish $\vecc$,
and use the terminology ``parabolic structure''
instead of ``$\vecc$-parabolic structure'',
when we do not have to care about a choice of $\vecc$.
The author hopes that there will be no risk of confusion
and that it will reduce unnecessary complexity of the description.

Relatedly we have the remark about the notation
to denote parabolic bundles.
We often use the notation
$(\prolongg{\vecc}{E},\vecF)$ or $\prolongg{\vecc}{E}_{\ast}$
to denote a $\vecc$-parabolic bundle,
when we would like to distinguish $\vecc$.
The notation ``$\prolongg{\vecc}{E}$'' is also appropriate and useful,
when we regard it as a prolongment of the locally free sheaf $E$ on $X-D$.
But, in some case,
a vector bundle is given not only on $X-D$
but also on $X$  from the beginning.
And, as is said above,
we will not care about a choice of $\vecc$.
In such a case,
we often prefer to using the notation
$(E,\vecF)$ or $E_{\ast}$ for simplicity
of the description.

\subsection{Generality for $\lambda$-connection
 in the $C^{\infty}$-category}

We will give some generality for $\lambda$-connections.
They are straightforward generalization
of the argument for Higgs bundles or flat bundles
given in Simpson's papers (for example \cite{s1} and \cite{s5}),
and hence we will often omit to give a detailed proof.
For simplicity, we will assume $\lambda\neq 0$.

\subsubsection{The induced operators}
\label{subsubsection;06.1.10.7}

Let $X$ be a complex manifold,
and $(E,\DDlambda)$ be a flat $\lambda$-connection
on $X$.
We have the decomposition of $\DDlambda$
into the $(0,1)$-part $d''_E$ and the $(1,0)$-part $d'_E$.
The holomorphic structure of $E$ is given by $d''_E$.
Recall that 
the twisted Leibniz rule
$d'_E(f\cdot v)=\lambda\cdot\del_X(f) v+f\cdot d'_Ev$ holds
for $f\in C^{\infty}(X)$ and $v\in C^{\infty}(X,E)$.
Let $h$ be a hermitian metric of $E$.
From $d''_E$ and $h$,
we obtain the $(1,0)$-operator $\delta_{E,h}'$
determined by
$ \delbar h(u,v)
=h(d''_Eu,v)+h(u,\delta'_{E,h}v)$.
From $d'_E$ and $h$,
we obtain the $(0,1)$-operator $\delta''_{E,h}$
determined by
$ \lambda\del h(u,v)
=h(d'_Eu,v)+h(u,\delta''_{E,h}v)$.
We remark 
$\delta''_{E,h}(f\cdot v)=\lambdabar\cdot \delbar_Xf\cdot v
+f\cdot \delta''_{E,h}(v)$.
We obtain the following operators:
\begin{equation}
 \label{eq;06.1.25.2}
\begin{split}
 \delbar_{E,h}:=\frac{1}{1+|\lambda|^2}
 (d_E''+\lambda\delta_{E,h}''),
\quad
 \del_{E,h}:=\frac{1}{1+|\lambda|^2}
 (\lambdabar d'_E+\delta_{E,h}'),
 \\
   \theta_{E,h}^{\dagger}:=\frac{1}{1+|\lambda|^2}
 (\lambdabar d''_E-\delta''_{E,h}),
\quad
 \theta_{E,h}:=\frac{1}{1+|\lambda|^2}
 (d'_E-\lambda\delta'_{E,h}).
\end{split}
\end{equation}
It is easy to see that the following Leibniz rule holds:
\[
 \delbar_{E,h}(fs)=\delbar_Xf\cdot s+f\cdot \delbar_{E,h}s,
\quad
 \del_{E,h}(fs)=\del_Xf\cdot s+f\cdot \del_{E,h}s.
\]
On the other hand,
$\theta$ and $\theta^{\dagger}$ give
the sections of $\End(E)\otimes\Omega^{1,0}$
and $\End(E)\otimes\Omega^{0,1}$
respectively.
We also have the formulas:
\[
 d''_E=\delbar_{E,h}+\lambda\theta_{E,h}^{\dagger},
\quad
 d'_E=\lambda\del_{E,h}+\theta_{E,h},
\quad
 \delta_{E,h}'=\del_{E,h}-\lambdabar\theta_{E,h},
\quad
 \delta_{E,h}''=\lambdabar\delbar_{E,h}-\theta^{\dagger}_{E,h}.
\]
\begin{rem}
The index ``$E,h$'' is attached
to emphasize the bundle $E$ and the metric $h$.
We will often omit them if there are no risk of confusion.
\hfill\qed
\end{rem}

\begin{rem}
We have the hermitian product
$(\cdot,\cdot)_h:\bigl(E\otimes\Omega^{\cdot}\bigr)\otimes
\bigl( E\otimes\Omega^{\cdot}\bigr)
\lrarr \Omega^{\cdot}$ induced by $h$.
For a section $A$ of $\End(E)\otimes\Omega^{p,q}$,
let $A^{\dagger}_h$ denote the section of
$\End(E)\otimes\Omega^{q,p}$
which is the adjoint of $A$ with respect to $h$
in the sense
$ \bigl(A\cdot u,v\bigr)_h
=\bigl(u, A^{\dagger}_hv\bigr)_h$.
The above $\theta^{\dagger}_h$
is the adjoint of $\theta_h$ in this sense.
\hfill\qed
\end{rem}

We put
$ \DDlambdastar_h:=\delta_h'-\delta_h''=
 \del_h+\theta^{\dagger}_h-\lambdabar(\delbar_h+\theta_h)$.
We have the following formula:
\[
 \delbar_h+\theta_h=
 \frac{\DDlambda-\lambda\DDlambdastar_h}{1+|\lambda|^2},
\quad
 \del_h+\theta^{\dagger}_h=
 \frac{\DDlambdastar_h+\lambdabar\DDlambda}{1+|\lambda|^2}.
\]
We recall that $h$ is called a pluri-harmonic metric if 
$(\delbar_h+\theta_h)^2=0$ holds,
i.e.,
$(E,\delbar_h,\theta_h)$ is a Higgs bundle.
The condition is equivalent to
$\bigl[\DDlambda,\DDlambdastar_h\bigr]=0$.
In the following,
a $\lambda$-flat bundle
with pluri-harmonic metric is called
a harmonic bundle.

Let us consider the case 
where $X$ is provided with a Kahler form $\omega$.
For a differential operator
$A$ of $E\otimes\Omega^{\cdot}$ of degree one,
i.e.,
$ A:C^{\infty}(X,E\otimes\Omega^{i})
 \lrarr C^{\infty}(X,E\otimes\Omega^{i+1})$,
let $A^{\ast}$ denote a formal adjoint with respect to $\omega$
and $h$,
i.e.,
$ \int_X (Au,v)_{h,\omega}\dvol_{\omega}=
 \int_X (u,A^{\ast}v)_{h,\omega}\dvol_{\omega}$
hold for any $C^{\infty}$-sections $u$ and $v$
with compact supports.
Here, $(\cdot,\cdot)_{h,\omega}$ denotes
the Hermitiann inner product of appropriate
vector bundles induced by $h$ and $\omega$.

\begin{lem}
$\bigl(\DDlambdastar\bigr)^{\ast}
=\sqrt{-1}\bigl[\Lambda_{\omega},\DDlambda\bigr]$
and
$\bigl(\DDlambda\bigr)^{\ast}
=-\sqrt{-1}\bigl[\Lambda_{\omega},\DDlambdastar\bigr]$.
\end{lem}
\pf
It follows from the relations
$\del^{\ast}=\sqrt{-1}[\Lambda_{\omega},\delbar_E]$,
$\delbar^{\ast}=-\sqrt{-1}[\Lambda_{\omega},\del_E]$,
$\theta^{\ast}=-\sqrt{-1}[\Lambda_{\omega},\theta^{\dagger}]$
and
$(\theta^{\dagger})^{\ast}=\sqrt{-1}[\Lambda_{\omega},\theta]$.
\hfill\qed

\vspace{.1in}
The Laplacian
$\Delta_{h,\omega}^{\lambda}:
C^{\infty}(X,E)\lrarr C^{\infty}(X,E)$
is defined by
$\Delta_{h,\omega}^{\lambda}
:=\sqrt{-1}\Lambda_{\omega} \DDlambda \DDlambdastar$.

\begin{rem}
 \label{rem;06.1.26.20}
For the differential operators of functions,
$\Delta^{\lambda}_{\omega}:=
\sqrt{-1}\Lambda
 (\delbar+\lambda\del)\circ(\del-\lambdabar\delbar)
=(1+|\lambda|^2)\sqrt{-1}\Lambda\delbar\del
=(1+|\lambda|^2)\Delta''_{\omega}$,
where $\Delta''_{\omega}$ denotes the usual
Laplacian $\sqrt{-1}\Lambda_{\omega}\delbar\del$.
\hfill\qed
\end{rem}

\begin{lem}
 \label{lem;06.1.8.1}
When $\lambda\neq 0$,
we have
$\lambdabar^{-1}\del_h^2+\lambda^{-1}\theta_h^2=0$
and
 $\lambda^{-1}\delbar_h^2
+\lambdabar^{-1}(\theta_h^{\dagger})^2=0$.
\end{lem}
\pf
From the flatness $(\DDlambda)^2=0$,
we obtain the following formulas:
\begin{equation}
 \label{eq;06.1.2.10}
 (\delbar_h+\lambda\theta_h^{\dagger})^2
=\delbar_h^2+\lambda\delbar_h\theta_h^{\dagger}
+\lambda^2(\theta^{\dagger}_h)^2=0,
\end{equation}
\begin{equation}
 \label{eq;06.1.2.11}
 (\lambda\del_h+\theta_h)^2
=\lambda^2\del_h^2+\lambda\del_h\theta_h+\theta_h^2=0,
\end{equation}
\begin{equation}
 \label{eq;06.1.2.12}
 \bigl[
 \delbar_h+\lambda\theta_h^{\dagger},\,
 \lambda\del_h+\theta_h
 \bigr]
=\lambda\Bigl(
 \bigl[\delbar_h\,,\del_h\bigr]
+\bigl[\theta_h^{\dagger},\theta_h\bigr]
\Bigr)
+\delbar_h\theta_h+\lambda^2\del_h\theta_h^{\dagger}=0.
\end{equation}
It is easy to see 
$(\delbar_h^2)_h^{\dagger}=-\del_h^2$,
$(\delbar_h\theta_h^{\dagger})^{\dagger}=\del_h\theta_h$
and
$(\theta_h^{\dagger})^2=-(\theta_h^2)^{\dagger}$.
Therefore, we obtain the following equality
from (\ref{eq;06.1.2.10}):
\begin{equation}
 \label{eq;06.1.2.13}
 -\del_h^2+\lambdabar \bigl(\del_h\theta_h\bigr)
-\lambdabar^2\theta_h^2=0.
\end{equation}
From (\ref{eq;06.1.2.11}) and (\ref{eq;06.1.2.13}),
we obtain
$ \bigl( \lambda+\lambdabar^{-1} \bigr)
\del_h^2
+\bigl(\lambda^{-1}+\lambdabar\bigr)\theta_h^2
=(1+|\lambda|^2)
\bigl(\lambdabar^{-1}\del_h^2
 +\lambda^{-1}\theta_h^2\bigr)
=0$,
which gives the first formula in the lemma.
The second formula can be obtained 
by taking the adjoint.
\hfill\qed

\begin{lem}
When $\lambda\neq 0$,
we have
$\lambdabar^{-1}\cdot\del_h\theta^{\dagger}_h
+\lambda^{-1}\cdot\delbar_h\theta_h=0$
and
$\bigl[\del_h\,,\delbar_h\bigr]
+\bigl[\theta_h\,,\theta_h^{\dagger}\bigr]=0$.
\end{lem}
\pf
It is easy to check
$[\del_h,\delbar_h]_h^{\dagger}=-[\del_h,\delbar_h]$,
$[\theta_h,\theta_h^{\dagger}]_h^{\dagger}
=-[\theta_h,\theta_h^{\dagger}]$
and $(\delbar_h\theta_h)_h^{\dagger}
 =\del_h\theta_h^{\dagger}$.
Hence, we obtain the following equality from (\ref{eq;06.1.2.12}):
\begin{equation}
\label{eq;06.1.2.14}
 -[\delbar_h,\del_h]-[\theta_h^{\dagger},\theta_h]
+\lambdabar^{-1}\cdot\del_h\theta_h^{\dagger}
+\lambdabar\cdot\delbar_h\theta_h=0.
\end{equation}
The claim of the lemma immediately follows from
(\ref{eq;06.1.2.12}) and (\ref{eq;06.1.2.14}).
\hfill\qed

\begin{cor}
 \label{cor;06.1.21.250}
When $\lambda\neq 0$,
the pluri-harmonicity of the metric $h$
is equivalent to
the vanishings $\theta_h^2=0$ and 
$\delbar_h\theta_h=0$.
\hfill\qed
\end{cor}

\subsubsection{Local expression}
\label{subsubsection;06.1.14.1}

Let $(E,\DDlambda)$ be a flat $\lambda$-connection,
and let $h$ be a $C^{\infty}$-metric.
Let $\vecv=(v_1,\ldots,v_r)$ be a holomorphic frame of $E$.
Let $H=H(h,\vecv)$ denote
the hermitian matrix valued function of $h$
with respect to $\vecv$,
i.e., $H_{i,j}=h(v_i,v_j)$.
Let us see the local expression of the induced operators.

Let $A$ denote the $M(r)$-valued $(1,0)$-form
of $\DDlambda$ with respect to $\vecv$,
i.e.,
$\DDlambda\vecv=\vecv\cdot A$,
in other words,
$\DDlambda v_i=\sum A_{j\,i}\cdot v_j$.
Let $B$ denote the $(1,0)$-form of $\delta_h'$
with respect to $\vecv$,
i.e.,
$\delta_h'\vecv=\vecv\cdot B$,
and then we have
$ \delbar h(v_i,v_j)
=h\bigl(v_i,\delta'_hv_j\bigr)
=\sum h\bigl(v_i,B_{k,j}v_k\bigr)$.
Hence, $\delbar H=H\cdot \overline{B}$,
i.e.,
we obtain $B=\overline{H}^{-1}\del \Hbar$.
Let $C$ denote the $(0,1)$-form of $\delta_h''$
with respect to $\vecv$,
i.e., $\delta''_h\vecv=\vecv\cdot C$,
and then we have
$ \lambda\cdot\del h(v_i,v_j)=h(d'v_i,v_j)+h(v_i,\delta_h''v_j)
=\sum_kh(A_{k,i}v_k,\,v_j)+\sum_kh(v_i,\,C_{k,j} v_k)$.
Hence, $\lambda\del H=\lefttop{t}AH+H\overline{C}$,
i.e.,
we obtain
$C=\lambdabar \cdot\Hbar^{-1}\delbar \Hbar
      -\Hbar^{-1}\lefttop{t}\Abar \Hbar$.
Thus, we obtain the following:
\[
 \theta_h\vecv=
 \vecv\cdot \frac{1}{1+|\lambda|^2}
 (A-\lambda\Hbar^{-1}\del\Hbar),
\quad
 \delbar_h\vecv=\vecv\cdot\frac{\lambda}{1+|\lambda|^2}
 (\lambdabar\cdot \Hbar^{-1}\delbar\Hbar-A^{\dagger}_h).
\]
Here, $A^{\dagger}$ denote the adjoint of $A$
with respect to $h$,
i.e.,
$A^{\dagger}_h=\Hbar^{-1}\cdot\lefttop{t}\overline{A}\cdot\Hbar$.

\subsubsection{Pseudo curvature and the Hermitian-Einstein condition}

Assume $\lambda\neq 0$.
For a flat $\lambda$-connection $(E,\DDlambda)$
with a hermitian metric $h$,
the pseudo curvature $G(h,\DDlambda)$ is defined as follows:
\[
 G(h,\DDlambda):=
\bigl[
 \DDlambda, \DDlambdastar_h
\bigr]
=-\frac{(1+|\lambda|^2)^2}{\lambda}
(\delbar_h+\theta_h)^2.
\]
Then, a hermitian metric $h$ is a pluri-harmonic metric for 
$(E,\DDlambda)$,
if and only if $G(h,\DDlambda)=0$ holds.
We will often use the notation
$G(h)$ or $G_h$ instead of $G(h,\DDlambda)$
if there are no risk of confusion.

When $X$ is provided with a Kahler form $\omega$,
a Hermitian-Einstein condition for $h$ 
is $\Lambda_{\omega}G(h,\DDlambda)^{\bot}=0$,
where ``$\bot$'' means the trace free part.

\subsubsection{Some relations between curvature and pseudo curvature}
\label{subsection;08.2.6.20}

By the construction of $\delta_h'$,
the operator $d''+\delta_h'$ is a unitary connection of $(E,h)$.
The curvature of $d''+\delta_h'$ is denoted by $R(d'',h)$.
We have the following expression of $R(d'',h)$
 due to $[d'',d']=0$:
\begin{equation}
 \label{eq;06.1.19.1}
 R(d'',h)=\bigl[d'',\delta_{h}'\bigr]
=\bigl[d'',\lambda^{-1}d'\bigr]
-\frac{1+|\lambda|^2}{\lambda}\bigl[d'',\theta_h\bigr]
=-\frac{1+|\lambda|^2}{\lambda}
 \bigl(\delbar_h\theta_h
 +\lambda[\theta_h^{\dagger},\theta_h]\bigr).
\end{equation}

\begin{lem}
 \label{lem;06.1.20.20}
The following equality holds:
\begin{equation}
 \label{eq;06.1.8.20}
 \tr R(d'',h)
=\frac{1}{1+|\lambda|^2}\tr G(\DDlambda,h)
=-\frac{1+|\lambda|^2}{\lambda} \delbar\tr\theta_h.
\end{equation}
\end{lem}
\pf
From (\ref{eq;06.1.19.1}), we obtain
$\tr R(d'',h)=
-(1+|\lambda|^2)\lambda^{-1} \cdot\delbar\tr\theta_h$.
On the other hand,
we have the following:
\[
 \tr G(h,\DDlambda)=
-\frac{\bigl(1+|\lambda|^2\bigr)^2}{\lambda}
 \tr\bigl(
 \delbar_h^2+\delbar_h\theta_h+\theta_h^2
 \bigr)
=
-\frac{\bigl(1+|\lambda|^2\bigr)^2}{\lambda}
 \delbar\tr \theta_h.
\]
Here, we have used  $\tr(\theta_h^2)=0$,
which implies $\tr(\delbar_h^2)=0$
due to Lemma \ref{lem;06.1.8.1}.
Thus we are done.
\hfill\qed

\begin{lem}
\label{lem;06.1.12.20}
In the case $\dim X=2$,
we have the following formula:
\[
 \tr \bigl(R(d'',h)^2\bigr)
=\frac{1}{(1+|\lambda|^2)^2}
 \tr \bigl(G(h,\DDlambda)^2\bigr)
-\frac{(1+|\lambda|^2)^2}{\lambda}
 \delbar\tr(\theta_h^2\cdot \theta_h^{\dagger}).
\]
\end{lem}
\pf
We have the following:
\[
 \tr\bigl(G(h,\DDlambda)^2\bigr)
=\frac{(1+|\lambda|^2)^4}{\lambda^2}
 \Bigl(
 \tr\bigl( (\delbar_h\theta_h)^2\bigr)
+2\tr\bigl(\delbar_h^2\cdot\theta_h^2\bigr)
 \Bigr)
\]
\[
 \tr\bigl(R(h,d'')^2\bigr)
=\frac{(1+|\lambda|^2)^2}{\lambda^2}
 \Bigl(
 \tr\bigl((\delbar_h\theta_h)^2\bigr)
+2\lambda\tr\bigl(\delbar_h\theta_h\cdot
 [\theta_h,\theta_h^{\dagger}]\bigr)
+\lambda^2\tr\bigl([\theta_h,\theta_h^{\dagger}]^2\bigr)
 \Bigr).
\]
Since we have
$ \tr\bigl([\theta_h,\theta_h^{\dagger}]^2\bigr)
=-2\tr\bigl(\theta_h^2\theta_h^{\dagger\,2}\bigr)$
and
$(\delbar_h+\lambda\theta_h^{\dagger})^2
 =\delbar_h^2+\lambda\delbar_h\theta_h^{\dagger}
 +\lambda^2\theta_h^{\dagger\,2}=0$,
we obtain the following:
\[
 \lambda^2\tr\bigl([\theta_h,\theta_h^{\dagger}]^2\bigr)
=-2\tr\bigl(\lambda^2\cdot\theta_h^2\cdot
 \theta_h^{\dagger\,2}\bigr)
=2\tr\Bigl(\delbar_h^2\cdot\theta_h^2
+\lambda\cdot\delbar_h\theta_h^{\dagger}\cdot\theta_h^2\Bigr).
\]
Hence, we have the following equality:
\[
 \tr\bigl(R(h,d'')^2\bigr)
=\left(
 \frac{1+|\lambda|^2}{\lambda} \right)^2
 \Bigl(
 \tr\bigl((\delbar_h\theta_h)^2\bigr)
+2\lambda\tr\bigl(\delbar_h\theta_h\cdot
 [\theta_h,\,\theta_h^{\dagger}]\bigr)
+2\tr\bigl(\delbar_h^2\cdot \theta_h^2\bigr)
+2\lambda\tr\bigl(
 \delbar_h\theta_h^{\dagger}\cdot \theta_h^2\bigr)
 \Bigr).
\]
We also remark the following:
\begin{multline}
 \tr\bigl(\delbar_h\theta_h\cdot
 [\theta_h,\theta_h^{\dagger}]\bigr)
+\tr(\theta_h^2\cdot 
 \delbar_h\theta_h^{\dagger})
=\tr\bigl((\delbar_h\theta_h)\cdot
 \theta_h\cdot \theta_h^{\dagger}\bigr)
+\tr\bigl(\delbar_h\theta_h\cdot
 \theta_h^{\dagger}\cdot \theta_h\bigr)
-\tr\bigl(\theta_h\cdot \delbar_h\theta_h^{\dagger}
 \cdot \theta_h\bigr) \\
=\delbar\tr\bigl(\theta_h\cdot \theta_h^{\dagger}
 \cdot \theta_h\bigr)
=-\delbar\tr(\theta_h^2\cdot \theta_h^{\dagger}).
\end{multline}
Then, the claim of the lemma immediately follows.
\hfill\qed

\subsubsection{Change of hermitian metrics}
\label{subsubsection;06.1.18.3}

Let $h_i$ $(i=1,2)$ be hermitian metrics of $E$.
The endomorphism $s$ is determined by
$h_2=h_1\cdot s$,
i.e., $h_2(u,v)=h_1\bigl(s\!\cdot\! u,v\bigr)
 =h_1\bigl(u,s\!\cdot\! v\bigr)$,
which is self-adjoint with respect to both of $h_i$.
Then, we have the relations
 $\delta'_{h_2}=\delta_{h_1}'+s^{-1}\delta_{h_1}'s$
and
 $\delta_{h_2}''=\delta_{h_1}''+s^{-1}\delta''_{h_1}s$.
Therefore, we have the following relations
from (\ref{eq;06.1.25.2}):
\[
 \delbar_{h_2}
=\delbar_{h_1}
+\frac{\lambda}{1+|\lambda|^2}s^{-1}\delta''_{h_1}s,
\quad
 \del_{h_2}=\del_{h_1}
+\frac{1}{1+|\lambda|^2}s^{-1}\delta'_{h_1}s,
\]
\[
 \theta^{\dagger}_{h_2}=\theta^{\dagger}_{h_1}
-\frac{1}{1+|\lambda|^2}s^{-1}\delta''_{h_2}s,
\quad
 \theta_{h_2}=\theta_{h_1}
-\frac{\lambda}{1+|\lambda|^2}s^{-1}\delta'_{h_1}s.
\]
We also have 
$\DDlambdastar_{h_2}=\DDlambdastar_{h_1}
+s^{-1}\DDlambdastar_{h_1}s$,
and thus
$\bigl[\DDlambda,\DDlambdastar_{h_2}\bigr]=
  \bigl[\DDlambda,\DDlambdastar_{h_1}\bigr]
+\DDlambda(s^{-1})\cdot \DDlambdastar_{h_1}s
+s^{-1}\DDlambda\DDlambdastar_{h_1}s$.
Then, we obtain the following formula:
\begin{equation}
 \label{eq;06.1.8.10}
 \Delta^{\lambda}_{h_1,\omega}s
=s\sqrt{-1}\bigl(
 \Lambda_{\omega}G(h_2)
-\Lambda_{\omega}G(h_1)
 \bigr)
+\sqrt{-1}\Lambda_{\omega}\DDlambda s\cdot s^{-1}
 \DDlambdastar s.
\end{equation}
In particular,
we obtain the following formula by taking the trace:
\begin{equation}
 \label{eq;06.1.8.11}
 \Delta^{\lambda}_{\omega}\tr(s)
=
\tr\Bigl(
s\sqrt{-1} \bigl(
 \Lambda_{\omega}G(h_2)
-\Lambda_{\omega}G(h_1)
 \bigr) \Bigr)
-\bigl|\DDlambda(s)s^{-1/2}
 \bigr|^2_{h_1,\omega}.
\end{equation}
As in Lemma 3.1 of \cite{s1}, 
we can derive the following inequality:
\begin{equation}
\label{eq;06.1.8.12}
 \Delta^{\lambda}_{\omega}\log\tr(s)
\leq
  \bigl|\Lambda_{\omega}G(h_1) \bigr|_{h_1}
+\bigl|\Lambda_{\omega}G(h_2) \bigr|_{h_2}
\end{equation}

\subsection{Review of existence result
 of a Hermitian-Einstein metric due to Simpson}

\subsubsection{Analytic stability of flat $\lambda$-bundle}
\label{subsubsection;06.1.25.5}

Let $X$ be a complex manifold
with a Kahler form $\omega$.
In this subsection, we impose the following condition
as in \cite{s1}.
\begin{condition}
 \label{condition;06.1.10.1}
 \mbox{{}}
\begin{enumerate}
\item
 The volume of $X$ with respect to $\omega$ 
 is finite.
\item
 There exists a $C^{\infty}$-function $\phi:X\lrarr\real_{\geq\,0}$
 with the following properties:
\begin{itemize}
\item $\{x\in X\,|\,\phi(x)\leq a\}$ is compact
 for any $a$.
\item
 $0\leq \sqrt{-1}\del\delbar\phi\leq C\cdot\omega$,
 and $\delbar\phi$ is bounded with respect to $\omega$.
\end{itemize}
\item
 There exists a continuous increasing function 
 $a:\closedopen{0}{\infty}\lrarr \closedopen{0}{\infty}$
 with the following properties:
\begin{itemize}
\item $a(0)=0$ and $a(t)=t$ for $t\geq 1$.
\item
 Let $f$ be a positive bounded function on $X$
 such that $\Delta_{\omega} f\leq B$ for some $B\in\real$.
 Then, there exists a constant $C(B)$, depending only on $B$,
 such that
$ \sup_X|f|\leq C(B)
 \cdot a\left(\int_X|f|\cdot \dvol_{\omega}\right)$. 
Moreover,
 $\Delta_{\omega}(f)\leq 0$ implies
 $\Delta_{\omega}(f)=0$.
\hfill\qed
\end{itemize}
\end{enumerate}
\end{condition}

Let $(E,\DDlambda)$ be a $\lambda$-flat bundle
on $X$.
There are two conditions on the finiteness of
the pseudo curvature
of $(E,\DDlambda,h)$.
The stronger one is the following:
\begin{equation}
 \label{eq;06.1.8.25}
 \sup \bigl|G(h,\DDlambda)\bigr|_{h,\omega}<\infty.
\end{equation}
The finiteness (\ref{eq;06.1.8.25}) implies the weaker one:
\begin{equation}
 \label{eq;06.1.8.26}
 \sup \bigl|\Lambda_{\omega}
 G(h,\DDlambda)\bigr|_{h,\omega}<\infty.
\end{equation}

When we are given a hermitian metric $h$ of $E$ 
satisfying the finiteness (\ref{eq;06.1.8.26}),
the degree $\deg_{\omega}(E,h)$ is defined as follows:
\[
 \deg_{\omega}(E,h):=
 \frac{\sqrt{-1}}{2\pi}
 \int_X\frac{\tr G(h,\DDlambda)}{1+|\lambda|^2}
 \cdot\omega^{n-1}
=\frac{\sqrt{-1}}{2\pi}
 \int_X\tr R(h,d'')\cdot\omega^{n-1}.
\]
Here, we have used (\ref{eq;06.1.8.20}).
For any $\lambda$-flat bundle
$(V,\DDlambda_V)\subset (E,\DDlambda)$,
the restriction $h_V:=h_{|V}$ induces $\deg_{\omega}(V,h_V)$.
As in Lemma 3.2 of \cite{s1},
we have the Chern-Weil formula.
The proof is same.
\begin{lem}
\label{lem;06.1.10.2}
Let $\pi_V$ denote the orthogonal projection of $E$
onto $V$.
Then, the following equality holds:
\[
 \deg_{\omega}(V,h_V)=
 \frac{1}{2\pi}
 \frac{1}{1+|\lambda|^2}
\left(
 \sqrt{-1}\int_X\tr\bigl(\pi_V \circ G(h,\DDlambda)\bigr)
 \cdot\omega^{n-1}
-\int_X\bigl|\DDlambda\pi_V\bigr|^2_{h,\omega}
\right)
\]
The value is finite or $-\infty$,
when {\rm(\ref{eq;06.1.8.26})} is satisfied.
\hfill\qed
\end{lem}

\begin{df}
$(E,\DDlambda,h)$ is defined to be
analytically stable with respect to $\omega$,
if the inequality 
\[
 \frac{\deg_{\omega}(V,h_V)} {\rank V}
<\frac{\deg_{\omega}(E,h)}{\rank E}
\]
holds for any $(V,\DDlambda_V)\subset (E,\DDlambda)$.
\hfill\qed
\end{df}

\subsubsection{Existence theorem of Simpson 
 and some consequence}

\begin{prop}
  \label{prop;06.1.13.250}
 Let $(X,\omega)$ be a Kahler manifold
 satisfying Condition {\rm\ref{condition;06.1.10.1}},
 and let $(E,\DDlambda,h_0)$ be 
 a metrized flat $\lambda$-connection satisfying
 {\rm(\ref{eq;06.1.8.25})}.
 Assume that $(E,\DDlambda,h_0)$ is analytically stable
 with respect to $\omega$.
 Then, there exists a hermitian metric $h=h_0\cdot s$
 satisfying the following conditions:
\begin{itemize}
\item
 $h$ and $h_0$ are mutually bounded.
\item
 $\det(h)=\det(h_0)$.
\item
 $\DDlambda(s)$ is $L^2$ with respect to $h_0$ and $\omega$.
\item
 It satisfies the Hermitian Einstein condition
 $\Lambda_{\omega} G(h)^{\bot}=0$,
 where $G(h)^{\bot}$ denotes the trace free part
 of $G(h)$.
\item
 The following equalities hold:
\[
 \int_{Y}\tr\Bigl(G(h)^2\Bigr)\cdot\omega^{n-2}
=\int_{Y}\tr\Bigl(G(h_0)^2\Bigr)\cdot\omega^{n-2},
\quad\quad
 \int_Y\tr\Bigl(G(h)^{\bot\,2}\Bigr)\cdot\omega^{n-2}
=\int_Y\tr\Bigl(G(h_0)^{\bot\,2}\Bigr)\cdot\omega^{n-2}.
\]
\end{itemize}
\end{prop}
We do not give a proof of this proposition,
because we need only minor modification of the proof
of Theorem 1, Proposition 3.5 and Lemma 7.4 of \cite{s1}.
Indeed, we have only to replace
$D''$, $D'$ and $F(h)$
with $\DDlambda$, $\DDlambdastar$ and $G(h)$,
and to make some obvious modification
of positive constant multiplications,
as was mentioned by Simpson himself.
(See the page 754 of \cite{s2}, for example.
 Remark that ``$D^c$'' corresponds to our $-\DDlambdastar$,
 and hence our $G(h)$ is slightly different from his.)
The author recommends the reader to read
a quite excellent discussion in \cite{s1}.
However,
we will use some results related with the Donaldson functional,
which are obtained from the proof.
Hence, we recall a brief outline of the proof of Proposition 
\ref{prop;06.1.13.250}.
We will use the notation in Subsection 
\ref{subsection;06.1.15.35}.

Let $h_{0}$ be a metric for $(E,\DDlambda)$
satisfying the finiteness (\ref{eq;06.1.8.26}).
Let us consider the heat equation 
for the self adjoint endomorphisms $s_t$
with respect to $h_{0}$:
\begin{equation}
 \label{eq;06.1.15.30}
 s_t^{-1}\frac{ds_t}{dt}=-\sqrt{-1}\Lambda_{\omega}G(h_t)^{\bot}.
\end{equation}
A detailed argument to solve (\ref{eq;06.1.15.30})
is given in Section 6 of \cite{s1}.
Moreover, $\Lambda_{\omega} G(h_t)$ is 
shown to be uniformly bounded.
We do not reproduce them here.

Then, we would like to show the existence
of an appropriate subsequence $t_i\to\infty$
such that $\{s_{t_i}\}$ converges to $s_{\infty}$
weakly in $L_2^p$ locally on $X$,
and we would like to show that $h_{\infty}=h_0\cdot s_{\infty}$ gives
the desired Hermitian-Einstein metric.
For that purpose, Simpson used the Donaldson functional
$M\bigl(h_0,h_0s_{t_i}\bigr)$.
(We recall the definition and some fundamental property
in Subsection \ref{subsection;06.1.15.35}, below.)
He showed that there exist positive constants $C_i$ $(i=1,2)$
such that the following holds:
(Proposition 5.3 of \cite{s1}.
We review it in Proposition \ref{prop;06.1.4.3}.
 We will use the notation there in the following.)
\begin{equation}
 \label{eq;06.1.15.40}
 \sup |\log s_t|\leq C_1+C_2\cdot M(h_0,h_0s_t).
\end{equation}
He also showed
(Lemma 7.1 of \cite{s1})
that $M(h_0,h_0s_t)$ is $C^1$ with respect to $t$,
and that the following formula holds:
\begin{equation}
\label{eq;06.1.15.41}
 \frac{d}{dt}M\bigl(h_0,h_0s_t\bigr)=
-\int_X\bigl|\Lambda_{\omega}G(h_t)^{\bot}\bigr|^2_{h_t,\omega}
\leq 0.
\end{equation}
Since we have $M(h_0,h_0)=0$ by definition,
we obtain $M(h_0,h_0s_t)\leq 0$ from (\ref{eq;06.1.15.41}).
Then, we obtain the boundedness of $s_t$
from (\ref{eq;06.1.15.40}).
For the solution of (\ref{eq;06.1.15.30}),
we have $\det(s_t)=1$.
Hence, we also obtain the boundedness of $s_t^{-1}$.
We also obtain the existence of a subsequence $\{t_i'\}$
such that $\bigl|\Lambda_{\omega}G(h'_{t_i})\bigr|_{L^2}\lrarr 0$.

From the uniform boundedness of $s_t$
and $\Lambda_{\omega}G(h_t)$,
we obtain the lower bound of $M\bigl(h_0,h_0s_t\bigr)$.
(See Corollary \ref{cor;06.1.4.1} in this paper, for example.)
Moreover,
we obtain the uniform bound  of
$\int_X \bigl|\DDlambda u_t\bigr|_{h_0}^2$
due to the positivity of $\Psi$ given in (\ref{eq;06.1.15.50}),
where $s_t=\exp(u_t)$.
Due to the boundedness of $s_t$ and $s_t^{-1}$,
we also obtain the boundedness of
$\int_X\bigl|\DDlambda s_t\bigr|^2_{h_0}$.
Then, we obtain the $L_1^2$ boundedness.
Hence, we can take a subsequence $\{t_i''\}$
such that $s_{t_i''}$ converges to some $s_{\infty}$
weakly in $L_1^2$ locally on $X-D$.
Due to some more excellent additional argument given
in the page 895 of \cite{s1},
it can be shown that the convergence is
weakly $L_2^p$ locally on $X-D$, for any $p$.
As a result, we obtain the Hermitian-Einstein metric.

By the above argument,
we can derive the following lemma,
which we would like to use in the later discussion.
\begin{lem}
 \label{lem;06.1.15.70}
Let $h_0$ be 
the hermitian metric satisfying {\rm(\ref{eq;06.1.8.25})}.
Let $h_{HE}$ be the Hermitian-Einstein metric 
obtained in Proposition {\rm\ref{prop;06.1.13.250}}.
Then, we have
$M\bigl(h_0,h_{HE}\bigr)\leq 0$.
\end{lem}
\pf
Recall that $h_{HE}$ is obtained as the limit
$h_0\cdot s_{\infty}$ of some sequence $\{h_0s_{t_i}\}$,
and we have $M(h_0,h_0\cdot s_{t_i})\leq 0$.
We use the formula (\ref{eq;06.1.3.1}).
Let $Z$ be any compact subset of $X$.
The sequence $\{s_{t_i}\}$ converges to $s_{\infty}$ 
in $C^0$ on $Z$.
The sequence $\{\Lambda_{\omega}G(h_{t_i})\}$
converges to $\Lambda_{\omega}G(h_{HE})$ 
weakly in $L^2$ on $Z$.
Therefore, we have the convergence:
\[
\lim_{t_i\to\infty} \int_Z 
 \tr \bigl(u_{t_i}\cdot \Lambda_{\omega}G(h_{t_i})\bigr)
 \dvol_{\omega}
=
\int_Z\tr\bigl(u_{\infty}\cdot \Lambda_{\omega}G(h_{HE})\bigr)
 \dvol_{\omega}.
\]
Here, $u_{t}$ are given by $\exp(u_{t})=s_t$.
Since $\sup_X|s_{t}|$ and $\sup_X|\Lambda G(h_t)|$ are 
bounded independently of $t$,
we can easily obtain the convergence:
\[
\underset{t_i\to\infty}{\lim}
 \int_X \tr \bigl(u_{t_i}\cdot \Lambda_{\omega}G(h_{t_i})\bigr)
  \dvol_{\omega}
=\int_X\tr\bigl(u_{\infty}\cdot \Lambda_{\omega}G(h_{HE})\bigr)
  \dvol_{\omega}.
\]
We have the $C^0$-convergence of the sequence
$\{\DDlambda u_{t_i}\}$ to $\DDlambda u_{\infty}$.
Hence, we have the following inequality due to Fatou's lemma:
\[
 \int_X\bigl(
 \Psi(u_{\infty})\DDlambda u_{\infty},\,\,
 \DDlambda u_{\infty}
 \bigr)\dvol_{\omega}
\leq\underline{\lim}
 \int_X \bigl(
\Psi(u_{t_i})\DDlambda u_{t_i},\,\,
 \DDlambda u_{t_i}
 \bigr)\dvol_{\omega}.
\]
Then, we obtain the desired inequality.
\hfill\qed

\subsubsection{Uniqueness}

The following proposition can be shown
by an argument similar to the proof of
Proposition 2.6 of \cite{mochi4}
via the method in \cite{s1}.
We state it for the reference in the later discussion.
\begin{prop}
Let $(X,\omega)$ be a complete Kahler manifold
satisfying Condition {\rm\ref{condition;06.1.10.1}},
and $(E,\DDlambda)$ be a $\lambda$-flat bundle
on $X$.
Let $h_i$ $(i=1,2)$ be hermitian metrics of $E$
such that
$\Lambda_{\omega}G(h_i)=0$.
We assume that $h_i$ $(i=1,2)$ are mutually bounded.
Then, the following holds:
\begin{itemize}
\item
 We have the decomposition of
 $\lambda$-flat bundles
 $(E,\DDlambda)=\bigoplus (E_a,\DDlambda_a)$
 which is orthogonal with respect to
 both of $h_i$ $(i=1,2)$.
\item
 The restrictions of $h_i$ to $E_a$ are denoted by
 $h_{i,a}$.
 Then, there exist positive numbers $b_a$
 such that $h_{1,a}=b_a\cdot h_{2,a}$.
\end{itemize}
\end{prop}
\pf
Let $s$ be determined by $h_2=h_1\cdot s$.
We can show $\DDlambda s=0$
by the argument explained in the proof of
Proposition 2.6 of \cite{mochi4}.
Note we are considering the case $\lambda\neq 0$.
Hence, the eigen decomposition of $s$
is $\DDlambda$-flat,
which gives the desired decomposition.
\hfill\qed

\subsection{Review of Donaldson functional}

\label{subsection;06.1.15.35}

We recall the Donaldson functional,
by following Donaldson and Simpson
(\cite{don3} and \cite{s1}).

\subsubsection{Functions of self-adjoint endomorphisms}

\label{subsubsection;06.1.9.2}

Let $V$ be a vector space over $\cnum$
with a hermitian metric $h$.
Let $S(V,h)$ denote the set of the endomorphisms of $V$
which are self-adjoint with respect to $h$.
Let $\varphi:\real\lrarr\real$ be a continuous function.
Then, $\varphi(s)$ is naturally defined for any $s\in S(V,h)$.
Namely,
let $v_1,\ldots,v_r$ be the orthogonal base 
which consists of the eigen vectors of $s$,
and let $v_1^{\lor},\ldots,v_r^{\lor}$ be 
the dual base.
Then, we have the description
$s=\sum \kappa_{i}\cdot v_i^{\lor}\otimes v_i$,
and we put
$\varphi(s):=\sum \varphi(\kappa_i)\cdot v_i^{\lor}\otimes v_i$.
Thus, we obtain the induced map $\varphi:S(V,h)\lrarr S(V,h)$,
which is well known to be continuous.
To see the continuity, for example,
we can argue as follows:
Let $U(h)$ denote the unitary group with respect to $h$.
Take $\vece=(e_1,\ldots,e_r)$  be an orthogonal base of $V$.
Let $T$ denote the set of endomorphisms of $V$
which is diagonal with respect to the base $\vece$.
Then, we have the continuous surjective map
$\pi:U(h)\times T\lrarr S(V,h)$
given by $(u,t)\longmapsto u\cdot t\cdot u^{-1}$.
It is easy to check the continuity of the composite
$\varphi\circ\pi$.
Since the topology of $S(V,h)$ is same
as the induced topology via $\pi$,
we obtain the continuity.
When $\varphi$ is real analytic
given by the convergent power series $\sum a_j\cdot t^j$,
then $\varphi(s)=\sum a_j\cdot s^j$.
The induced map is real analytic in this case.

Let $\Psi:\real\times\real\lrarr\real$
be a continuous function.
For a self-adjoint map $s\in S(V,h)$, let $v_1,\ldots,v_r$
and $v_1^{\lor},\ldots,v_{r}^{\lor}$ be as above.
Then, 
we put
$\Psi(s)(A)=\sum
 \Psi(\kappa_i,\kappa_j)\cdot A_{i,j}\cdot v_i^{\lor}\otimes v_j$
for any endomorphism
$A=\sum A_{i,j}\cdot v_i^{\lor}\otimes v_j $ of $V$.
Thus, we obtain $\Psi:S(V,h)\lrarr S(\End(V),h)$,
which is also well known to be continuous.
Here, $S(\End(V),h)$ denotes the set of the self-adjoint
endomorphisms of $\End(V)$ with respect to the metric
induced by $h$.
To see the continuity, we can use the same argument
as above.
When $\Psi$ is real analytic given by a power series,
$\sum b_{m,n}t_1^mt_2^n$,
then we have $\Psi(s)(A)=\sum b_{m,n}s^m \cdot A \cdot s^n$,
and the induced map is real analytic.

Let $\varphi:\real\lrarr\real$ be $C^1$,
and let $d\varphi:\real^2\lrarr\real^2$
denote the continuous function given by
$ d\varphi(t_1,t_2)
=(t_1-t_2)^{-1}\bigl(\varphi(t_1)-\varphi(t_2)\bigr)$
$(t_1\neq t_2)$
and $d\varphi(t_1,t_1)=\varphi'(t_1)$.
In this case,
the induced map
$\varphi:S(V,h)\lrarr S(V,h)$ is also $C^1$,
and the derivative at $s$ is given by $d\varphi(s)$.
To see it, we can argue as follows:
When $\varphi$ is real analytic,
the claim can be checked by a direct calculation.
In general, we can take an approximate sequence
$\varphi_i\lrarr \varphi$ by real analytic functions
on an appropriate compact neighbourhoods
of the eigenvalues of $s\in S(V,h)$.
The induced maps
$\varphi_i:S(V,h)\lrarr S(V,h)$
and $d\varphi_i:S(V,h)\lrarr S(\End(V),h)$
uniformly converge  to $\varphi$ and $d\varphi$
on an appropriate compact neighbourhoods of $s$.
Then, we can derive that $\varphi$
is the integral of the form $d\varphi$ by a general fact.

The construction can be done on manifolds.
Namely, let $E$ be a $C^{\infty}$-vector bundle
with a hermitian metric $h$.
Let $S_h(E)$ (or simply $S_h$)
be the bundle of the self-adjoint endomorphisms
of $(E,h)$,
and let $S_h(\End(E))$ be the bundle of 
the self-adjoint endomorphisms of $(\End(E),h)$.
Then, a continuous function $\varphi:\real\lrarr\real$
induces $\varphi:S_h(E)\lrarr S_h(E)$,
and $\Psi:\real^2\lrarr\real$ induces
$\Psi:S_h(E)\lrarr S_h(\End(E))$.
We have
$ \DDlambda \varphi(s)
=d\varphi(s)\bigl(\DDlambda s\bigr)$,
when $\varphi$ is $C^1$.

\subsubsection{A closed one form}

Let $(X,\omega)$ and $(E,\DDlambda)$ be as in Subsection
\ref{subsubsection;06.1.25.5}.
Following Simpson \cite{s1},
we introduce the space $P(S_h)$,
which consists of sections $s$ of $S_h(E)$
satisfying the following finiteness:
\[
 \|s\|_{h,\omega,P}:=
 \sup_X|s|_{h}+\|\DDlambda s\|_{2,h,\omega}
+\|\Delta^{\lambda}_{h,\omega}s\|_{1,h,\omega}
<\infty.
\]
Here, $\|\cdot\|_{p,h,\omega}$ denote the $L^p$-norm
with respect to $(h,\omega)$.
We will omit to denote $\omega$ and $h$,
when there are no risk of confusion.
The following lemma corresponds to
Proposition 4.1 (d) in \cite{s1}.
The proof is same.
\begin{lem}
\label{lem;06.1.9.1}
Let $\varphi$ and $\Psi$ are analytic functions on $\real$
with infinite radius of convergence.
Then, $\varphi:P(S_h)\lrarr P(S_h)$
and $\Psi:P(S_h)\lrarr P(S_h(\End(E)))$ are analytic.
\hfill\qed
\end{lem}

Let $h$ be a metric satisfying (\ref{eq;06.1.8.26}).
Let $\nbigp_+(S_h)$ denote the set of
the self-adjoint positive definite endomorphisms $s$
with respect to $h$ 
such that $\|s\|_{h,P}<\infty$ and $\|s^{-1}\|_{h,P}<\infty$.
Note $\|s\|_{h,P}<\infty$
and $\sup|s^{-1}|_h<\infty$
imply $\|s^{-1}\|_{h,P}<\infty$.
We put
$\nbigp_h:=\bigl\{h\cdot s\,\big|\,s\in \nbigp_+(S_h)\bigr\}$.
It is easy to see that
any $h_1\in\nbigp_h$ also satisfies (\ref{eq;06.1.8.26})
due to (\ref{eq;06.1.8.10}).
It is also easy to see
$\nbigp_h=\nbigp_{h_1}$ for $h_1\in\nbigp_h$.

Let $\nbigp(S_h)$ denote the space
of the self-adjoint endomorphisms $s$ with respect to $h$
such that $\|s\|_{P,h}<\infty$.
It is easy to see that
$\nbigp_+(S_h)$ is open in $\nbigp(S_h)$.
In particular, we obtain the Banach manifold structure
of $\nbigp_+(S_h)$.
By the natural bijection $\nbigp_h\simeq \nbigp_+(S_{h_1})$
 for $h_1\in\nbigp_h$,
we also obtain the Banach manifold structure of $\nbigp_h$,
which is independent of a choice of $h_1\in\nbigp_h$.
We have the map $\nbigp(S_{h_1})\lrarr \nbigp_+(S_{h_1})$
given by $s\longmapsto e^s$
(Lemma \ref{lem;06.1.9.1}).
It gives a diffeomorphism around 
$0\in \nbigp(S_{h_1})$ and $1\in\nbigp_+(S_{h_1})$.
Therefore,
the map $\nbigp(S_{h_1})\lrarr \nbigp_h$ 
by $s\longmapsto h_1\cdot e^s$
gives a diffeomorphism around $0$ and $h_1$.
In particular, 
the tangent space $T_{h_1}\nbigp_h$ can be naturally identified
with $\nbigp(S_{h_1})$ for any $h_1\in \nbigp_h$.
We also have the natural isomorphism
$\nbigp(S_{h_1})\simeq\nbigp(S_{h})$
given by $t\longmapsto u\cdot t$
for  $h_1=h\cdot u\in \nbigp_h$,
which gives the local trivialization of the tangent bundle.

For any $h_1\in \nbigp_h$ and $s\in T_{h_1}\nbigp_h$,
we put as follows:
\[
 \Phi_{h_1}(s):=
 \int_X\Phi'_{h_1}(s)\dvol_{\omega}\in\cnum,
\quad
 \Phi'_{h_1}(s):=
 \sqrt{-1}\tr\bigl(s\cdot \Lambda_{\omega} G(\DDlambda,h_1)\bigr).
\]
Then, $\Phi'$ gives the $L^1(X,\Omega_X^{1,1})$-valued
one form on $\nbigp_h$,
and $\Phi$ gives the one form of $\nbigp_h$.
The differentiability of $\Phi$ is easy to see.

\begin{lem}
$\Phi$ is a closed one form.
\end{lem}
\pf
In the following argument,
we use the notation $\DDlambdastar$ instead of 
$\DDlambdastar_h$.
Let $k_1,k_2\in \nbigp_h$.
They naturally give the vector field by addition.
At any point $h_1\in\nbigp_{h}$, they give
the tangent vectors $\sigma=h_1^{-1}k_1$ and
$\tau=h_1^{-1}k_2$ in $T_{h_1}\nbigp_h=\nbigp(S_{h_1})$.
Hence, we have the following at $h+\epsilon k_1$:
\[
 \Phi_{h+\epsilon k_1}(k_2)=
 \sqrt{-1}\int \tr\Bigl(
 (h+\epsilon k_1)^{-1}\cdot
 k_2\cdot G(h+\epsilon k_1)
 \Bigr)\cdot\omega^{n-1}.
\]
We have
$(h+\epsilon k_1)^{-1}k_2
=(1+\epsilon \sigma)^{-1}\tau
=\tau-\epsilon\sigma\tau
+(1+\epsilon\sigma)^{-2}\epsilon^2\sigma^2\tau$.
Remark $\sigma^2\tau$ is bounded.
We also have the following:
\begin{multline}
(1+\epsilon\sigma)
\bigl(
 G(h+\epsilon k_1)-G(h)
\bigr)
=\DDlambda\DDlambdastar(1+\epsilon\sigma)
-\DDlambda(1+\epsilon\sigma)\cdot
(1+\epsilon\sigma)^{-1}
 \DDlambdastar(1+\epsilon\sigma) \\
=\epsilon\DDlambda\DDlambdastar\sigma
-\epsilon^2\DDlambda\sigma\cdot
(1+\epsilon\sigma)^{-1}\DDlambdastar\sigma.
\end{multline}
Hence, we have
$G(h+\epsilon k_1)-G(h)=
\epsilon \DDlambda\DDlambdastar \sigma+
\epsilon^2 R_0(\epsilon,\sigma,\tau)$,
where $R_0(\epsilon,\sigma,\tau)$ is an $L^1$-section
of $\End(E)\otimes\Omega^2$, and
the $L^1$-norm is bounded independently from $\epsilon$.
Therefore, we obtain the following:
\begin{multline}
 \Phi_{h+\epsilon k_1}(k_2)-\Phi_{h}(k_2)
=\sqrt{-1}
\int \tr\bigl(
 (h+\epsilon k_1)^{-1}\cdot k_2\cdot
 G(h+\epsilon k_1)\bigr)
 \cdot\omega^{n-1}
-\sqrt{-1}\int \tr\bigl(
 h^{-1}\cdot k_2 \cdot G(h)\bigr)\cdot\omega^{n-1}
\\
=\sqrt{-1}
 \int\tr\bigl(
 \tau G(h+\epsilon k_1)-\tau G(h)
 \bigr)\cdot\omega^{n-1}
-\epsilon\sqrt{-1}\int
 \tr\bigl(\sigma\tau G(h+\epsilon k_1)\bigr)\cdot\omega^{n-1}
+\epsilon\cdot R_1(\epsilon,\sigma,\tau) \\
=
 \epsilon\left(
 \sqrt{-1}
 \int \tr\bigl(\tau \DDlambda\DDlambdastar\sigma\bigr)
 \cdot\omega^{n-1}
-\sqrt{-1}\int \tr\bigl(\sigma\cdot \tau\cdot G(h)\bigr)
 \cdot\omega^{n-1}
 \right)
+\epsilon R_2(\epsilon,\sigma,\tau).
\end{multline}
Here, we have $R_i(\epsilon,\sigma,\tau)\lrarr 0$ $(i=1,2)$
in $\epsilon\to 0$,
due to $\|\sigma\|_P<\infty$ and $\|\tau\|_P<\infty$.
Hence, we obtain the following equality:
\[
  d_h\Phi(\sigma,\tau)
=\sqrt{-1}\int \Bigl(
 \tr\bigl(\tau\DDlambda\DDlambdastar \sigma\bigr) 
-\tr\bigl(\sigma\DDlambda\DDlambdastar\tau\bigr)
 \Bigr)\cdot\omega^{n-1}
-\sqrt{-1}\int\tr\bigl([\sigma,\tau]\cdot G(h)\bigr)
 \cdot\omega^{n-1}.
\]
We have the following equality,
due to $[\DDlambda,\DDlambdastar]=G(h)$:
\begin{multline}
(-\lambdabar\delbar+\del)\tr(\tau\DDlambda\sigma)
+(\lambda\del+\delbar)\tr(\sigma\DDlambdastar\tau)  
=\tr(\DDlambdastar\tau\DDlambda\sigma)
+\tr(\tau\DDlambdastar\DDlambda\sigma)
+\tr(\DDlambda\sigma\DDlambdastar\tau)
+\tr(\sigma\DDlambda\DDlambdastar\tau)  \\
=-\tr(\tau\DDlambda\DDlambdastar\sigma)
+\tr(\tau\cdot[G(h),\sigma])
+\tr(\sigma\DDlambda\DDlambdastar\tau) 
=
-\tr\bigl(\tau\DDlambda\DDlambdastar\sigma\bigr)
+\tr(\sigma\DDlambda\DDlambdastar\tau)
+\tr\bigl([\sigma,\tau]\cdot G(h)\bigr)
\end{multline} 
Hence, we obtain
$ d_h\Phi(\sigma,\tau)=
-\sqrt{-1}\int_X
\Bigl(
 (-\lambdabar\delbar+\del)\tr(\tau\DDlambda\sigma)
+(\lambda\del+\delbar)\tr(\sigma\DDlambdastar\tau)
\Bigr)\cdot\omega^{n-1}$.
By using $\|\sigma\|_P<\infty$
and $\|\tau\|_{P}<\infty$,
we obtain the vanishing of $d_h\Phi(\sigma,\tau)$,
due to Lemma 5.2 of \cite{s1}.
\hfill\qed

\subsubsection{Donaldson functional}

For $h_1,h_2\in \nbigp_h$,
take a differentiable path $\gamma:[0,1]\lrarr\nbigp_h$ 
such that $\gamma(0)=h_1$ and $\gamma(1)=h_2$,
and the Donaldson functional is defined to be
$ M(h_1,h_2):=\int_{\gamma}\Phi$.
It is independent of a choice of a base metric $\omega$,
in the case $\dim X=1$.
We have
$M(h_1,h_2)+M(h_2,h_3)=M(h_1,h_3)$
by the construction.

\begin{lem}
When $h_2=h_1\cdot e^s$ 
for $s\in \nbigp(S_{h_1})$, we have the following formula:
\begin{equation}
 \label{eq;06.1.3.1}
 M\bigl(h_1,h_2\bigr)=
 \sqrt{-1}\int_X\tr\bigl(s\Lambda_{\omega} G(h_1)\bigr)
 \dvol_{\omega}
+\int_X\bigl(\Psi(s)\DDlambda s,\DDlambda s\bigr)_{\omega,h_1}
 \dvol_{\omega}.
\end{equation}
Here, $(\cdot,\cdot)_{\omega,h_1}$ denotes
the hermitian product induced by $\omega$ and $h_1$,
and $\Psi$ is given  as follows:
\begin{equation}
 \label{eq;06.1.15.50}
 \Psi(t_1,t_2)
=\frac{e^{t_2-t_1}-(t_2-t_1)-1}{(t_2-t_1)^2}.
\end{equation}
See Subsection {\rm\ref{subsubsection;06.1.9.2}}
for the meaning of $\Psi(s)(\DDlambda s)$.
\end{lem}
\pf
Let $M'(h_1,h_2)$ denote the right hand side of 
(\ref{eq;06.1.3.1}).
The following formula immediately follows from the definition:
\[
 \frac{\del}{\del u}
 M'\bigl(h_1e^{ts},h_1e^{(t+u)s}\bigr)_{|u=0}
=\int_X\sqrt{-1}\tr\bigl(s\Lambda_{\omega} G(h_1e^{ts})\bigr).
\]
We also have the following equalities:
\begin{equation}
 \label{eq;06.1.19.10}
 \frac{\del^2}{\del t\del u}
 M'\bigl(h_1e^{ts},h_1e^{(t+u)s}\bigr)_{|u=0}
=\frac{\del^2}{\del t^2}
 M'\bigl(h_1,h_1e^{ts}\bigr)_{|u=0}
=\frac{\del^2}{\del t\del u}
 M'\bigl(h_1,h_1e^{(t+u)s}\bigr)_{|u=0}.
\end{equation}
The second equality can be shown formally.
The first equality can be shown by the argument
in the page 883 of \cite{s1}.
We also have the obvious equality:
\[
 \frac{\del}{\del u}M'(h_1e^{ts},h_1e^{(t+u)s})_{|t=0,u=0}
=\frac{\del}{\del u} M'(h_1,h_1e^{(t+u)s})_{|t=0,u=0}.
\]
Hence, we obtain the following:
\[
 \frac{\del}{\del t} M'(h_1,h_1e^{ts})
=\int_X\sqrt{-1}\tr\bigl(s\Lambda_{\omega} G(h_1e^{ts})\bigr).
\]
Thus, $M'(h_1,h_1e^{s})$ is the integral of $\Phi'$
along the path 
$\gamma(t)=h_1e^{ts}$,
and hence $M'(h_1,h_2)=M(h_1,h_2)$.
\hfill\qed

\begin{rem}
In {\rm\cite{s1}},
the formula {\rm(\ref{eq;06.1.3.1})}
is adopted to be the definition of the functional.
We follow the original definition of Donaldson
{\rm \cite{don3}}.
\hfill\qed
\end{rem}

We obtain the following corollary
due to the positivity of the function $\Psi$.

\begin{cor}
\label{cor;06.1.4.1}
If $\sup |\Lambda_{\omega} G(h)|_h<B$ is satisfied,
we have the following inequality:
\[
 M(h,he^{s})\geq
 \sqrt{-1}\int \tr\bigl(s\Lambda_{\omega} G(h)\bigr) \cdot\dvol_{\omega}
 \geq -B\int |s|_h\cdot\dvol_{\omega}.
\]
In particular, the upper bound of $s$
gives the lower bound of $M(h,he^s)$.
\hfill\qed
\end{cor}

\subsubsection{Main estimate}

The following key estimate is the counterpart of 
Proposition 5.3 in \cite{s1}.
The proof is same.

\begin{prop}
 \label{prop;06.1.4.3}
Fix  $B>0$.
Let $(E,\DDlambda)$ be a flat $\lambda$-connection.
Let $h$ be a hermitian metric of $E$ such that
$\sup\bigl| \Lambda_{\omega} G(h,\DDlambda)\bigr|_h\leq B$.
Let $(E,\DDlambda,h)$ be analytically stable
with respect to $\omega$.
Then, there exist positive constants $C_i$ $(i=1,2)$
with the following property:
\begin{itemize}
\item
Let $s$ be any self-adjoint endomorphism
satisfying
$ \|s\|_{P,h}<\infty$,
$ \tr(s)=0$
and $\sup\bigl|\Lambda_{\omega} G(h\cdot e^s,\DDlambda)\bigr|\leq B$.
Then, the following inequality holds:
\[
 \sup_X |s|_h\leq C_1+C_2\cdot M(h,he^s)
\]
\end{itemize}
\end{prop}
{\bf (Sketch of the proof)}\,\,
The excellent argument given in \cite{s1} 
works in the case of $\lambda$-connection
without any essential change.
Since we would like to use some minor variants of the proposition
(Subsections
\ref{subsubsection;06.1.19.5}--\ref{subsubsection;06.1.19.6}),
we recall an outline of the proof
for the convenience of the reader.
To begin with,
we remark that we have only to show the following inequality
due to Corollary \ref{cor;06.1.4.1}:
\[
 \sup_X|s|_h\leq
 C_1'+C_2'\cdot\max\bigl\{
 0,M(h,he^s) \bigr\},
\]
As is noticed in Subsection \ref{subsubsection;06.1.18.3},
the inequality
$ \Delta^{\lambda}_{\omega}
 \log \tr(e^s)\leq 
 \bigl|\Lambda G(h)\bigr|_h
+\bigl|\Lambda G(he^s)\bigr|_{he^s}
 \leq 2B$
holds.
Hence, there exist some constants $C_i$ $(i=3,4)$
such that the inequality
$\log\tr(e^s)\leq C_3+C_4\cdot \int\log \tr (e^s)$
holds for any $s$ as above,
due to Condition \ref{condition;06.1.10.1}.
Since we have 
$ C_5+C_6\cdot |s|_h
\leq 
 \log \tr e^s
 \leq
 C_7+C_8\cdot |s|_h$
for some positive constants $C_i$ $(i=5,6,7,8)$,
there exist some constants $C_i$ $(i=9,10)$
such that the following holds for any $s$ as above:
\begin{equation}
 \label{eq;06.1.9.3}
 \sup|s|_h\leq C_{9}+C_{10}\cdot \int|s|_h.
\end{equation}

Assume that the claim of the proposition does not hold,
and we will derive a contradiction.
Under the assumption,
either one of the following occurs:
\begin{description}
\item[Case 1.]
 There exists a sequence 
 $\{s_i\in\nbigp(S_h)\,|\,i=1,2,\cdots,\}$ such that
 $\sup |s_i|_{h}\lrarr \infty$ and
 $M(h,he^{s_i})\leq 0$.
\item[Case 2.]
 There exist sequences $\{s_i\in\nbigp(S_h)\}$ 
 and $\{C_{2,i}\in\real\}$ with the following properties: 
\[
 \sup_X |s_i|\lrarr \infty,\quad
 C_{2,i}\lrarr \infty,\quad
 (i\lrarr\infty)
\]
\[
 M(h,he^{s_i})>0,\quad
 \sup|s_i|_h\geq C_{2,i}M(h,he^{s_i})
\]
\end{description}

In both cases, we have $\|s_i\|_{L^1}\lrarr \infty$
due to (\ref{eq;06.1.9.3}).
We put $\ell_i:=\|s_i\|_{L^1}$ and $u_i:=s_i/\ell_i$.
Clearly we have $\|u_i\|_{L^1}=1$,
and uniform boundedness $\sup_X |u_i|<C$
due to (\ref{eq;06.1.9.3}).
In the following,
let $L^2(S_h)$ (resp. $L_1^2(S_h)$)
denote the space of $L^2$-sections 
(resp. $L_1^2$-sections) of $S_h$.
The following lemma is one of the keys
in the proof of Proposition \ref{prop;06.1.4.3}.

\begin{lem}
 \label{lem;06.1.4.2}
After going to an appropriate subsequence,
$\{u_i\}$ weakly converges to some $u_{\infty}\neq 0$
in $L_1^2(S_h)$.
Moreover,
we have the following inequality,
for any $C^{\infty}$-function 
$\Phi:\real\times\real\lrarr\real_{\geq\,0}$
such that
$\Phi(y_1,y_2)
 \leq(y_1-y_2)^{-1}$
for $y_1>y_2$:
\[
 \sqrt{-1}\int\tr\bigl(u_{\infty}\Lambda_{\omega} G(h)\bigr)
+\int_X\bigl(\Phi(u_{\infty})\DDlambda u_{\infty},
 \DDlambda u_{\infty}\bigr)_{h,\omega}\leq 0.
\]
\end{lem}
\pf
By considering $\Phi-\epsilon$ for any small positive number
$\epsilon$,
we have only to consider the case
$\Phi(y_1,y_2)<(y_1-y_2)^{-1}$
for $y_1>y_2$.
In the both cases, we have the inequalities
for some positive constant $C$,
from the formula (\ref{eq;06.1.3.1}):
\[
 \ell_i\sqrt{-1}\int_X\tr\bigl(u_i\Lambda_{\omega} G(h,\DDlambda)\bigr)
+\ell_i^2\int\bigl(\Psi(\ell_iu_i)\DDlambda u_i,\DDlambda u_i\bigr)_h
\leq \ell_i\cdot \frac{C}{C_{2,i}}.
\]
(In the case 1, we take any sequence $\{C_{2,i}\}$
such that $C_{2,i}\lrarr \infty$).
Let $\Phi$ be as above.
Due to the uniform boundedness of $u_i$,
we may assume that $\Phi$ has the compact support.
Then, if $\ell$ is sufficiently large,
we have
$\Phi(\lambda_1,\lambda_2)
<\ell \Psi(\ell\lambda_1,\ell\lambda_2)$.
Therefore, we obtain the following inequality:
\[
 \sqrt{-1}\int_X\tr\bigl(u_i\Lambda_{\omega} G(h,\DDlambda)\bigr)
+\int_X\bigl(\Phi(u_i)\DDlambda u_i,\DDlambda u_i\bigr)_{h,\omega}
\leq \frac{C}{C_{2,i}}.
\]
Since $\sup_X|u_i|$ is bounded independently of $i$,
there exists a function $\Phi$ as above
which satisfies $\Phi(u_i)=c\cdot \id$, moreover,
for some small positive number $c>0$.
Therefore, we obtain the boundedness
of $\{u_i\}$ in $L_1^2$.
By taking an appropriate subsequence,
$\{u_i\}$ is weakly convergent in $L_1^2$.
Let $u_{\infty}$ denote the weak limit.
Let $Z$ be any compact subset of $X$.
Then, $\{u_i\}$ is convergent to $u_{\infty}$ on $Z$
in $L^2$,
and hence $\int_Z|u_i|\to \int_Z|u_{\infty}|$.
Since $\sup|u_i|$ are uniformly bounded,
we obtain $\int_Z|u_{\infty}|\neq 0$,
if the volume of $X-Z$ is sufficiently small.
Thus, $u_{\infty}\neq 0$.
Similarly, we can show the convergence
$\int \tr \bigl(u_i\Lambda G(h,\DDlambda)\bigr)
\lrarr \int \tr \bigl(u_{\infty}\Lambda G(h,\DDlambda)\bigr)$.
Since $\{u_i\}$ are weakly convergent to $u_{\infty}$
in $L_1^2$,
we have the almost everywhere convergence
of $\{u_i\}$ and $\{\DDlambda u_i\}$ 
to $u_{\infty}$ and $\DDlambda u_{\infty}$
respectively.
Therefore, 
the sequence $\{\Phi(u_{i})\DDlambda u_i \}$ converges
to $\Phi(u_{\infty})\DDlambda u_{\infty}$
almost everywhere.
Hence, we have
\[
 \int \bigl(\Phi(u_{\infty})\DDlambda u_{\infty},
 \DDlambda u_{\infty}\bigr)_{h,\omega}
\leq 
\underline{\lim}
 \int \bigl(\Phi(u_{i})\DDlambda u_{i},
 \DDlambda u_{i}\bigr)_{h,\omega}
\]
due to Fatou's lemma.
Thus, we obtain the desired inequality,
and the proof of Lemma \ref{lem;06.1.4.2} is finished.
\hfill\qed

\vspace{.1in}
We reproduce the rest of the excellent argument
given in \cite{s1} just for the completeness.
We do not use it in the later argument.
The point is that we can derive a contradiction
from the existence of the non-trivial section $u_{\infty}$
as in Lemma \ref{lem;06.1.4.2}.

\begin{lem}
The eigenvalues of $u_{\infty}$ are constant,
and $u_{\infty}$ has at least two distinct eigenvalues.
\end{lem}
\pf
To show the constantness of the eigenvalues,
we have only to show the constantness 
of $\tr\bigl(\varphi(u_{\infty})\bigr)$
for any $C^{\infty}$-function $\varphi:\real\lrarr\real$.
We have
$\bigl( \delbar+\lambda\del \bigr)
 \tr\varphi(u_{\infty})
=\tr\bigl(\DDlambda \varphi (u_{\infty})\bigr)
=\tr\bigl(d\varphi(u_{\infty})\DDlambda u_{\infty}\bigr)$.
Let $N$ be any large number.
We can take a $C^{\infty}$-function
$\Phi:\real\times\real\lrarr\real$
such that
$\Phi(y_1,y_1)=d\varphi(y_1,y_1)$
and
$N\Phi^2(y_1,y_2)<(y_1-y_2)^{-1}$
for $y_1>y_2$.
We obtain
$\tr \bigl(
 d\varphi(u_{\infty})(\DDlambda u_{\infty})\bigr)
=\tr\bigl(\Phi(u_{\infty})\DDlambda u_{\infty}\bigr)$
due to the first condition.
We obtain the following inequality from Lemma \ref{lem;06.1.4.2}:
\[
 \int_X|\Phi(u_{\infty})\DDlambda u_{\infty}|^2
\leq 
-\frac{\sqrt{-1}}{N}
 \int_X\tr\bigl(u_{\infty}\Lambda G(h)\bigr).
\]
Therefore,
$\bigl|(\delbar+\lambda\del)\tr\varphi(u_{\infty})
 \bigr|_{L^2}^2 =0$.
Thus, the eigenvalues of $u_{\infty}$ are constant.
Since $\tr(u_{\infty})=0$ and $u_{\infty}\neq 0$,
$u_{\infty}$ has at least two distinct eigenvalues.
\hfill\qed

\vspace{.1in}
Let $\kappa_1<\kappa_2<\cdots<\kappa_w$ denote 
the constant distinct eigenvalues of $u_{\infty}$.
Then, $\varphi(u_{\infty})$ and $\Phi(u_{\infty})$
depend only on the values
$\varphi(\kappa_i)$ and $\varphi(\kappa_i,\kappa_j)$
respectively.

\begin{lem}
 \label{lem;06.1.9.4}
Let $\Phi:\real^2\lrarr\real$ be a $C^{\infty}$-function
such that $\Phi(\kappa_i,\kappa_j)=0$ 
for $\kappa_i>\kappa_j$.
Then, $\Phi(u_{\infty})(\DDlambda u_{\infty})=0$.
\end{lem}
\pf
We may replace $\Phi$ with $\Phi_1$
satisfying $\Phi_1(\kappa_i,\kappa_j)=0$
for $\kappa_i>\kappa_j$
and $N\Phi_1^2(y_1,y_2)<(y_1-y_2)^{-1}$ for $y_1>y_2$.
Then, we obtain
$\bigl\|\Phi_1(u_{\infty})\DDlambda u_{\infty}\bigr\|_{L^2}^2
\leq C/N$
due to Lemma \ref{lem;06.1.4.2},
and hence we obtain
$\Phi(u_{\infty})\DDlambda u_{\infty}
=\Phi_1(u_{\infty})\DDlambda u_{\infty}=0$.
\hfill\qed

\vspace{.1in}

Let $\gamma_i$ denote the open interval
$\openopen{\kappa_i}{\kappa_{i+1}}$.
Let $p_{\gamma}:\real\lrarr [0,1]$ be
any decreasing $C^{\infty}$-function such that
$p_{\gamma}(\kappa_i)=1$  and 
$p_{\gamma}(\kappa_{i+1})=0$.
We put $\pi_{\gamma}=p_{\gamma}(u_{\infty})$.
It is easy to see that $\pi_{\gamma}$ is $L_1^2$.
Due to $p_{\gamma}^2=p_{\gamma}$,
we have $\pi_{\gamma}^2=\pi_{\gamma}$.
We have
$\DDlambda \pi_{\gamma}
=dp(u_{\infty})\DDlambda u_{\infty}$.
We put
$\Phi_{\gamma}(y_1,y_2)
=(1-p_{\gamma})(y_2)\cdot dp_{\gamma}(y_1,y_2)$,
and then we have
$(1-\pi_{\gamma})\circ\DDlambda\pi_{\gamma}
=\Phi_{\gamma}(u_{\infty})\circ\DDlambda u_{\infty}$.
On the other hand,
since we have $\Phi_{\gamma}(\kappa_i,\kappa_j)=0$
$(\kappa_i>\kappa_j)$,
we obtain $\Phi_{\gamma}(u_{\infty})\DDlambda u_{\infty}=0$
due to Lemma \ref{lem;06.1.9.4}.
Therefore, we obtain
$(1-\pi_{\gamma})\circ \DDlambda\pi_{\gamma}=0$.

From $(1-\pi_{\gamma})d''\pi_{\gamma}=0$,
we obtain a saturated coherent subsheaf $V_{\gamma}$
such that $\pi_{\gamma}$ is the orthogonal projection on
$V_{\gamma}$
due to the result of Uhlenbeck-Yau
\cite{uy}.
From $(1-\pi_{\gamma})d'\pi_{\gamma}=0$,
the bundle $V_{\gamma}$  is $\DDlambda$-invariant.
Since we consider the case $\lambda\neq $0,
it is easy to see that $V_{\gamma}$ is indeed a subbundle of $E$.
Namely, we obtain the $\lambda$-flat subbundle
$(V_{\gamma},\DDlambda_{V_{\gamma}})
\subset (E,\DDlambda)$.

Let us show
$\deg_{\omega}(V_{\gamma},
 h_{\gamma})/\rank V_{\gamma}
\geq\deg_{\omega}(E,h)/\rank E$ for some $\gamma$,
which contradicts the stability assumption of $(E,\DDlambda,h)$,
where $h_{\gamma}:=h_{|V_{\gamma}}$.
From Lemma \ref{lem;06.1.10.2}, we have
\[
 \deg(V_{\gamma})=
\frac{1}{2\pi}
\frac{1}{1+|\lambda|^2}
\left(
 \sqrt{-1}\int \tr\bigl(\pi_{\gamma}G(h)\bigr)
-\int |\DDlambda \pi_{\gamma}|^2
\right).
\]
We have 
$u_{\infty}=
 \kappa_w\cdot\id_E-\sum |\gamma|\cdot \pi_{\gamma}$,
where $|\gamma|$ denotes the length of $\gamma$.
We put
\[
 W=\kappa_w\deg(E)-\sum |\gamma|\cdot \deg(V_{\gamma})
=\frac{1}{2\pi}
\frac{1}{1+|\lambda|^2}
\left(
\sqrt{-1}\int \tr \bigl(u_{\infty} \Lambda G(h)\bigr)
+\int \sum |\gamma|\cdot \bigl|\DDlambda\pi_{\gamma}\bigr|^2
\right).
\]
Since 
$\DDlambda\pi_{\gamma}
=dp_{\gamma}(u_{\infty}) \DDlambda u_{\infty}$,
we have
\[
 W=
\frac{1}{2\pi}
\frac{1}{1+|\lambda|^2}
\left(
\sqrt{-1}\int \tr \bigl(u_{\infty} \Lambda G(h)\bigr)
+\int
 \Bigl(\sum |\gamma|\cdot dp_{\gamma}(u_{\infty})^2\cdot
 \DDlambda u_{\infty},\,\,
 \DDlambda u_{\infty}
 \Bigr)
\right).
\]
We can check
$ \sum |\gamma| (dp_{\gamma})(\kappa_i,\kappa_j)^2
=(\kappa_i-\kappa_j)^{-1}$ for $\kappa_i>\kappa_j$
by a direct argument.
Therefore, we obtain $W\leq 0$,
due to Lemma \ref{lem;06.1.4.2}.
Namely we obtain
$\kappa_w\cdot \deg E\leq \sum |\gamma|\cdot \deg (V_{\gamma})$.
On the other hand,
we have 
$0=\tr (u_{\infty})=
\kappa_{w}\cdot\rank E-\sum |\gamma|\cdot\rank V_{\gamma}$.
Therefore, we obtain 
$\deg(V_{\gamma})/\rank V_{\gamma}
\geq \deg(E)/\rank E$ for at least one of $\gamma$,
which contradicts with the stability of
$(E,\DDlambda,h)$.
Thus, the proof of Proposition \ref{prop;06.1.4.3}
is finished.
\hfill\qed

\subsubsection{Variant 1 of Proposition \ref{prop;06.1.4.3}}
\label{subsubsection;06.1.19.5}

Let $C$ be a smooth projective curve,
and $D$ be a simple divisor.
Let $(E,\DDlambda,\vecF)$ be 
a $\lambda$-flat bundle on $(C,D)$.
Let $\eta$ be a sufficiently small positive number
such that $10\cdot\eta<\gap(E,\vecF)$.
Let $\epsilon_0$ be a sufficiently smaller number
than $\eta$,
for example $10\rank(E)\epsilon_0<\eta$.
Let $\omega_{\epsilon}$ $(0\leq \epsilon<\epsilon_0)$
be a Kahler metric of $C-D$ with the following conditions:
\begin{itemize}
\item
 Let $P\in D$.
 Let $(U,z)$ be a holomorphic coordinate around $P$
 such that $z(P)=0$.
 Then,  the following holds
 for some positive constants $C_i$ $(i=1,2)$:
\[
 C_1\cdot\omega_{\epsilon}
\leq 
 \epsilon^2|z|^{2\epsilon}\frac{dz\cdot d\zbar}{|z|^2}
+\eta^2|z|^{2\eta}\frac{dz\cdot d\zbar}{|z|^2}
\leq
 C_2\cdot\omega_{\epsilon}
\]
\item
$\omega_{\epsilon}\lrarr\omega_0$ 
for $\epsilon\to 0$
in the $C^{\infty}$-sense 
locally on $C-D$.
\end{itemize}

Let $\vecF^{(\epsilon)}$ be an $\epsilon$-perturbation
of $\vecF$.
 See Subsection \ref{subsubsection;06.1.19.10}
 for the notion of $\epsilon$-perturbation.
 We discuss the surface case there, but it can be applied
 in the curve case.
Suppose that we are given
hermitian metrics $h^{(\epsilon)}$ for
$(E,\vecF^{(\epsilon)})$ with the following properties:
\begin{itemize}
\item
$ \bigl|\Lambda_{\omega_{\epsilon}}
 G(h^{(\epsilon)},\DDlambda)
 \bigr|_{h^{(\epsilon)}}
\leq C_1$,
where the constant $C_1$ is independent of $\epsilon$.
\item
$\{h^{(\epsilon)}\}$ converges to $h^{(0)}$
for $\epsilon\to 0$
in the $C^{\infty}$-sense locally on $C-D$.
\end{itemize}

\begin{lem}
\label{lem;06.1.10.6}
Let $s^{(\epsilon)}$ be self-adjoint endomorphisms
of $(E,h^{(\epsilon)})$
satisfying $\tr s^{(\epsilon)}=0$
and the following properties:
\begin{itemize}
\item
$\|s^{(\epsilon)}\|_{P,h^{(\epsilon)},\omega_{\epsilon}}<\infty$.
But we do not assume the uniform boundedness.
\item
$\bigl|
 \Lambda_{\omega_{\epsilon}}
 G(h^{(\epsilon)}e^{s^{(\epsilon)}},\DDlambda)
\bigr|_{h^{(\epsilon)}}\leq C_1$.
The constant $C_1$ is independent of $\epsilon$.
\end{itemize}

Then, there exist constants $C_i$ $(i=3,4)$,
which are independent of $\epsilon$,
with the following property:
\[
 \sup |s^{(\epsilon)}|_{h^{(\epsilon)}}
\leq C_3+C_4\cdot 
 M(h^{(\epsilon)},\,\,h^{(\epsilon)}e^{s^{(\epsilon)}}).
\]
\end{lem}
{\bf (Sketch of a proof)}\,\,
The argument is essentially same as the proof of
Proposition \ref{prop;06.1.4.3}.
We assume that the claim does not hold,
and we will derive a contradiction.
After going to an appropriate subsequence,
either one of the following holds:
\begin{description}
\item[Case 1.]
 $M(h^{(\epsilon)},h^{(\epsilon)}e^{s^{(\epsilon)}})\leq 0$
and 
 $\sup_{C-D} |s^{(\epsilon)}|_{h^{(\epsilon)}}
\lrarr \infty$ for $\epsilon\to 0$.
\item[Case 2.]
$M\bigl(h^{(\epsilon)},
 h^{(\epsilon)}e^{s^{(\epsilon)}}\bigr)>0$,\,\,
$\sup|s^{(\epsilon)}|\geq
 C_2^{(\epsilon)}
 M\bigl(h^{(\epsilon)},
 h^{(\epsilon)}e^{s^{(\epsilon)}}\bigr)$,\,\,
 $\sup_{C-D} \bigl|s^{(\epsilon)}\bigr|_{h^{(\epsilon)}}
 \lrarr\infty$ and
$C_{2}^{(\epsilon)}\lrarr\infty$ for $\epsilon\to 0$.
\end{description}

By using Lemma \ref{lem;06.1.10.3} (given below)
and the argument given in the first part of
Proposition \ref{prop;06.1.4.3},
we can show that 
there exist positive constants $C_i$ $(i=5,6)$,
which are independent of $\epsilon$, 
with the following property:
\[
 \sup_{C-D}|s^{(\epsilon)}|_{h^{(\epsilon)}}
\leq C_5+C_6\cdot \int |s^{(\epsilon)}|_{h^{(\epsilon)}}
 \dvol_{\omega_{\epsilon}}.
\]

We put $\ell^{(\epsilon)}:=\|s^{(\epsilon)}\|_{L^1}$
and $u^{(\epsilon)}:=s^{(\epsilon)}/\ell^{(\epsilon)}$.
The following lemma is the counterpart of Lemma \ref{lem;06.1.4.2}.
\begin{lem}
 \label{lem;06.1.10.4}
We have a non-trivial $L_1^2$-section $u_{\infty}$ of $S_{h^{(0)}}$
with the following property:
\begin{itemize}
\item
The following inequality holds
for any $C^{\infty}$-function
$\Phi:\real\times\real\lrarr\real_{\geq\,0}$
such that
$\Phi(y_1,y_2)
 \leq(y_1-y_2)^{-1}$
for $y_1>y_2$:
\[
 \sqrt{-1}\int_{C-D}\tr\bigl(u_{\infty}\Lambda_{\omega_0} G(h^{(0)})\bigr)
 \dvol_{\omega_0}
+\int_{C-D}\bigl(\Phi(u_{\infty})\DDlambda u_{\infty},
 \DDlambda u_{\infty}\bigr)_{h^{(0)},\omega_0}\dvol_{\omega_0}
\leq 0.
\]
\end{itemize}
\end{lem}
\pf
The argument is essentially same as the proof of 
Lemma \ref{lem;06.1.4.2}.
We have the following
for some positive constant $C_5$:
\[
\sqrt{-1}\int_{C-D}
 \tr\bigl(u^{(\epsilon)}\Lambda_{\omega_{\epsilon}}
 G(h^{(\epsilon)})\bigr)\dvol_{\omega_{\epsilon}}
+\int_{C-D}
 \bigl(\Phi(u^{(\epsilon)})\DDlambda u^{(\epsilon)},
 \DDlambda u^{(\epsilon)}
 \bigr)_{h^{(\epsilon)},\omega_{\epsilon}}
\dvol_{\omega_{\epsilon}}
\leq \frac{C_5}{C_{2}^{(\epsilon)}}.
\]
(In the case $1$, we take any sequence $\{C^{(\epsilon)}_{2}\}$
 such that $C_{2}^{(\epsilon)}\lrarr \infty$.)
From this,
we obtain the following boundedness
as in the proof of Lemma \ref{lem;06.1.4.2}:
\[
 \int_{C-D}\bigl|\DDlambda u^{(\epsilon)}\bigr|_{h^{(\epsilon)}}^2
 \dvol_{\omega_{\epsilon}}<C_{10}.
\]

Let us take a sequence of $C^{\infty}$-isometries
$F_{\epsilon}:(E,h^{(\epsilon)})\lrarr (E,h^{(0)})$
which converges to the identity of $E$,
in the $C^{\infty}$-sense locally on $C-D$.
Remark that the sequence
$\{F_{\epsilon}(\DDlambda)\}$
converges to $\DDlambda$ for $\epsilon\to 0$
in the $C^{\infty}$-sense locally on $C-D$.
The sequence $\{F_{\epsilon}(u^{(\epsilon)})\}$ 
is bounded on $L_1^2$ locally on $C-D$.
By going to an appropriate subsequence,
we may assume that
the sequence $\{u^{(\epsilon)}\}$
is weakly convergent in $L_1^2$ locally on $C-D$,
and hence it is convergent in $L^2$ on 
any compact subset $Z\subset C-D$.
Let $u_{\infty}$ denote the weak limit.
We have $\int_{Z}|u^{(\epsilon)}|\lrarr \int_Z|u_{\infty}|$.
Hence
$\int_Z|u_{\infty}|\neq 0$,
when the volume of $C-Z\cup D$ is sufficiently small,
due to the boundedness of
$\bigl\{\sup|u^{(\epsilon)}|\,\big|\,\epsilon>0\bigr\}$.
In particular, $u_{\infty}\neq 0$.
Similarly, we obtain
$ \int_{C-D}
 \tr(u^{(\epsilon)}G(h^{(\epsilon)}))
\lrarr 
 \int_{C-D}
 \tr(u_{\infty}G(h^{(0)}))$.
Since we can derive the almost everywhere convergence
$\Phi(u^{(\epsilon)})\DDlambda u^{(\epsilon)}
\lrarr \Phi(u_{\infty})\DDlambda u_{\infty}$
and $u^{(\epsilon)}\lrarr u_{\infty}$,
we obtain 
$ \int_{C-D}\bigl(
 \Phi(u_{\infty})\DDlambda u_{\infty},
 \DDlambda u_{\infty}
\bigr)
\leq
 \underline{\lim}
 \int_{C-D}\bigl(
 \Phi(u^{(\epsilon)})\DDlambda u^{(\epsilon)},
 \DDlambda u^{(\epsilon)}
 \bigr)$
due to Fatou's lemma.
Thus, the proof of Lemma \ref{lem;06.1.10.4}
is finished.
\hfill\qed

\vspace{.1in}
The rest of the proof of Lemma \ref{lem;06.1.10.6}
is completely same as the argument for Proposition \ref{prop;06.1.4.3}.
\hfill\qed

\vspace{.1in}
We have used the following lemma
in the proof.
\begin{lem}
\label{lem;06.1.10.3}
For any positive number $B$,
there exist positive constants $C_i$  $(i=1,2)$
with the following property:
\begin{itemize}
\item
 Let $\epsilon$ be any positive number such that
$\epsilon<1/2$.
 Let $f$ be any non-negative bounded $C^{\infty}$-function 
 on $C-D$
 such that
 $\Delta_{\omega_{\epsilon}}f\leq B$.
Then, the inequality
$\sup(f)\leq C_1+C_2\int f\cdot\dvol_{\omega_{\epsilon}}$
holds.
\end{itemize}
\end{lem}
\pf
Let $(U_P,z)$ be as above for $P\in D$,
and $U_P^{\ast}:=U_P-\{z=0\}$.
On $U_P^{\ast}$, the inequality
$ \Delta_{\omega_{\epsilon}} f \leq B$
is equivalent to the following:
\begin{equation}
 \label{eq;08.2.4.1}
 \Delta_{g_0} f
 \leq B\cdot
 \left(\epsilon^2\frac{|z|^{2\epsilon}}{|z|^2}
 +\eta^2\frac{|z|^{2\eta}}{|z|^2}\right).
\end{equation}
Here, $g_0:=dz\cdot d\zbar$.
Because of the boundedness of $f$,
(\ref{eq;08.2.4.1}) holds on $U_P$.
(See the proof of Proposition 2.2 of \cite{s1}.)
Then, we obtain the following inequality on $U_P$:
\[
 \Delta_{g_0}
 \bigl(
 f-B\cdot\phi
 \bigr)\leq 0,\quad
\phi=|z|^{2\epsilon}+|z|^{2\eta}.
\]
For any point $Q\in\Delta(P,1/2)$,
we have the following:
\[
 \bigl(f-B\cdot\phi\bigr)(Q)\leq
 \frac{4}{\pi}\int_{\Delta(Q,1/2)}
 \bigl(f-B\cdot\phi\bigr)\cdot\dvol_{g_0}.
\]
Therefore, there exist some constants $C_i$ $(i=3,4)$
which are independent of $\epsilon$,
such that the following holds:
\[
 f(Q)\leq C_3+C_4\int f\cdot\dvol_{\omega_{\epsilon}}.
\]
Thus, we obtain the upper bound of $f(Q)$,
when $Q$ is close to a point of $D$.
We can obtain such an estimate when $Q$ is far from $D$,
similarly and more easily.
\hfill\qed

\subsubsection{Variant 2 of Proposition \ref{prop;06.1.4.3}}
\label{subsubsection;06.1.19.6}

We will use another variant.
Let $\pi:\nbigc\lrarr \Delta$ be a holomorphic family
of smooth projective curves.
Let $\nbigd\subset\nbigc$ be a relative divisor.
Let $(E,\DDlambda,\vecF)$ be a 
logarithmic parabolic $\lambda$-flat bundle on $(\nbigc,\nbigd)$.
We denote the fiber $\pi^{-1}(t)$ by
$\nbigc_t$ for $t\in\Delta$.
The restriction
$(E,\DDlambda,\vecF)_{|\nbigc_t}$ is denoted
by $(E_t,\DDlambda_t,\vecF_t)$.
Let $\omega$ be a metric of the relative tangent bundle
of $\nbigc/\Delta$ such that 
$\omega\sim \eta^2|z|^{2\eta-2}dz\cdot d\zbar$
around $\nbigd$.
Here, $\eta$ denotes a small positive number
such that $10\rank(E)\cdot\eta<\gap(E,\vecF)$,
and $z$ is holomorphic function such that
$z^{-1}(0)=\nbigd$ and $dz\neq 0$.
The restriction $\omega_{|\nbigc_t}$ is denoted by
$\omega_t$ for $t\in\Delta$.
Let $h$ be a $C^{\infty}$-hermitian metric 
of $E$ adapted to $\vecF$
such that
$ \bigl|
 \Lambda_{\omega_{t}}G(\DDlambda_t,h_t)
 \bigr|_{h_{t}}\leq C_1$
for any $t\in \Delta$,
where a constant $C_1$ is independent of $t$,
and $h_t$ denotes the restriction $h_{|\nbigc_t}$.
The following lemma can be shown 
by an argument similar to the proof of
Lemma \ref{lem;06.1.10.3}.
\begin{lem}
There exist positive constants
$C_i$ $(i=3,4)$, which are independent of $t$,
with the following property.
\begin{itemize}
\item
Let $s^{(t)}$ be an element of 
$\nbigp_{h_{t}}(E_{t})$
satisfying
$\tr s^{(t)}=0$,
$ \|s^{(t)}\|_{h_t,P}<\infty$ and
$\bigl|
 \Lambda_{\omega_t}G(\DDlambda_t,h_te^{s^{(t)}})
 \bigr|\leq C_1$.
Then, the inequality
$ \sup|s^{(t)}|
\leq C_3+C_4\cdot M(h_{t},h_t e^{s^{(t)}})$
holds.
\hfill\qed
\end{itemize}
\end{lem}

\subsection{Regular filtered $\lambda$-flat bundles
associated to tame harmonic bundles}
\label{subsection;06.2.4.50}

\subsubsection{Tame pluri-harmonic metric}

Let $X$ be a complex manifold
with a simple normal crossing divisor $D$.
Let $(E,\DDlambda)$ be a $\lambda$-flat bundle on $X-D$.
Let $h$ be a pluri-harmonic metric of  $(E,\DDlambda)$.
Then, we have the induced Higgs bundle
$(E,\delbar_h,\theta_h)$.
Recall the tameness
of pluri-harmonic metric.
Let $P$ be any point of $X$,
and let $(U_P,z_1,\ldots,z_n)$ be a holomorphic coordinate
around $P$ such that $D\cap U_P=\bigcup_{i=1}^l\{z_i=0\}$.
Then, we have the expression:
\[
 \theta=\sum_{i=1}^l f_i\cdot \frac{dz_i}{z_i}
+\sum_{j=l+1}^ng_j\cdot dz_j.
\]
The pluri-harmonic metric $h$ is called tame,
if the coefficients of the characteristic polynomials
$\det(t-f_i)$ and $\det(t-g_j)$ are holomorphic on $U_P$
for any $P$.
A $\lambda$-flat bundle
with tame pluri-harmonic metric
is called a tame harmonic bundle.
Recall that the curve test for tameness is valid.
\begin{prop}
[Corollary 8.7 of \cite{mochi2}]
\label{prop;08.2.10.1}
A pluri-harmonic metric $h$ for $(E,\DDlambda)$
is tame, if and only if $h_{|C}$ is tame
for any closed curve $C\subset X$ transversal with $D$.
\hfill\qed
\end{prop}

From a holomorphic vector bundle $E$
with a hermitian metric $h$,
we obtain the filtered sheaf
$\vecE_{\ast}(h):=
 \bigl(
 \prolongg{\vecc}{E}\,\big|\,
 \vecc\in\real^S \bigr)$
as explained in Subsection 3.5 of \cite{mochi4}.
We recall the following proposition.
\begin{prop}
[Theorem 8.58, Theorem 8.59
 and Corollary 8.89 of \cite{mochi2}]
Let $(E,\DDlambda,h)$ be a tame harmonic bundle
on $X-D$.
Then, $\bigl(\vecE_{\ast}(h),\DDlambda\bigr)$
is a regular filtered $\lambda$-flat bundle.
\hfill\qed
\end{prop}

\subsubsection{One dimensional case}

In the one dimensional case,
Simpson established the Kobayashi-Hitchin correspondence for
parabolic flat bundles and the parabolic Higgs bundles,
i.e.,
$\lambda$-flat bundles in the case $\lambda=0,1$.
His result can be generalized for any $\lambda$.

\begin{prop}
[Simpson, \cite{s2}]
Let $X$ be a smooth projective curve,
and $D$ be a simple divisor of $X$.
Let $\bigl(\vecE_{\ast},\DDlambda\bigr)$ be a
regular filtered $\lambda$-flat bundle on $(X,D)$.
We put $E=\prolongg{\vecc}{E}_{|X-D}$.
The following conditions are equivalent:
\begin{itemize}
\item
 $(\vecE_{\ast},\DDlambda)$ is poly-stable
 with $\pardeg(\vecE_{\ast})=0$.
\item
 There exists a harmonic metric $h$ of
 $(E,\DDlambda)$, which is adapted to the parabolic structure of
 $\vecE_{\ast}$,
 i.e.,
 $\vecE_{\ast}\simeq \vecE_{\ast}(h)$.
\end{itemize}
Moreover,
such a metric is unique up to obvious ambiguity.
Namely, let $h_i$ $(i=1,2)$ be two harmonic metrics
as above.
Then, we have the decomposition of Higgs bundles
$(E,\DDlambda)=\bigoplus (E_a,\DDlambda_a)$ 
satisfying the following:
\begin{itemize}
\item
 The decomposition is orthogonal with respect to
 both of $h_i$.
\item
 The restrictions of $h_i$ to $E_a$ are denoted by $h_{i,a}$.
 Then, there exist positive numbers $b_a$ 
 such that $h_{1,a}=b_{a}\cdot h_{2,a}$.
\hfill\qed
\end{itemize}
\end{prop}

\subsubsection{The projective case}
\label{subsubsection;06.1.26.100}

Let $X$ be a smooth projective variety
with an ample line bundle $L$,
and let $D$ be a simple normal crossing divisor of $X$
with the irreducible decomposition
$D=\bigcup_{i\in S}D_i$.
Let $(E,\DDlambda,h)$ be a tame harmonic bundle
on $X-D$.

\begin{prop}
\label{prop;06.1.23.60}
Let $(\vecE_{\ast},\DDlambda)$ be as above.
\begin{itemize}
\item
 $(\vecE_{\ast},\DDlambda)$
 is $\mu_L$-polystable
 with $\pardeg_L(\vecE_{\ast})=0$.
\item
 Let
 $(\vecE_{\ast},\DDlambda)
=\bigoplus_j (\vecE_{j\,\ast},\DDlambda_j)
 \otimes\cnum^{p(j)}$ be the canonical decomposition
 of $\mu_L$-polystable regular filtered $\lambda$-flat bundle.
 Then, we have the corresponding decomposition
 of the metric $h=\bigoplus h_i\otimes g_i$,
 where $h_i$ denote pluri-harmonic metrics
 of $(E_i,\DDlambda_i)$ adapted to the parabolic structure,
 and $g_i$ denote metrics of $\cnum^{p(i)}$.
\item
We have the vanishings of characteristic numbers:
\[
 \int_X\parch_{2,L}(\vecE_{\ast})
=\int_X\parchern_{1,L}^2(\vecE_{\ast})=0.
\]
\end{itemize}
\end{prop}
\pf
The first two claims can be shown by the same argument
as the proof of Proposition 5.1 of \cite{mochi4}.
The third claim can be shown by an argument
similar to the proof of Proposition 5.3 of \cite{mochi4},
which we explain briefly.
We have only to consider the case $\dim X=2$.
Since $h$ is pluri-harmonic,
we have the equalities
$\tr R(d'',h)=
(1+|\lambda|^2)^{-1}
\tr G(h,\DDlambda)=0$
and $\tr\bigl(R(d'',h)^2\bigr)=
 (1+|\lambda|^2)^{-2}\cdot\tr\bigl(G(h,\DDlambda)^2\bigr)=0$,
due to Lemma \ref{lem;06.1.20.20}
and Lemma \ref{lem;06.1.12.20} on $X-D$.
We also have the norm estimate for 
the holomorphic sections of $\prolongg{\vecc}{E}$.
(It is explained in Subsection 2.5 of \cite{mochi4} for $\lambda=0$.
 Similar claims hold for any $\lambda$,
 as shown in Subsection 13.3 of \cite{mochi2}.)
Then, the same argument
as the proof of Proposition 5.3 works.
\hfill\qed

\begin{prop}
 \label{prop;06.1.23.50}
Let $(\vecE_{\ast},\DDlambda)$
be a regular filtered $\lambda$-flat bundle.
We put $(E,\DDlambda):=(\vecE_{\ast},\DDlambda)_{|X-D}$.
Let $h_a$ $(a=1,2)$ be pluri-harmonic metrics
of $(E,\DDlambda)$ on $X-D$
which is adapted to the parabolic structure.
Then, we have the decomposition
$(E,\DDlambda)=\bigoplus (E_i,\DDlambda)$
with the following properties:
\begin{itemize}
\item
 The decomposition is orthogonal with respect to
 both of $h_a$ $(a=1,2)$.
 Hence, we have the decomposition $h_a=\bigoplus_i h_{a,i}$.
\item
 There exist positive numbers $b_i$
 such that $h_{1,i}=b_i\cdot h_{2,i}$.
\end{itemize}
The decomposition on $X-D$
is prolonged to the decomposition
$(\vecE_{\ast},\DDlambda)
=\bigoplus (\vecE_{i\,\ast},\DDlambda)$ on $X$.
\end{prop}
\pf
Similar to Proposition 5.2 of \cite{mochi4}.
\hfill\qed

\subsection{Some integral for
non-flat $\lambda$-connection on a curve}

Let $Y$ be a smooth projective curve,
and let $D$ be a divisor.
Let $(E,\vecF)$ be a parabolic bundle on $(Y,D)$.
Let $\DDlambda$ be 
a $C^{\infty}$ $\lambda$-connection on $E_{|Y-D}$.
In this subsection,
we do not assume $\DDlambda$ is flat,
i.e., $(\DDlambda)^2$ may not be $0$.
 But, 
 it is assumed to be 
 flat around an appropriate neighbourhood $U_P$
 of each $P\in D$,
 and $(E,\vecF,\DDlambda)_{|U_P}$ is
 a parabolic $\lambda$-flat bundle.
 In particular, 
 we have $\Res_P(\DDlambda)\in \End(E_{|P})$.
 We assume that it is graded semisimple, for simplicity,
i.e.,
the induced endomorphism 
on $\Gr^F(E_{|P})$ is semisimple
for each $P\in D$.
(By using $\epsilon$-perturbation in 
 Subsection \ref{subsubsection;06.1.19.10},
 we can drop the condition.)

For each $P\in D$,
we have the generalized eigen decomposition
$E_{|P}:=\bigoplus \lefttop{P}\EE_{\alpha}$
of $\Res_P(\DDlambda)$.
We also have the filtration $\lefttop{P}F$ of $E_{|P}$.
Let us take a holomorphic frame $\vecv$
of $E_{|U_P}$, which is compatible with 
$(\lefttop{P}\EE, \lefttop{P}F)$.
We put $\alpha(v_i):=\deg^{\EE}(v_i)$
and $a(v_i):=\deg^{F}(v_i)$.
Let $h$ be a $C^{\infty}$-metric of $E_{|Y-D}$
such that 
$h(v_i,v_j)=|z|^{-2a(v_i)}$ $(i=j)$
and $0$ $(i\neq j)$.
Let us decompose $\DDlambda=d''+d'$.
Let us take a $(1,0)$-operator $d_0'$ such that
$d''+d_0'$ is $C^{\infty}$ $\lambda$-connection of $E$ 
on $Y$, not only on $Y-D$.
We also assume $d_0'\vecv=0$.
We put $A:=d'-d_0'$, which is a $C^{\infty}$-section
of $\End(E)\otimes\Omega^{1,0}(\log D)$
on $Y$, and holomorphic around $D$.
We have $\tr\Res_P(A)=\tr\Res_P(\DDlambda)$.

Let $h_0$ be a $C^{\infty}$-metric of $E$ on $Y$
such that $h_0(v_i,v_j)$ is $1$ $(i=j)$ or $0$ $(i\neq j)$
on $U_P$ $(P\in D)$.
Let $s$ be the endomorphism determined by
$h=h_0\cdot s$.
Then, $s$ is described by the diagonal matrix
$\diag\bigl(|z|^{-2a(v_1)},\ldots,|z|^{-2a(v_r)}\bigr)$
with respect to the frame $\vecv$ on $U_P$.

Although $\DDlambda$ is not necessarily flat,
we obtain the operators
$\delta_h'$, $\delta_{h}''$, $\delbar_h$, $\del_h$,
$\theta_h$ and $\theta^{\dagger}_h$
as in Subsection \ref{subsubsection;06.1.10.7}.
We put 
$\wt(E,\vecF,P):=
 \sum_{a\in\Par(E,P)}
 a\cdot\rank\bigl(\lefttop{P}\Gr^F_a(E)\bigr)$.

\begin{lem}
We have the following formula:
\begin{equation}
 \label{eq;06.1.10.11}
 \frac{\sqrt{-1}}{2\pi}\int_Y\delbar\tr\theta
=\frac{\lambda}{1+|\lambda|^2}
 \sum_P
 \Bigl(\lambda^{-1}\cdot\tr\Res_P\DDlambda
 +\wt(E,\vecF,P)\Bigr).
\end{equation}
\end{lem}
\pf
Let $\delta'_{h_0}$ denote the $(1,0)$-operator obtained from
$d''$ and $h_0$ as in Subsection \ref{subsubsection;06.1.10.7}.
Then, we have
\[
 \theta_h=\frac{1}{1+|\lambda|^2} (d'-\lambda\cdot\delta_h')=
\frac{1}{1+|\lambda|^2}(d_0'-\lambda\cdot\delta_{h_0}')
+\frac{1}{1+|\lambda|^2}
 (A-\lambda\cdot s^{-1}\delta_{h_0}'s).
\]
We would like to apply the Stokes formula
to the integral of $\delbar\tr\theta_h$.
If we do so,
$d_0'-\lambda\delta_{h_0}'$ does not contribute,
because it is the $C^{\infty}$-section of 
$\End(E)\otimes\Omega^{1,0}$.
We have
\[
 \frac{\sqrt{-1}}{2\pi} \int_Y\delbar\tr(A)
=\sum_P \tr\Res_P\DDlambda.
\]
Since $s^{-1}\delta_{h_0}'s$ is described
by $\diag\bigl(-a(v_1),\ldots,-a(v_r)\bigr)\!\cdot\! dz/z$
with respect to $\vecv$ on $U_P$ $(P\in D)$,
we have
\[
 \frac{\sqrt{-1}}{2\pi}
 \int_Y\delbar\tr (s^{-1}\delta_{h_0}'s)
=\sum_P\sum_{i=1}^{\rank E}-a(v_i)
=-\sum_P\wt(E,\vecF,P).
\]
Therefore, we obtain the following formula:
\begin{equation}
 \label{eq;06.1.10.16}
\frac{\sqrt{-1}}{2\pi}\frac{1+|\lambda|^2}{\lambda}
 \int\delbar\tr\theta_h
=\sum_P \Bigl(\lambda^{-1}\tr\Res_P\DDlambda
 +\wt(E,\vecF,P)\Bigr).
\end{equation}
Thus,
we obtain (\ref{eq;06.1.10.11}).
\hfill\qed

%% file: 3.tex
We give a construction of an ordinary metric
for a graded semisimple
parabolic $\lambda$-flat bundle on a surface
satisfying SPW-condition,
and we give the estimate for the pseudo curvature.
Then, we obtain the existence result
of Hermitian-Einstein metric 
if such a parabolic $\lambda$-flat bundle is stable.

\subsection{Around the intersection of the divisor}
\label{subsection;08.2.7.10}

\subsubsection{Some estimates}
\label{subsection;08.2.6.30}

We put $X:=\Delta^2_z$,
$D_i:=\{z_i=0\}$,
and $D:=D_1\cup D_2$.
Let $(\vecE_{\ast},\DDlambda)$ be 
a {\em graded semisimple} 
regular filtered $\lambda$-flat bundle on $(X,D)$.
Let $\vecc=(c_1,c_2)\in\real^2$
such that $c_i\not\in\Par(\vecE_{\ast},i)$.
We assume the following:
\begin{description}
\item[(SPW)]
We have a positive integer $m$
and $\gamma_i\in\real$ with $-1/m<\gamma_i\leq 0$,
such that 
$\Par\bigl(\prolongg{\vecc}{E}_{\ast},i\bigr)$
is contained in 
$\bigl\{
 c_i+\gamma_i+p/m\,\big|\,
 p\in\seisuu,\,
 -1<\gamma_i+p/m<0
 \bigr\}$.
(The condition $-1/m<\gamma_i\leq 0$
 is not essential.)
\end{description}

We put $\Xtilde:=\Delta^2_{\zeta}$,
$\Dtilde_i:=\{\zeta_i=0\}$
and $\Dtilde=\Dtilde_1\cup\Dtilde_2$.
Let $\varphi:\Xtilde\lrarr X$ be
the ramified covering given by
$\varphi(\zeta_1,\zeta_2)=(\zeta_1^m,\zeta_2^m)$.
Let $\Gal(\Xtilde/X)$ denote the Galois group
of $\Xtilde/X$.
Recall the construction in \cite{i-s}.
For any $\veca\in\real^2$,
let $\prolongg{\veca}{\Etilde}$
denote the subsheaf of 
$\vecEtilde:=\varphi^{\ast}(\vecE)$
given as follows:
\[
 \prolongg{\veca}{\Etilde}:=
 \bigcup_{\vecn+m\vecd\leq\veca}
 \varphi^{\ast}\bigl(\prolongg{\vecd}{E}\bigr)
 \cdot \prod_{i=1,2}\zeta_i^{-n_i}
\]
Then, it is easy to see that
$\vecEtilde_{\ast}=\bigl(
 \prolongg{\veca}{\Etilde}\,\big|\,
 \veca\in\real^2
 \bigr)$ is a filtered bundle,
and the induced flat $\lambda$-connection
$\DDtilde^{\lambda}$ is regular.
We put $\ctilde_i:=m\cdot (\gamma_i+c_i)$.
By the assumption,
we have 
$ \Par\bigl(\vecEtilde_{\ast},i\bigr)
=\bigl\{
 p+\ctilde_i\,\big|\,p\in\seisuu
 \bigr\}$.

We have the generalized eigen decompositions
$\prolongg{\vecc}{E}_{|D_i}
=\bigoplus\lefttop{i}\EE_{\alpha}$
with respect to $\Res_i(\DDlambda)$.
We have the parabolic filtration $\lefttop{i}F$
of $\prolongg{\vecc}{E}$.
Let $\vecv$ be a frame of $\prolongg{\vecc}{E}$
compatible with $\lefttop{i}F$ and $\lefttop{i}\EE$
$(i=1,2)$.
We put as follows:
\[
 a_i(v_j):=\lefttop{i}\deg^F(v_j)-(c_i+\gamma_i)
 \in \frac{1}{m}\cdot \seisuu_{\leq 0}
\]
Let $\alpha_i(v_j)\in\cnum$
denote the complex number
determined by
$v_{j|D_i}\in \lefttop{i}\EE_{\alpha_i(v_j)}$.
We put 
$\vtilde_j:=
 \varphi^{\ast}(v_j)
 \cdot \prod_{i=1,2}\zeta_i^{ma_i(v_j)}$.
Then, $\vecvtilde=(\vtilde_j)$ gives 
the frame of $\prolongg{\vecctilde}{\Etilde}$.
We put $\beta_i(v_j):=
 m\bigl(\lambda\cdot a_i(v_j)+\alpha_i(v_j)\bigr)$.
Let $\Gamma$ be the diagonal matrix
whose $(j,j)$-entry is
$\sum_{i=1,2}\beta_i(v_j)\cdot d\zeta_i/\zeta_i$.
Let $A$ be determined by
$\DDtildelambda\vecvtilde
=\vecvtilde\cdot A$,
and let $A_0:=A-\Gamma$.
In the following, 
let $F_{\Gamma}\in
 \End\bigl(\prolongg{\vecctilde}{\Etilde}\bigr)
 \otimes\Omega^1(\log \Dtilde)$ be determined by
$F_{\Gamma}(\vecvtilde)=\vecvtilde\cdot \Gamma$.
We put $\DDtildelambda_0:=\DDtildelambda-F_{\Gamma}$.

Let $A_0=\sum_{i=1,2}A_0^i\cdot d\zeta_i$.
If $m$ is sufficiently large,
we may assume the following:
\begin{description}
\item[(A):]
$A_0^i=O(\zeta_i^2)$.
Moreover,
$(A_0^1)_{j,k}=O\bigl(\zeta_1^2\cdot \zeta_2^2\bigr)$
in the case $\beta_2(v_j)\neq\beta_2(v_k)$,
and
$(A_0^2)_{j,k}=O\bigl(\zeta_1^2\cdot \zeta_2^2\bigr)$
in the case $\beta_1(v_j)\neq\beta_1(v_k)$.
\end{description}
Let $\htilde$ be the hermitian metric of 
$\prolongg{\vecctilde}{\Etilde}$ determined by
$\htilde(\vtilde_i,\vtilde_j)=
 \delta_{i,j}\cdot|\zeta_1|^{-2\ctilde_1}
 \cdot |\zeta_2|^{-2\ctilde_2}$.

Let $\thetatilde$ (resp. $\thetatilde_0$)
be the section of 
$\End(\Etilde)\otimes\Omega^1$
on $\Xtilde-\Dtilde$
induced by $\htilde$ and $\DDtildelambda$
(resp. $\DDtildelambda_0$)
as in Subsection \ref{subsubsection;06.1.10.7}.
Let $\thetatilde^{\dagger}$
and $\thetatilde_0^{\dagger}$ denote
the adjoint of $\thetatilde$ and $\thetatilde_0$,
respectively.
Let $\gtilde$ denote the Euclidean metric of $\Xtilde$.

\begin{lem}
\mbox{{}}\label{lem;08.2.6.21}
\begin{itemize}
\item
$\bigl[
 \thetatilde,\thetatilde^{\dagger}
\bigr]$ is bounded
with respect to $\htilde$ and $\gtilde$.
\item
$\thetatilde^2
=O(z_1\cdot z_2)\cdot dz_1\cdot dz_2$.
\end{itemize}
\end{lem}
\pf
We have the relations
$\thetatilde=\thetatilde_0
 +(1+|\lambda|^2)^{-1}F_{\Gamma}$
and $\thetatilde^{\dagger}=\thetatilde_0^{\dagger}
+(1+|\lambda|^2)^{-1}F_{\Gamma}^{\dagger}$.
Hence,
we have the following:
\begin{equation}
 \label{eq;08.2.6.35}
 \bigl[
 \thetatilde,\thetatilde^{\dagger}
 \bigr]
=\bigl[\thetatilde_0,\thetatilde_0^{\dagger}\bigr]
+(1+|\lambda|^2)^{-1}
 \bigl[\thetatilde_0,F_{\Gamma}^{\dagger}\bigr]
+(1+|\lambda|^2)^{-1}
 \bigl[\thetatilde_0^{\dagger},F_{\Gamma}\bigr]
\end{equation}
We have the decomposition of
$\thetatilde_0$ into the sum 
$\lambda(1+|\lambda|^2)^{-1}
 \sum\ctilde_i\cdot d\zeta_i/\zeta_i
+\thetatilde_0'$,
where $\thetatilde_0'$ is the $C^{\infty}$-section 
of $\End(\prolongg{\vecctilde}{\Etilde})
 \otimes\Omega^1_{\Xtilde}$ on $\Xtilde$.
Hence, 
$\bigl[\thetatilde_0,\thetatilde_0^{\dagger}\bigr]$
is the $C^{\infty}$-section of
$\End(\prolongg{\vecctilde}{\Etilde})
 \otimes\Omega^2$ on $\Xtilde$.
By Condition (A),
$\bigl[\thetatilde_0,F_{\Gamma}^{\dagger}\bigr]$ 
and $\bigl[F_{\Gamma},\thetatilde_0^{\dagger}\bigr]$
is also bounded.
We have
$\thetatilde^2
=\thetatilde_0^2
+2\bigl[\thetatilde_0,F_{\Gamma}\bigr]$.
Then, we obtain the desired estimate for
$\thetatilde^2$ by Condition (A).
\hfill\qed

\begin{lem}
\label{lem;08.2.6.41}
We have the boundedness of
$G(\DDtildelambda,\htilde)$ and
$\thetatilde^2\cdot\thetatilde^{\dagger}$
with respect to $\htilde$ and $\gtilde$.
\end{lem}
\pf
The boundedness of 
$\thetatilde^2\cdot\thetatilde^{\dagger}$
follows from the estimate for $\thetatilde^2$.
We have the following equality
(See Subsection \ref{subsection;08.2.6.20}):
\[
 G(\DDtildelambda,\htilde)
=(1+|\lambda|^2)\cdot R(\htilde,d'')
-\frac{(1+|\lambda|^2)^2}{\lambda}
 \bigl(\delbar^2_{\htilde}+\thetatilde^2
 -\lambda\bigl[\thetatilde,\thetatilde^{\dagger}\bigr]\bigr)
\]
We have the vanishing of the curvature
$R(\htilde,d'')=0$,
and the relation
$\lambda^{-1}\delbar^2_{\htilde}
=\lambdabar^{-1}(\thetatilde^{\dagger})^2$.
Hence, we obtain the boundedness
of $G(\DDtildelambda,\htilde)$
from Lemma \ref{lem;08.2.6.21}.
\hfill\qed

\vspace{.1in}

Since $\htilde$ is $\Gal(\Xtilde/X)$-equivariant,
we obtain the induced metric $h$ of $E$ on $X-D$.
Clearly,
$h$ is given by
$h(v_i,v_j)=
 \delta_{i,j}\cdot |z_1|^{-2a(v_i)}
 \cdot |z_2|^{-2a(v_j)}$.
Let $\theta$ be the section of
$\End(E)\otimes\Omega^{1,0}$ on $X-D$
induced by $\DDlambda$ and $h$,
and let $\theta^{\dagger}$ denote the adjoint
of $\theta$.
Let $g_m$ denote the metric of $X-D$
given by
$g_m=\sum
 |z_i|^{2(-1+1/m)}\cdot dz_i\cdot d\zbar_i$.
\begin{cor}
\label{cor;08.2.7.15}
We have the boundedness of
$G(\DDlambda,h)$ and $\theta^2\cdot\theta^{\dagger}$
with respect to
$g_m$ and $h$.
\hfill\qed
\end{cor}

\subsubsection{The induced metric and the $\lambda$-connection
on the divisors}
\label{subsection;08.2.7.201}

For simplicity,
we assume $c_i=\gamma_i=0$ ($i=1,2$)
in this subsection.
Let $(a,\alpha)\in\KMS(\prolong{E},i)$.
Let $\rho$ be a $C^{\infty}$-function on $X$
such that $\rho>0$.
We put $\chi:=\rho\cdot |z_1|^2$.
Let $D_i^{\circ}:=D_i-(D_1\cap D_2)$.
We discuss the induced hermitian metric 
and the induced $\lambda$-connection
of $\lefttop{i}\Gr^{F,\EE}_{a,\alpha}
 \bigl(\prolong{E}\bigr)_{|D_i^{\circ}}$ 
$(i=1,2)$,
depending on the choice of $\rho$.
Let us consider the case $i=1$.
Let $u_j$ $(j=1,2)$ be sections of
$\lefttop{1}\Gr^{F,\EE}_{a,\alpha}(\prolong{E})$.
We take sections 
$u_j'$ $(j=1,2)$ of $\prolong{E}$ which induce $u_j$.
Then, it can be shown that
$\bigl(\chi^a\cdot h_0(u_1',u_2')\bigr)_{|D_1^{\circ}}$
is independent of the choice of $u_j'$,
which is denoted by
$h_{a,\alpha}(u_1,u_2)$.

We have the frame 
$\vecv_{(a,\alpha)}$
induced by $\vecv$ above.
By construction,
$h_{a,\alpha}(v_i,v_j)
=\rho^a\cdot\delta_{i,j}\cdot |z_2|^{-2a_2(v_i)}$.
Hence, the following equality can be checked 
by a direct calculation:
\begin{equation}
\label{eq;08.2.7.31}
 \tr R(h_{a,\alpha})
-a\cdot\rank\Gr^{F,\EE}_{a,\alpha}(\prolong{E})
\cdot \delbar\del\log\rho=0
\end{equation}

Let $F_0$ denote the $C^{\infty}$-section of
$\End(\prolong{E})\otimes\Omega^{1,0}_X(\log D)$
determined by
$F_0(v_j)=v_j\cdot
 \alpha_1(v_j)\cdot \del\log\chi$.
Then, $\DDlambda-F_0$ is $C^{\infty}$
around $D_1^{\circ}$,
whose restriction preserves
the filtration $\lefttop{1}F$
and the decomposition $\lefttop{1}\EE$.
Hence, we obtain the induced $\lambda$-connection
$\DDlambda_{a,\alpha}$
of $\lefttop{1}\Gr^{F,\EE}_{(a,\alpha)}(\prolong{E})$.
We have $\theta_{a,\alpha}$ induced by
$\DDlambda_{a,\alpha}$ and $h_{a,\alpha}$.
\begin{lem}
\label{lem;08.2.7.30}
The following holds:
\begin{equation}
 \label{eq;08.2.7.1}
 \delbar\tr\theta_{a,\alpha}
+\frac{\lambda a+\alpha}{1+|\lambda|^2}
 \rank\Bigl(
 \lefttop{1}\Gr^{F,\EE}_{a,\alpha}(\prolong{E})
\Bigr)\cdot
 \delbar\del\log\rho=0
\end{equation}
\end{lem}
\pf
Let $\DDlambda_{a,\alpha,1}$
and $\theta_{a,\alpha,1}$ denote the operator,
and let $h_{a,\alpha,1}$ denote the metric
in the case where $\rho$ is constantly $1$.
Since $\theta_{a,\alpha,1}$ is holomorphic,
we have $\delbar\tr\theta_{a,\alpha,1}=0$.
Note that we have
$\DDlambda_{a,\alpha}
=\DDlambda_{a,\alpha,1}
-\alpha\cdot \del\log\rho$
and $h_{a,\alpha}=h_{a,\alpha,1}\cdot \rho^a$.
Then, we obtain
$\theta_{a,\alpha,1}
=\theta_{a,\alpha}
+(1+|\lambda|^2)^{-1}
 \bigl(\lambda\cdot a+\alpha\bigr)
\cdot \del\log\rho$.
Thus, we obtain (\ref{eq;08.2.7.1}).
\hfill\qed

\subsection{Around the smooth part of the divisor}
\label{subsection;08.2.7.11}

\subsubsection{Construction of the metric
and some estimates}

Let $X=\Delta^2$ and $D=\{z_1=0\}$.
Let $\rho$ be a positive $C^{\infty}$-function on $X$,
and we put $\chi:=\rho\cdot|z_1|^2$.
Let $(\vecE_{\ast},\DDlambda)$
be a graded semisimple regular filtered $\lambda$-flat bundle 
on $(X,D)$
with rational weights.
We take $c\in\real$ such that 
$c\not\in\Par(\vecE_{\ast})$.
We assume the following:
\begin{description}
\item[(SPW)]
We have a positive integer $m$
and $\gamma\in\real$ with $-1/m<\gamma\leq 0$
such that $\Par\bigl(\prolongg{c}{E}_{\ast}\bigr)$
is contained in 
$\bigl\{
 c+\gamma+p/m\,\big|\,
 p\in\seisuu,\,
 -1<\gamma+p/m<0
 \bigr\}$.
\end{description}

Let $\Xtilde:=\Delta^2_{\zeta}$ and $D=\{\zeta_1=0\}$.
Let $\varphi:\Xtilde\lrarr X$ be given by
$\varphi(\zeta_1,\zeta_2)=(\zeta_1^m,\zeta_2)$.
We have the induced filtered $\lambda$-flat bundle
$(\vecEtilde_{\ast},\DDtildelambda)$ on 
$(\Xtilde,\Dtilde)$
as in Subsection \ref{subsection;08.2.6.30}.
We put $\ctilde:=m\cdot (c+\gamma)$.
Then,
$\Par\bigl(\prolongg{\ctilde}{\Etilde}_{\ast}\bigr)$
is contained in 
$\bigl\{
 \ctilde+p\,\big|\,p\in\seisuu
 \bigr\}$.

We have the generalized eigen decomposition
$\prolongg{c}{E}_{|D}=\bigoplus\EE_{\alpha}$.
We have the filtration $F$ of $\prolongg{c}{E}_{|D}$.
Let $\vecv$ be a frame of $\prolongg{c}{E}$ 
compatible with $F$ and $\EE$.
We put 
$a(v_j):=\deg^F(v_j)-(c+\gamma)$.
Let $\alpha(v_j)\in\cnum$ be determined 
by $v_{j|D}\in\EE_{\alpha(v_j)}$.
We put 
$\vtilde_j:=
 \varphi^{\ast}(v_j)
 \cdot\zeta_1^{m\cdot a(v_j)}$.
Then, $\vecvtilde=(\vtilde_j)$ gives 
the frame of $\prolongg{\ctilde}{\Etilde}$.
Let $\Gamma$ be the diagonal matrix
whose $(j,j)$-entry is given by the following:
\[
 \alpha(v_j)\cdot \del\log(\varphi^{\ast}\chi)
+\lambda \cdot m\cdot a(v_j)\cdot 
\frac{d\zeta_1}{\zeta_1}
\]
Let $A$ be determined by
$\DDtildelambda\vecvtilde
=\vecvtilde\cdot A$,
and let $A_0:=A-\Gamma$.
Let $F_{\Gamma}$ be the 
$C^{\infty}$-section of $\End(\Etilde)\otimes\Omega^{1}$
on $\Xtilde-\Dtilde$,
determined by
$F_{\Gamma}\vecvtilde=\vecvtilde\Gamma$.
We put $\DDtilde_0:=\DDtilde-F_{\Gamma}$.

Let $A_0=A_0^1\cdot d\zeta_1+A_0^2\cdot d\zeta_2$.
If $m$ is sufficiently large,
the following holds:
\begin{description}
\item[(A)]
$A_0^1=O\bigl(|\zeta_1|^2\bigr)$.
Moreover,
$(A_0^2)_{k,l}=O\bigl(|\zeta_1|^2\bigr)$
in the case $\bigl(a(v_k),\alpha(v_k)\bigr)
\neq \bigl(a(v_l),\alpha(v_l)\bigr)$.
\end{description}
Let $\htilde_1$ be the $\Gal(\Xtilde/X)$-equivariant
hermitian metric of $\prolong{\Etilde}$ such that
$\htilde_1(v_i,v_j)=O\bigl(|\zeta_1|^2\bigr)$
in the case $\bigl(a(v_i),\alpha(v_i)\bigr)
\neq \bigl(a(v_j),\alpha(v_j)\bigr)$.
Then,
let $\htilde:=
 \varphi^{\ast}(\chi^{-c-\gamma})\cdot \htilde_1$.

Let $\thetatilde$ (resp. $\thetatilde_0$)
be the section of 
$\End(\Etilde)\otimes\Omega^1$
on $\Xtilde-\Dtilde$
induced by $\htilde$ and $\DDtildelambda$
(resp. $\DDtildelambda_0$)
as in Subsection \ref{subsubsection;06.1.10.7}.
Let $\thetatilde^{\dagger}$
and $\thetatilde_0^{\dagger}$ denote
the adjoint of $\thetatilde$ and $\thetatilde_0$,
respectively.
Let $\gtilde$ denote the Euclidean metric of $\Xtilde$.

\begin{lem}
\mbox{{}}\label{lem;08.2.6.40}
\begin{itemize}
\item
$\bigl[
 \thetatilde,\thetatilde^{\dagger}
\bigr]$ is bounded
with respect to $\htilde$ and $\gtilde$.
\item
$\thetatilde^2
=O\bigl(|z_1|\bigr)\cdot dz_1\cdot dz_2$.
\end{itemize}
\end{lem}
\pf
Similar to Lemma \ref{lem;08.2.6.21}.
\hfill\qed

\begin{lem}
We have the boundedness of
$G\bigl(\DDtildelambda,\htilde\bigr)$ and
$\thetatilde^2\cdot\thetatilde^{\dagger}$
with respect to $\htilde$ and $\gtilde$.
\end{lem}
\pf
It follows from Lemma \ref{lem;08.2.6.40}.
See the proof of Lemma \ref{lem;08.2.6.41}.
\hfill\qed

\vspace{.1in}

We have the induced hermitian metric $h$ of $E$ on $X-D$.
It is adapted to the parabolic structure of $E$.
Let $\theta$ denote the section of
$\End(E)\otimes\Omega^1_{X-D}$
induced by $h$ and $\DDlambda$,
and let $\theta^{\dagger}$ denote the adjoint.
Let $g_m$ denote the metric of $X-D$
given by
$g_m=|z_1|^{-2+2/m}\cdot dz_1\cdot d\zbar_1
+dz_2\cdot d\zbar_2$.

\begin{cor}
\label{cor;08.2.7.16}
We have the boundedness of
$G(\DDlambda,h)$ and 
$\theta^2\cdot\theta^{\dagger}$
with respect to $h$ and $g_m$.
\hfill\qed
\end{cor}

\subsubsection{The induced metric
and the $\lambda$-connections}
\label{subsection;08.2.7.200}

For simplicity,
we assume $c=\gamma=0$ in this subsection.
Let $(a,\alpha)\in\KMS(\prolong{E}_{\ast})$.
We discuss the induced hermitian metric 
and the induced $\lambda$-connection
of $\Gr^{F,\EE}_{a,\alpha}(\prolong{E})$.
Let $u_j$ $(j=1,2)$ be sections of
$\Gr^{F,\EE}_{a,\alpha}(\prolong{E})$.
We take sections 
$u_j'$ $(j=1,2)$ of $\prolong{E}$ which induce $u_j$.
Then, it can be shown that
$\big(\chi^a\cdot h_0(u_1',u_2')\bigr)_{|D}$
is independent of the choice of $u_j'$,
which is denoted by
$h_{a,\alpha}(u_1,u_2)$.

On the other hand,
let $U_{a,\alpha}$ be the subbundle of
$\prolong{\Etilde}$ generated by
$\vtilde_j$ such that 
$\bigl(a(v_j),\alpha(v_j)\bigr)=(a,\alpha)$.
It is easy to see that
the restriction
$U_{a,\alpha|\Dtilde}$ are independent 
of the choice of the frame $\vecv$,
and $U_{a,\alpha|\Dtilde}$ are orthogonal
with respect to $\htilde_{|\Dtilde}$.
The induced metric of $U_{a,\alpha|\Dtilde}$
is denoted by $h'_{a,\alpha}$.

\begin{lem}
\label{lem;08.2.6.101}
Let $R(h_{a,\alpha})$ and 
$R(h'_{a,\alpha})$ denote the curvatures
of $\bigl(\Gr^{F,\EE}_{a,\alpha}(E),h_{a,\alpha}\bigr)$
and $\bigl(U_{a,\alpha|\Dtilde},h'_{a,\alpha}\bigr)$.
Then, we have the following relation:
\begin{equation}
 \label{eq;08.2.6.50}
 \tr\bigl( R(h'_{a,\alpha}) \bigr)
=\tr\bigl(R(h_{a,\alpha})\bigr)
-a\cdot\rank \Gr^{F,\EE}_{a,\alpha}(E)
 \cdot\delbar\del\log\rho
\end{equation}
\end{lem}
\pf
We take the isomorphism
$\Phi:\Gr^{F,\EE}_{a,\alpha}(E)
\simeq
 U_{a,\alpha|\Dtilde}$ given as follows.
Let $v$ be a section of 
$\Gr^{F,\EE}_{a,\alpha}(E)$.
Let $v'$ be a section of $\prolong{E}$
which induces $v$.
Then, $\Phi(v):=\bigl(
 \varphi^{\ast}(v')\cdot z_1^{m\cdot a}
\bigr)_{|\Dtilde}$ is contained in
$U_{a,\alpha|\Dtilde}$,
and independent of the choice of $v'$.
Under the isomorphism,
we have $h'_{a,\alpha}=h_{a,\alpha}\cdot \rho^{-a}$.
Then, (\ref{eq;08.2.6.50}) follows.
\hfill\qed

\vspace{.1in}

We have the induced $\lambda$-connection,
once we fix $\chi$.
(See \cite{b}.)
Let $f$ be any section of 
$\Gr^{F,\EE}_{a,\alpha}(E)$.
Let $\ftilde$ be a lift of $f$ to $\prolong{E}$.
We put 
$\DDlambda f-\alpha\cdot \log\chi\cdot f
=:G_1\cdot (dz_1/z_1)+G_2\cdot dz_2$.
Then, $G_{1|D}$ is contained in $F_{<a}(E)$.
Hence, $G_2\cdot dz_2$ induces
the well defined section of
$\Gr^{F,\EE}_{a,\alpha}(E)\otimes\Omega_{D}^1$,
which is $\DD_{a,\alpha}(f)$.
We have the induced section
$\theta_{a,\alpha}$ of 
$\End\bigl(\Gr^{F,\EE}_{a,\alpha}(E)\bigr)
 \otimes\Omega_D^1$.

\begin{lem}
\label{lem;08.2.6.115}
We have the following relation:
\begin{equation}
 \label{eq;08.2.6.120}
 \tr\bigl(R(h'_{a,\alpha})\bigr)
=-\frac{1+|\lambda|^2}{\lambda}
 \left(\delbar\tr\theta_{a,\alpha}
 +\frac{(a\lambda+\alpha)\cdot
  \delbar\del\log \rho}{1+|\lambda|^2}
 \rank\Gr^{F,\EE}_{a,\alpha}(E)
 \right)
\end{equation}
\end{lem}
\pf
We have the relation:
\begin{equation}
 \label{eq;08.2.6.100}
 R\bigl(\htilde\bigr)
=-\frac{1+|\lambda|^2}{\lambda}d''\thetatilde
=-\frac{1+|\lambda|^2}{\lambda} 
\left(d''\thetatilde_0
 +\frac{1}{1+|\lambda|^2}d''F_{\Gamma}
\right)
\end{equation}
Let $\DD^{\lambda\prime}_{a,\alpha}$
be the induced $\lambda$-connection
of $U_{a,\alpha|\Dtilde}$.
Let $\theta_{a,\alpha}'$ be the section
of $\End(U_{a,\alpha|\Dtilde})\otimes\Omega^1_{\Dtilde}$
induced by $\DD^{\lambda\prime}_{a,\alpha}$
and $h'_{a,\alpha}$.
Then, we obtain the following equality
from (\ref{eq;08.2.6.100}):
\begin{equation}
 \label{eq;08.2.6.110}
 \tr\bigl( R(h'_{a,\alpha}) \bigr)
=-\frac{1+|\lambda|^2}{\lambda}
\left(
 \delbar\tr\theta'_{a,\alpha}
+\frac{1}{1+|\lambda|^2}
 \cdot \alpha \cdot\delbar\del\log\rho
 \cdot\rank\Gr^{F,\EE}_{a,\alpha}(E)
\right)
\end{equation}
Under the isomorphism $\Phi$ in the proof of
Lemma \ref{lem;08.2.6.101},
we have the 
$\DD^{\lambda\prime}_{a,\alpha}
=\DDlambda_{a,\alpha}$.
Because of $h'_{a,\alpha}=h_{a,\alpha}\cdot \rho^{-a}$,
we have
$\theta'_{a,\alpha}
=\theta_{a,\alpha}
+a\lambda(1+|\lambda|^2)^{-1}\del\log\rho$.
Therefore, the right hand side of (\ref{eq;08.2.6.110})
is the same as (\ref{eq;08.2.6.120}).
\hfill\qed

\subsection{An ordinary metric}
\label{subsection;08.2.7.25}

\subsubsection{Setting}

Let $X$ be a smooth projective surface,
and let $D$ be a simple normal crossing divisor
with the irreducible decomposition $D=\bigcup_{i\in S}D_i$.
Let $L$ be an ample line bundle on $X$,
and $\omega$ be a Kahler form which represents $c_1(L)$.
We take a hermitian metric $g_i$ of $\nbigo(D_i)$.
The canonical section $\nbigo\lrarr\nbigo(D_i)$
is denoted by $\sigma_i$.

Let $\epsilon$ be any number such that
$0<\epsilon<1/2$.
Let us fix a sufficiently large number $N$,
for example $N>10$.
We put as follows,
for some positive number $C>0$:
\begin{equation}
 \label{eq;06.1.20.100}
 \omega_{\epsilon}:= \omega+
\sum_i C\cdot \epsilon^N\cdot
       \sqrt{-1}\del\delbar |\sigma_i|_{g_i}^{2\epsilon}.
\end{equation}
It can be shown that
$\omega_{\epsilon}$ are Kahler metrics of $X-D$
for any $0<\epsilon<1/2$,
if $C$ is sufficiently small.

\begin{rem}
The factor $\epsilon^N$ is added
for our later discussion
(Subsection {\rm\ref{subsection;06.1.21.31}}).
\hfill\qed
\end{rem}

\begin{rem}
\label{rem;08.2.7.101}
Let $\tau$ be a closed $2$-form on $X-D$
which is bounded with respect to $\omega_{\epsilon}$.
Then, the following formula holds:
\[
 \int_{X-D}\omega\cdot\tau
=\int_{X-D}\omega_{\epsilon}\cdot\tau.
\]
In particular,
we also have
$\int_{X-D}\omega^2=\int_{X-D}\omega_{\epsilon}^2$.
\hfill\qed
\end{rem}

In the case $\epsilon=1/m$ for 
some positive integer $m$,
it can be shown that
the metric $\omega_{\epsilon}$ satisfies
Condition {\rm\ref{condition;06.1.10.1}}.
The Kahler forms $\omega_{\epsilon}$ behave as follows
around any point of $D$,
which is clear from the construction:
\begin{itemize}
\item
Let $P$ be any point of $D_i\cap D_j$,
and $(U_P,z_i,z_j)$ be a holomorphic coordinate
around $P$ such that 
$D_i\cap U_P=\{z_i=0\}$
and $D_j\cap U_P=\{z_j=0\}$.
Then, there exist positive constants $C_i$ $(i=1,2)$
such that the following holds on $U_P$,
for any $0<\epsilon<1/2$:
\[
 C_1\cdot \omega_{\epsilon}
\leq
\sqrt{-1}\cdot\epsilon^{N+2}\cdot
 \left(
 \frac{dz_i\cdot d\bar{z}_i}{|z_i|^{2-2\epsilon}}
+\frac{dz_j\cdot d\bar{z}_j}{|z_j|^{2-2\epsilon}}
 \right)
+\sqrt{-1}\bigl(
dz_i\cdot d\bar{z}_i
+dz_j\cdot d\bar{z}_j\bigr)
\leq
 C_2\cdot\omega_{\epsilon}.
\]
\item
Let $Q$ be any point of $D_i\setminus \bigcup_{j\neq i}D_j$,
and $(U,w_1,w_2)$ be a holomorphic coordinate around $Q$
such that $U\cap D_i=\{w_1=0\}$.
Then, there exist positive constants $C_i$ $(i=1,2)$
such that the following holds for any $0<\epsilon<1/2$ on $U$:
\[
 C_1\cdot \omega_{\epsilon}
\leq
\sqrt{-1}\cdot\epsilon^{N+2}\cdot
 \left(
 \frac{dw_1\cdot d\bar{w}_1}{|w_1|^{2-2\epsilon}}
 \right)
+\sqrt{-1}\bigl(
  dw_1\cdot d\bar{w}_1
+dw_2\cdot d\bar{w}_2\bigr)
\leq
 C_2\cdot\omega_{\epsilon}.
\]
\end{itemize}

\subsubsection{Construction and some property}
\label{subsection;08.2.7.26}

Let $(\vecE_{\ast},\DDlambda)$ be a 
graded semisimple parabolic $\lambda$-flat bundle.
For simplicity, we consider only the case $\lambda\neq 0$.
We take $\vecc\in\real^S$ 
such that $c_i\not\in\Par(\vecE_{\ast},i)$
for each $i\in S$.
We assume the following:
\begin{description}
\item[(SPW)]
We have a positive integer $m$
and $\gamma_i\in\real$ $(i\in S)$
with $-1/m<\gamma_i\leq 0$,
such that 
$\Par\bigl(\prolongg{\vecc}{E}_{\ast},i\bigr)$
is contained in 
$\bigl\{
 c_i+\gamma_i+p/m\,\big|\,
 p\in\seisuu,\,
 -1<\gamma_i+p/m<0
 \bigr\}$.
\end{description}
Let $\epsilon=m^{-1}$.
Let $h_0$ be a $C^{\infty}$-hermitian metric of $E$ on $X-D$
as in Subsection \ref{subsection;08.2.7.10} around
the intersection points of $D$,
and as in Subsection \ref{subsection;08.2.7.11}
around the smooth points of $D$.
Let $\theta_0$ denote the section of
$\End(E)\otimes\Omega^{1,0}$ on $X-D$
induced by $\DDlambda$ and $h_0$,
and let $\theta_0^{\dagger}$ denote 
the adjoint.

\begin{lem}
\label{lem;08.2.7.20}
We have the boundedness of
$G(\DDlambda,h_0)$ 
and $\theta_0^2\cdot\theta_0^{\dagger}$
with respect to
$h_0$ and $\omega_{\epsilon}$.
\end{lem}
\pf
It follows from Corollary \ref{cor;08.2.7.15}
and Corollary \ref{cor;08.2.7.16}.
\hfill\qed

\begin{cor}
The following equality holds:
\[
 \int_{X-D} \tr \bigl( R(h_0)^2\bigr)
=\frac{1}{(1+|\lambda|^2)^2}
 \int_{X-D}\tr \bigl(G(h_0)^2\bigr).
\]
As a result, we have the following formula:
\begin{equation}
 \left(
 \frac{\sqrt{-1}}{2\pi}
 \right)^2
 \frac{1}{(1+|\lambda|^2)^2}
\int _{X-D} \tr \bigl(G(h_0)^2\bigr)
=2\int_X\parch_2(\vecE_{\ast}).
\end{equation}
\end{cor}
\pf
The second equality follows from the first equality
and the equality (36)
in the proof of Proposition 4.18 of \cite{mochi4}.
Due to Lemma \ref{lem;06.1.12.20},
we have only to show the vanishing
$\int \delbar
 \tr\bigl(\theta_{0}^2\cdot\theta^{\dagger}_{0}\bigr)=0$,
which follows from the estimate
of $\theta_{0}^2\cdot \theta^{\dagger}_0$
in Lemma \ref{lem;08.2.7.20}.
\hfill\qed

\vspace{.1in}
We can also show the following equality
by using Lemma 4.16 of \cite{mochi4}
and the equality $\tr G(h_0)=(1+|\lambda|^2)\cdot \tr R(h_0)$:
\[
\left(
 \frac{\sqrt{-1}}{2\pi}
\right)^2\int_{X-D}
\left(
 \frac{\tr G(h_0)}{1+|\lambda|^2}
\right)^2
=\left(
 \frac{\sqrt{-1}}{2\pi}
\right)^2\int_{X-D}
 \bigl( \tr R(h_0)\bigr)^2
=\int_X \parchern_1(\vecE_{\ast})^2.
\]

Let $V\subset E$ be a $\lambda$-flat subbundle.
Because of $\lambda\neq 0$ and the regularity,
we have the saturated filtered $\lambda$-flat subsheaf
$\vecV_{\ast}\subset\vecE_{\ast}$.
Let $h_V$ be the metric of $V$ induced by $h_0$.
\begin{lem}
\label{lem;08.2.7.100}
We have
$\deg_{\omega_{\epsilon}}(V,h_V)
=\pardeg_{\omega}(\vecV_{\ast})$.
In particular,
$\deg_{\omega_{\epsilon}}(E,h_0)
=\pardeg_{\omega}(\vecE_{\ast})$.
\end{lem}
\pf
It can be shown 
by the same argument 
as the proof of Lemma 4.20 of \cite{mochi4}.
\hfill\qed

\subsubsection{The induced metric and 
the $\lambda$-connection on $D_i^{\circ}$}

For simplicity,
we assume $c_i=\gamma_i=0$ $(i\in S)$
in this subsection.
We put $\nbigs(D_i):=
 D_i\cap \bigcup_{j\neq i}D_j$
and $D_i^{\circ}:=D_i\setminus \nbigs(D_i)$.
Let $(a,\alpha)\in
 \KMS\bigl(\prolong{E},\vecF,i\bigr)$.
We have the naturally induced parabolic flat bundle
$\lefttop{i}\Gr^{F,\EE}_{a,\alpha}
 (\prolong{E})_{\ast}$
on $\bigl(D_i,\nbigs(D_i)\bigr)$.
By using the functions $|\sigma_i|_{g_i}^2$,
we obtain the induced hermitian metric
$\lefttop{i}h_{a,\alpha}$ 
and the $\lambda$-connection
$\lefttop{i}\DDlambda_{a,\alpha}$ of 
$\lefttop{i}\Gr^{F,\EE}_{a,\alpha}
 \bigl(\prolong{E}\bigr)_{|D_i^{\circ}}$,
as explained in Subsection \ref{subsection;08.2.7.200}.
(See also Subsection \ref{subsection;08.2.7.201}.)
Let $\tau_i:=\delbar\del\log|\sigma_i|_{g_i}^2$.

\begin{lem}
\label{lem;08.2.7.40}
We have the following equality:
\[
 \tr\bigl(
 R(h_{a,\alpha})
 \bigr)
=-\frac{1+|\lambda|^2}{\lambda}
 \delbar\bigl( \lefttop{i}\theta_{a,\alpha}\bigr)
-\lambda^{-1}\alpha\cdot
 \tau_i\cdot
 \rank \lefttop{i}\Gr^{F,\EE}_{a,\alpha}(\prolong{E})
\]
\end{lem}
\pf
It follows from 
(\ref{eq;08.2.7.31}),
(\ref{eq;08.2.7.1}),
(\ref{eq;08.2.6.50})
and (\ref{eq;08.2.6.120}).
\hfill\qed

\vspace{.1in}

\begin{cor}
We have the following equalities:
\begin{multline}
\label{eq;08.2.8.1}
 \pardeg_{D_i}
\Bigl(
 \lefttop{i}\Gr^{F,\EE}_{a,\alpha}(\prolong{E})_{\ast}
\Bigr) 
=-\sum_{P\in\nbigs(D_i)}
 \Bigl(
 \Re\bigl(
 \lambda^{-1}
 \tr\bigl(\Res_P(\lefttop{i}\DDlambda_{a,\alpha})
 \bigr)\bigr) 
+\wt\bigl(\lefttop{i}\Gr^{F,\EE}_{a,\alpha}
 (\prolong{E})_{\ast},P \bigr)
 \Bigr) \\
-\Re(\lambda^{-1}\alpha)\cdot
\rank\lefttop{i}
 \Gr^{F,\EE}_{a,\alpha}(\prolong{E})
 \cdot [D_i]^2
\end{multline}
\begin{equation}
\label{eq;08.2.8.2}
0=\sum_{P\in\nbigs(D_i)}
 \Image\Bigl(
 \lambda^{-1}
 \tr\bigl(\Res_P(\lefttop{i}\DDlambda_{a,\alpha})
 \bigr)\Bigr) 
+\Image(\lambda^{-1}\alpha)\cdot
\rank\lefttop{i}
 \Gr^{F,\EE}_{a,\alpha}(\prolong{E})
 \cdot [D_i]^2
\end{equation}
\end{cor}
\pf
It follows from Lemma \ref{lem;08.2.7.40}
and (\ref{eq;06.1.10.11}).
\hfill\qed

\begin{rem}
Although we have assumed that
graded semisimplicity and (SPW)-condition 
for $(\vecE_{\ast},\DDlambda)$,
the formulas {\rm(\ref{eq;08.2.8.1})} 
and {\rm(\ref{eq;08.2.8.2})}
without the assumption,
because 
the general case can be reduced to the above special case
by perturbation explained in Subsection 
{\rm\ref{subsubsection;06.1.19.10}}.
\hfill\qed
\end{rem}

\subsection{Preliminary existence result
 of a Hermitian-Einstein metric}

\subsubsection{Hermitian-Einstein metric
 for graded semisimple $\lambda$-flat parabolic bundle
 on surface}

We use the setting 
in Subsection \ref{subsection;08.2.7.25}.
Let $X$ be a smooth projective surface
with an ample line bundle $L$
and a simple normal crossing divisor $D$.
Let $\omega$ be a Kahler form representing $c_1(L)$.
Let $(\vecE_{\ast},\DDlambda)$ be a 
graded semisimple 
regular filtered $\lambda$-flat bundle on $(X,D)$.
We assume 
the (SPW)-condition in Subsection \ref{subsection;08.2.7.26}.
Let $\epsilon=m^{-1}$,
and let $\omega_{\epsilon}$ be the Hermitian metric 
given in (\ref{eq;06.1.20.100}).
We have an ordinary metric $h_0$ constructed
in Subsection \ref{subsection;08.2.7.26}.

\begin{lem}
 \label{lem;06.1.13.200}
We can construct a hermitian metric $h_{in}$
for $E_{|X-D}$ which satisfies the following conditions:
\begin{itemize}
\item
 $h_{in}$ is adapted to $\vecE_{\ast}$.
  More strongly, $h_{in}=h_0\cdot e^{\chi}$ 
 for some function $\chi$
 such that $\chi$, $\del\chi$ and $\delbar\del\chi$
 are bounded with respect to $\omega_{\epsilon}$.
\item
 $G(h_{in},\DDlambda)$
 is bounded with respect to $h_{in}$ and $\omega_{\epsilon}$.
\item
 Let $\vecV_{\ast}$ be a $\lambda$-flat
 filtered subsheaf of $\vecE_{\ast}$.
 Let $V:=\vecV_{|X-D}$ and 
 let $h_{in,V}$ denote the induced metric of $V$.
 Then, we have
 $\pardeg_{\omega}(\vecV_{\ast})
 =\deg_{\omega_{\epsilon}}(V,h_{in,V})$.
\item
 $\tr G(h_{in},\DDlambda)\cdot\omega_{\epsilon}
 =(1+|\lambda|^2)\cdot a\cdot\omega_{\epsilon}^2$ 
 for some constant $a$.
 The constant $a$ is determined by the following condition:
\begin{equation}
 \label{eq;06.1.13.201}
 a\cdot\frac{\sqrt{-1}}{2\pi}\int_{X-D}\omega_{\epsilon}^2
=a\cdot\frac{\sqrt{-1}}{2\pi}\int_X\omega^2
=\pardeg_{\omega}(\vecE_{\ast}).
\end{equation}
\item The following equalities hold:
\[
 \left(\frac{\sqrt{-1}}{2\pi}\right)^2
 \int_{X-D}\frac{\tr \bigl(G(h_{in})^2\bigr)}{(1+|\lambda|^2)^2}
=\int_X2\parch_{2}(\vecE_{\ast}),
\quad
 \left(\frac{\sqrt{-1}}{2\pi}\right)^2
 \int_{X-D}\frac{ \tr\bigl(G(h_{in})\bigr)^2}{(1+|\lambda|^2)^2}
=\int_X\parchern_{1}^2(\vecE_{\ast}).
\]
\item
Let $s$ be determined by $h_{in}=h_0\cdot s$.
Then, $s$ and $s^{-1}$ are bounded,
and $\DDlambda s$ is $L^2$ with respect to $h_0$
and $\omega_{\epsilon}$.
\end{itemize}
Due to the third condition,
 $(E,h_{in},\theta)$ is analytic stable 
with respect to $\omega_{\epsilon}$,
 if and only if $(\vecE_{\ast},\DDlambda)$ is $\mu_L$-stable.
The metric $h_{in}$ is called an initial metric.
\end{lem}
\pf
Let $\chi$ be a positive-valued function $\chi$
such that
$\tr G(h_0)\cdot\omega_{\epsilon}
=a\cdot\omega_{\epsilon}^2$ holds.
We put $h_{in}:=h_0\cdot e^{\chi}$.
By construction, the fourth condition is satisfied.
The other property can be reduced
to the property for $h_0$,
as in Lemma 6.3 of \cite{mochi4}.
\hfill\qed

\begin{prop}
\label{prop;06.1.18.5}
There exists a hermitian metric $h_{HE}$ of $(E,\DDlambda)$
with respect to $\omega_{\epsilon}$
satisfying the following properties:
\begin{itemize}
\item
Hermitian-Einstein condition
$\Lambda_{\omega_{\epsilon}}G(h_{HE})=a$ holds
for the constant $a$ determined by {\rm(\ref{eq;06.1.13.201})}.
\item
$\pardeg_{L}(\vecE_{\ast})
 =\deg_{\omega}(E,h_{HE})$.
\item
We have the following formulas:
\begin{equation}
 \label{eq;06.1.20.150}
\left(
\frac{\sqrt{-1}}{2\pi}
\right)^2
\int_{X-D}
 \frac{\tr\bigl(G(h_{HE})^{\bot\,2}\bigr)}
 {(1+|\lambda|^2)^2}
=\int_X\left(
2\parch_2(\vecE_{\ast})
-\frac{\parchern_1^2(\vecE_{\ast})}{\rank E}
\right)
\end{equation}
\begin{equation}
 \left(\frac{\sqrt{-1}}{2\pi}\right)^2
 \int_{X-D}\frac{\tr \bigl(G(h_{HE})^2\bigr)}{(1+|\lambda|^2)^2}
=\int_X2\parch_2(\vecE_{\ast}). 
\end{equation}
\item
$h_{HE}$ is adapted to $\vecE_{\ast}$,
i.e.,
$\vecE_{\ast}(h_{HE})\simeq \vecE_{\ast}$.
More strongly,
let $s$ be determined by $h_{HE}=h_{in}\cdot s$.
Then, $s$ and $s^{-1}$ are bounded with respect to $h_{in}$,
and $\DDlambda s$ is $L^2$ 
with respect to $h_{in}$ and $\omega_{\epsilon}$.
\end{itemize}
\end{prop}
\pf
It follows from Lemma \ref{lem;06.1.13.200}
and Proposition \ref{prop;06.1.13.250}.
\hfill\qed

\subsubsection{Bogomolov-Gieseker inequality}

Let $Y$ be a smooth projective variety
of {\em any} dimension.
Let $L$ be an ample line bundle on $Y$,
and let $D$ be a simple normal crossing divisor.

\begin{cor}
 \label{cor;06.1.13.150}
Let $(\vecE_{\ast},\DDlambda)$ be a $\mu_L$-stable 
regular filtered $\lambda$-flat bundle on $(Y,D)$
in codimension two.
Then, Bogomolov-Gieseker inequality
holds for $\vecE_{\ast}$.
Namely, we have the following inequality:
\[
 \int_Y\parch_{2,L}(\vecE_{\ast})
\leq
 \frac{\int_Y\parchern_{1,L}^2(\vecE_{\ast})}{2\rank E}.
\]
\end{cor}
\pf
Similar to Theorem 6.1 of \cite{mochi4}.
Namely, since we have the Mehta-Ramanathan type theorem
(Proposition \ref{prop;06.1.18.6}),
we have only to prove the claim in the case $\dim Y=2$.
Due to the method of perturbation of parabolic structure,
we have only to prove the inequality in the case
$(\vecE_{\ast},\DDlambda)$ is a graded semisimple
$\mu_L$-stable regular filtered $\lambda$-flat bundle
on $(Y,D)$, satisfying (SPW)-condition.
Then we can take a Hermitian-Einstein metric $h_{HE}$
as in Proposition \ref{prop;06.1.18.5},
for which we have the standard inequality
(See Proposition 3.4 of \cite{s1}):
\begin{equation}
 \label{eq;06.1.20.151}
 \int_{Y-D}\tr\bigl(G(h_{HE},\DDlambda)^{\bot\,2}\bigr)\geq 0.
\end{equation}
Here $G(h_{HE},\DDlambda)^{\bot}$ denotes the trace free part
of $G(h_{HE},\DDlambda)$.
Hence we obtain the desired inequality from 
(\ref{eq;06.1.20.151}).
\hfill\qed

\subsection{Some formula
 and vanishing of characteristic numbers}

\subsubsection{Formula for
 $\int_X\parch_2(\vecE_{\ast})$}

Let $X$ be a smooth projective surface,
and let $D$ be a simple normal crossing divisor of $X$.
We will derive some formulas and vanishings
for the characteristic numbers
of $(\vecE_{\ast},\DDlambda)$.

\begin{rem}
\label{rem;06.1.20.50}
To begin with, we remark that 
we have only to show such formulas
for regular filtered $\lambda$-flat bundles
satisfying the following conditions
due to the method of perturbation of the parabolic structure
(Subsection {\rm\ref{subsubsection;06.1.19.10}}).
\begin{itemize}
\item {\em graded semisimple}.
\item $\Par\bigl(\vecE_{\ast},i\bigr)\subset\rnum$.
\item  $0\not\in \Par\bigl(\vecE_{\ast},i\bigr)$
 for any $i\in S$.
\end{itemize}
We will use it without mention
in the following argument.

We restrict ourselves to the case $\dim X=2$
just for simplicity.
The formula can be obviously generalized
for $\int_X\parch_{2,L}(\vecE_{\ast})$ of
regular $\lambda$-flat parabolic bundles 
$(\vecE_{\ast},\DDlambda)$ on $(X,D)$
in codimension two for $\dim X>2$,
where $L$ denotes a line bundle on $X$.
\hfill\qed
\end{rem}

For simplicity of the description,
we put as follows, for 
$u=\in
 \KMS(i):=\KMS\bigl(\prolong{E},i\bigr)$:
\[
 r(i,u):=\rank_{D_i}\bigl(
 \lefttop{i}\Gr^{F,\EE}_{u}(\prolong{E})
 \bigr)
\]
For any point $P\in D_i\cap D_j$
and $(u_i,u_j)\in\KMS(P):=\KMS\bigl(\prolong{E},P\bigr)$,
we put as follows:
\[
 r(P,u_i,u_j):=\rank\bigl(
 \lefttop{P}\Gr^{F,\EE}_{u_i,u_j}\bigl(E\bigr)
 \bigr)
\]

\begin{prop}
 \label{prop;06.1.13.70}
We have the following equality:
\begin{multline}
\label{eq;08.2.7.61}
 \int_X2\parch_2(\vecE_{\ast})=
\sum_{i\in S}\sum_{u\in\KMS(i)}
\bigl(
 \Re(\lambda^{-1}\alpha)+a
 \bigr)^2\cdot
 r(i,u)\cdot [D_i]^2\\
+\sum_{i\in S}
 \sum_{\substack{j\neq i\\ P\in D_i\cap D_j}}
 \sum_{(u_i,u_j)\in\KMS(P)}
 \bigl(\Re \lambda^{-1}\alpha_i+a_i\bigr)
 \bigl(\Re\lambda^{-1}\alpha_j+a_j\bigr)
\cdot r(P,u_i,u_j).
\end{multline}
We also have the following equalities:
\begin{equation} 
\label{eq;08.2.7.60}
 2\parch_2(\vecE_{\ast})=
\sum_{i\in S}\sum_{u\in\KMS(i)}
 \Re\bigl(\lambda^{-1}\alpha+a\bigr)\cdot
 \Bigl(
-\pardeg\bigl(\lefttop{i}\Gr^{F,\EE}_{a,\alpha}
 \bigl(\prolong{E}\bigr)_{\ast}\bigr)
+a\cdot r(i,u)\cdot [D_i]^2
 \Bigr).
\end{equation}
\begin{equation}
\label{eq;08.2.7.62}
0=
\sum_{i\in S}\sum_{u\in \KMS(i)}
 \Image(\lambda^{-1}\alpha)
\cdot
 \Bigl(
-\pardeg_{D_i}\bigl(
 \lefttop{i}\Gr^{F,\EE}_{u}(\prolong{E})_{\ast}
\bigr)
+a\cdot r(i,u)\cdot [D_i]^2
 \Bigr)
\end{equation}
\end{prop}
\pf
Let $X_{\delta}:=\bigcap\bigl\{|\sigma_i|\geq\delta\bigr\}$
and $Y_{\delta,i}:=X_{\delta}\cap \{|\sigma_i|=\delta\}$.
We have 
$R(h_0)=-\lambda^{-1}(1+|\lambda|^2)\cdot d''\theta$.
Hence, we have the following equality:
\begin{equation}
 \label{eq;08.2.7.50}
\lim_{\delta\to 0}
\left(\frac{\sqrt{-1}}{2\pi}\right)^2
 \int_{X_{\delta}}
 \tr\bigl(R(h_0)^2\bigr)
=-\frac{1+|\lambda|^2}{\lambda}
 \lim_{\delta\to 0}
\left(\frac{\sqrt{-1}}{2\pi}\right)^2
 \int_{\del X_{\delta}}
 d\tr\bigl(\theta\cdot R(h_0)\bigr)
\end{equation}
By using the estimates in Subsections 
\ref{subsection;08.2.7.10}--\ref{subsection;08.2.7.11},
the contribution of $Y_{\delta,i}$ to (\ref{eq;08.2.7.50})
is the following:
\begin{multline}
-\frac{1}{m}
 \sum_{(a,\alpha)\in\KMS(i)}
 m\big(a+\lambda^{-1}\alpha\bigr)
 \frac{\sqrt{-1}}{2\pi}
 \int_{D_i}\Bigl(
 \tr\bigl(R(\lefttop{i}h_{a,\alpha})\bigr)
-a\cdot r\bigl(i,(a,\alpha)\bigr)\cdot \tau_i
\Bigr) \\
=-\sum_{(a,\alpha)\in\KMS(i)}
 (a+\lambda^{-1}\alpha)
\cdot
 \bigl(
 \deg_{D_i}\bigl(
 \lefttop{i}\Gr^{F,\EE}_{a,\alpha}(\prolong{E})_{\ast}
 \bigr)
-a\cdot r\bigl(i,(a,\alpha)\bigr) \cdot [D_i]^2
 \bigr)
\end{multline}
By taking the real part,
we obtain (\ref{eq;08.2.7.60}).
By taking the imaginary part,
we obtain (\ref{eq;08.2.7.62}).
The equality (\ref{eq;08.2.7.61})
follows from (\ref{eq;08.2.7.60})
and Lemma \ref{lem;08.2.7.40}
by a formal calculation.
\hfill\qed

\begin{lem}
\label{lem;08.2.7.105}
For any $C^{\infty}$ $2$-form $\tau$,
we have the following:
\begin{equation}
\label{eq;08.2.7.120}
\int_X\parchern_1(\vecE_{\ast})\cdot\tau
=\frac{\sqrt{-1}}{2\pi}
 \int_{X}\tr R(h_0)\cdot \tau
=-\sum_{i\in S}
 \sum_{(a,\alpha)\in\KMS(i)}
 \Re(\lambda^{-1}\alpha+a)
\cdot r\bigl(i,(a,\alpha)\bigr)\cdot (D_i,\tau)
\end{equation}
Namely,
the cohomology class of $\tr R(h_0)$
is $\parchern_1(\vecE_{\ast})$.
In particular,
we also have the following equality:
\begin{equation}
\label{eq;08.2.7.102}
 \pardeg_{\omega}(\vecE_{\ast})
=-\sum_{i\in S}
 \sum_{(a,\alpha)\in\KMS(i)}
 \Re(\lambda^{-1}\alpha+a)
\cdot r\bigl(i,(a,\alpha)\bigr)\cdot (D_i,\omega)
\end{equation}
\end{lem}
\pf
Recall we have
$R(h_0)=
 \lambda^{-1}(1+|\lambda|^2)\cdot
 d''\theta_0$.
Then, we obtain (\ref{eq;08.2.7.120})
by using the estimates in Subsections 
\ref{subsection;08.2.7.10}--\ref{subsection;08.2.7.11}.
\hfill\qed

\begin{rem}
We considered the KMS-spectra of
$\prolong{E}$.
But, we have the following equality
for any $\vecc\in\real^S$ and $i\in S$:
\[
 \bigl\{
 \Re(\lambda^{-1}\alpha)+a\,\big|\,
 (a,\alpha)\in
 \KMS(\prolong{E},i)
 \bigr\}
=\bigl\{
 \Re(\lambda^{-1}\alpha)+a\,\big|\,
 (a,\alpha)\in
 \KMS(\prolongg{\vecc}{E},i)
 \bigr\}
\]
We also have such comparison
for $\KMS\bigl(\prolong{E},P\bigr)$
and $\KMS\bigl(\prolongg{\vecc}{E}\bigr)$
for $\vecc\in\real^S$
and $P\in D_i\cap D_j$.
Namely,
the choice $\prolong{E}$ is not essential.
(See also Section {\rm\ref{section;06.2.3.150}}.)
\hfill\qed
\end{rem}

\subsubsection{Remark on the vanishing 
of the parabolic Chern character numbers}

Recall the formulas for $\int_X\parch_2(\vecE_{\ast})$
in Proposition \ref{prop;06.1.13.70}
and the formula for
$\pardeg_{\omega}(\vecE_{\ast})$
in Lemma \ref{lem;08.2.7.105}.
Then, we immediately obtain the following corollary.
\begin{cor}
 \label{cor;06.1.13.120}
If $a+\Re\lambda^{-1}\alpha=0$ holds
for any KMS-spectrum $(a,\alpha)$ of 
$(\vecE_{\ast},\DDlambda)$,
the characteristic numbers
$\pardeg_{\omega}(\vecE_{\ast})$
and $\int_X\parch_2(\vecE_{\ast})$
automatically vanish.
\hfill\qed
\end{cor}

\begin{rem}
 \label{rem;06.2.3.1}
Let $E$ be a vector bundle on $X-D$
with a flat connection $\nabla$.
We have the Deligne extension
$(\widetilde{E},\nabla)$.
(See Subsection {\rm\ref{subsubsection;06.2.4.30}},
for example.)
We have 
the canonically defined parabolic structure $\vecF$
such that 
$\Re\alpha+a=0$ for any KMS-spectrum.
In that case,
the stability of $(\widetilde{E},\vecF,\nabla)$
and the semisimplicity of $(E,\nabla)$ is equivalent.
The corollary means
$\int_X\parchern_2(\widetilde{E},\vecF)
=\pardeg_{\omega}(\widetilde{E},\vecF)=0$.

When $(E,\nabla)$ is semisimple,
we know that there exists the Corlette-Jost-Zuo metric
of $(E,\nabla)$ 
which is a pure imaginary tame pluri-harmonic metric
adapted to the parabolic structure $\vecF$
(See {\rm \cite{corlette}} for the case $D=\emptyset$
and {\rm\cite{JZ2}} for the general case.
See also {\rm\cite{mochi2}}.)
To show such an existence theorem
from the Kobayashi-Hitchin correspondence,
we have to show the vanishing of the characteristic numbers
which is ``the obstruction 
on the way from harmonicity to pluri-harmonicity''.
Corollary {\rm\ref{cor;06.1.13.120}} 
clarifies the point.
\hfill\qed
\end{rem}

%% file: 4.tex
\subsection{Statements}

In this section,
we will show continuity of two kinds of families
of harmonic metrics on curves,
i.e.,
Proposition \ref{prop;06.1.21.1}
and Proposition \ref{prop;06.1.18.80}.
We will give a detailed proof of the first one.
Because the second one can be proved similarly and more easily,
we just give some remarks in the end of this section.

\subsubsection{Continuity for $\epsilon$-perturbation}

Let $C$ be a smooth projective curve
with a simple divisor $D$.
Let $(E,\vecF,\DDlambda)$ be a parabolic flat
$\lambda$-connection over $(C,D)$,
which is stable and $\pardeg(E,\vecF)=0$.
Let $\vecF^{(\epsilon)}$ be the $\epsilon$-perturbation
of the parabolic structures,
explained in (II) of 
Subsection \ref{subsubsection;06.1.19.10}.
We remark $\det(E,\vecF)=\det(E,\vecF^{(\epsilon)})$.
Let $h^{(\epsilon)}$ be the harmonic metric for 
$(E,\vecF^{(\epsilon)},\DDlambda)$.
Let $\theta^{(\epsilon)}$ denote the Higgs fields
for the harmonic bundles $(E,\DDlambda,h^{(\epsilon)})$.

\begin{prop}
 \label{prop;06.1.21.1}
The sequences $\{h^{(\epsilon)}\,|\,\epsilon>0\}$
and $\bigl\{\theta^{(\epsilon)}\bigr\}$ converge
to $h^{(0)}$ and $\theta^{(0)}$
respectively, in the $C^{\infty}$-sense locally on $C-D$.
\end{prop}

The proof is given in Subsection
\ref{subsection;06.1.21.3}
after the preparation given in Subsections
\ref{subsection;06.1.21.2}--\ref{subsection;06.1.21.30}.

\subsubsection{Continuity for a holomorphic family}

Before going into the proof of Proposition \ref{prop;06.1.21.1},
we give a similar statement for another family.
Let $\nbigc\lrarr \Delta$ be a holomorphic family
of smooth projective curve,
and $\nbigd\lrarr\Delta$ be a relative divisor.
Let $(E,\vecF,\DDlambda)$ be a parabolic flat bundle
on $(\nbigc,\nbigd)$.
Let $t$ be any point of $\Delta$.
We denote the fibers by $\nbigc_t$ and $\nbigd_t$,
and the restriction of $(E,\vecF,\DDlambda)$
to $(\nbigc_t,\nbigd_t)$ is denoted by
$(E_t,\vecF_t,\DDlambda_t)$.
We assume 
$\pardeg(E_t,\vecF_t)=0$ 
and that $(E_t,\vecF_t,\DDlambda_t)$ is stable
for each $t$.
For simplicity, we also assume that
we are given a pluri harmonic metric $h_{\det (E)}$
of $\det(E,\DDlambda)_{|\nbigc-\nbigd}$
which is adapted to the induced parabolic structure.

Let $h_{H,t}$ be a harmonic metric
of $(E_t,\vecF_t,\DDlambda_t)$
such that
$\det(h_{H,t})=h_{\det(E)\,|\,\nbigc_t}$.
They give the metric $h_H$ of $E$.
Let $\theta_{H,t}$ be the Higgs field
obtained from
$(E_t,\DDlambda,h_{H,t})$,
which is a section of
$\End(E_t)\otimes\Omega^{1,0}_{\nbigc_t}(\log\nbigd_t)$.
They give the section $\theta_H$
of $\End(E)\otimes\Omega^{1,0}_{\nbigc/\Delta}(\log\nbigd)$,
where $\Omega^{1,0}_{\nbigc/\Delta}(\log \nbigd)$ denotes
the sheaf of the logarithmic relative $(1,0)$-forms.

\begin{prop}
 \label{prop;06.1.18.80}
$h_H$ and $\theta_H$ are continuous.
Their derivatives of any degree along the fiber directions
are continuous.
\end{prop}

Since Proposition \ref{prop;06.1.18.80}
can be proved similarly and more easily,
we will not give a detailed proof.
See Remark \ref{rem;06.1.21.10}.

\subsection{Preliminary from elementary calculus}
\label{subsection;06.1.21.2}

For any
$z\in\Delta^{\ast}=
 \bigl\{z\in\cnum\,\big|\,|z|<1\bigr\}$
and $\epsilon>0$,
we put as follows:
\[
 L_{\epsilon}(z):=\frac{|z|^{-\epsilon}-|z|^{\epsilon}}{\epsilon},
\quad
 K_{\epsilon}(z):=\frac{|z|^{-\epsilon}+|z|^{\epsilon}}{2},
\quad
 M_{\epsilon}(z):=|z|^{4\epsilon}(1-\log|z|^{4\epsilon}).
\]
We also put
$L_0(z):=-\log|z|^2$, $K_0(z)=1$ and $M_0(z)=1$.
Then, they are continuous with respect to 
$(z,\epsilon)\in\Delta^{\ast}\times\real_{\geq\,0}$.

\begin{lem} 
\label{lem;06.1.14.10}
For any $(z,\epsilon)\in \Delta^{\ast}\times\real_{\geq\,0}$,
we have $L_0(z)\leq L_{\epsilon}(z)$.
\end{lem}
\pf
We put $g(\epsilon):=
 a^{-\epsilon}-a^{\epsilon}+\epsilon\cdot \log a^2$
for $0<a<1$ and $0\leq\epsilon$.
Taking the derivative with respect to $\epsilon$,
we obtain the following:
\[
 g'(\epsilon)=-\bigl(a^{-\epsilon}+a^{\epsilon}\bigr)\log a
 +\log a^2,
\quad
g''(\epsilon)=(a^{-\epsilon}-a^{\epsilon})(\log a)^2\geq 0.
\]
Since we have $g(0)=g'(0)=0$,
the claim of the lemme follows.
\hfill\qed

\begin{lem}
\label{lem;06.1.2.7}
$(K_{\epsilon}(z)-1)\cdot 
\bigl(
L_{\epsilon}(z)^2\cdot\epsilon^2\cdot |z|^{\epsilon}
\bigr)^{-1}$ are bounded on $\Delta^{\ast}$,
independently of $\epsilon$.
We also have $K_{\epsilon}(z)-1\geq 0$.
\end{lem}
\pf
The second claim is clear.
Let us check the first claim.
We put as follows, for 
$0<a<1$ and $0\leq \epsilon\leq 1$:
\[
 g_1(\epsilon):=
a^{-\epsilon}-2+a^{\epsilon},
\quad
g_{2}(\epsilon):=
 (a^{-\epsilon}-a^{\epsilon})^2a^{\epsilon}
=a^{-\epsilon}-2a^{\epsilon}+a^{3\epsilon}.
\]
We have only to show that
$g_2(\epsilon)\geq g_1(\epsilon)$.
We put
$g(\epsilon):=g_2(\epsilon)-g_1(\epsilon)
=2+a^{3\epsilon}-3a^{\epsilon}$.
By taking the derivative with respect to $\epsilon$,
we obtain the following:
\[
 g'(\epsilon)
=3a^{3\epsilon}\cdot\log a-3a^{\epsilon}\cdot\log a
=3(-a^{3\epsilon}+a^{\epsilon})(-\log a)\geq 0.
\]
Since we have $g(0)=0$,
we obtain $g(\epsilon)\geq 0$.
Thus we are done.
\hfill\qed

\begin{lem}
 \label{lem;06.1.2.6}
$\bigl(1-M_{\epsilon}(z)\bigr)\cdot 
\bigl(
L_{\epsilon}(z)^2\cdot \epsilon^2\cdot |z|^{\epsilon}
\bigr)^{-1}$ are bounded on $\Delta^{\ast}$,
independently of $\epsilon$.
We also have $1-M_{\epsilon}(z)\geq  0$.
\end{lem}
\pf
We have only to show the following inequalities
for $0<a<1$ and $0\leq \epsilon<1$:
\[
 0\leq 1-a^{4\epsilon}(1-\log a^{4\epsilon})
\leq 3(a^{-\epsilon}-a^{\epsilon})^2 a^{\epsilon}.
\]
To show the left inequality,
we put
$h(\epsilon):=1-a^{4\epsilon}(1-\log a^{4\epsilon})$.
By taking the derivative with respect to $\epsilon$,
we have
$h'(\epsilon)=-a^{4\epsilon}\log a^4 (1-\log a^{4\epsilon})
+a^{4\epsilon}\log a^4
=\epsilon a^{4\epsilon}(\log a^4)^2\geq 0$.
We also have $h(0)=0$.
Hence, we obtain $h(\epsilon)\geq 0$.
To show the right inequality,
we put as follows:
\[
 g(\epsilon):=
a^{-4\epsilon}
\Bigl(
 3(a^{-\epsilon}-a^{\epsilon})^2 a^{\epsilon}
-\bigl(1-a^{4\epsilon}(1-\log a^{4\epsilon})\bigr)
 \Bigr)
=3(a^{-5\epsilon}-2a^{-3\epsilon}+a^{-\epsilon})
+(1-\log a^{4\epsilon})-a^{-4\epsilon}.
\]
By taking the derivative with respect to $\epsilon$,
we obtain the following:
\[
 g'(\epsilon)=3\bigl(a^{-5\epsilon}(-5\log a)
 -2a^{-3\epsilon}(-3\log a)
 +a^{-\epsilon}(-\log a)
 \bigr)
-4\log a-a^{-4\epsilon}(-4\log a)
\]
\[
 g''(\epsilon)
=\bigl(
75 a^{-5\epsilon}
-16 a^{-4\epsilon}
-54a^{-3\epsilon}
+3a^{-\epsilon}
 \bigr)\cdot(\log a)^2.
\]
It is easy to check $g''(\epsilon)\geq 0$
by using $a^{-5\epsilon}\geq a^{-k\epsilon}$ $(k=3,4)$.
Since we have $g'(0)=g(0)=0$,
we obtain $g(\epsilon)\geq 0$.
Thus we are done.
\hfill\qed

\begin{lem}
 \label{lem;06.1.14.11}
Let $P(t)$ be a polynomial with variable $t$,
and let $b$ be any fixed positive number.
Then, we have the boundedness of
$|z|^{b\epsilon}\cdot P\bigl(\epsilon L_0(z)\bigr)$
on $\Delta^{\ast}$,
independently of $0\leq \epsilon\leq 1/2$.
\end{lem}
\pf
We put $u:=|z|^{\epsilon}$,
and then $|z|^{b\epsilon}P(\epsilon L_0(z))
=u^b\cdot P\bigl(L_0(u)\bigr)$.
Hence, we have only to show the boundedness
of $u^b\cdot P\bigl(L_0(u)\bigr)$
when $0<u<1$,
but it is easy.
\hfill\qed

\subsection{A family of the metrics
 for a logarithmic $\lambda$-flat bundle of rank two
 on a disc}

\subsubsection{Construction of a family of metrics}

We put $X=\Delta=\{z\,\big|\,|z|<1\}$.
Let $O$ denote the origin,
and we put $X^{\ast}:=X-\{O\}$.
We use the Kahler form
$\omega_{\epsilon}:=
 (\epsilon^2|z|^{\epsilon-2}+1)\cdot dz\cdot d\zbar$
in this subsection.
We will use the notation
in Subsection \ref{subsection;06.1.21.2}.

To begin with, we recall an example 
of a harmonic bundle on a punctured disc.
Let $E=
 \nbigo_{X}\cdot v_1\oplus\nbigo_{X}\cdot v_2$
be a holomorphic vector bundle on a disc.
Let $\theta$ be a Higgs bundle such that
$\theta \cdot v_1=v_2\cdot dz/z$
and $\theta \cdot v_2=0$.
Let $h$ be the metric of $E_{|X^{\ast}}$
such that
$ h(v_1,v_1)=L_0$,
$h(v_2,v_2)=L_0^{-1}$
and $h(v_i,v_j)=0$ $(i\neq j)$.
Recall that the tuple
$(E,\delbar_E,\theta,h)$ is a harmonic bundle.
Let us see the associated $\lambda$-connection.
We put $u_1:=v_1$ and
$u_2:=v_2-\lambda\cdot L_0^{-1}\cdot v_1$.
Then, we can show
$(\delbar_E+\lambda\theta^{\dagger})u_i=0$
$(i=1,2)$,
$\DDlambda u_1=u_2\cdot dz/z$
and $\DDlambda u_2=0$
by a direct calculation.
We also have the following:
\[
 h(u_1,u_1)=L_0,\quad
 h(u_2,u_2)=(1+|\lambda|^2)\cdot L_0^{-1},\quad
 h(u_1,u_2)=-\overline{\lambda},
\quad 
 h(u_2,u_1)=-\lambda.
\]

Motivated by this example,
we consider the following family of the metrics
$h_{\epsilon}$ on the $\lambda$-connection $(E,\DDlambda)$
given as follows:
\[
 h_{\epsilon}(u_1,u_1)=L_{\epsilon},
\quad
h_{\epsilon}(u_2,u_2)
 =\bigl(1+|\lambda|^2\bigr)\cdot L_{\epsilon}^{-1},
\quad
h_{\epsilon}(u_1,u_2)=-\lambdabar\cdot M_{\epsilon},
\quad
h_{\epsilon}(u_2,u_1)=-\lambda\cdot M_{\epsilon}.
\]
The $\lambda$-connection $\DDlambda$
and the metric $h_{\epsilon}$ induce
the operators
$\delbar_{\epsilon}$ and $\theta_{\epsilon}$
(Subsection \ref{subsubsection;06.1.10.7}).
The main purpose of this subsection
is to show the following proposition.

\begin{prop}
 \label{prop;06.1.14.2}
There exists a some positive constant $C$
such that
$\bigl|
\delbar_{\epsilon}\theta_{\epsilon}
\bigr|_{h_{\epsilon},\omega_{\epsilon}}
\leq C$
for any $0\leq \epsilon<1/2$.
\end{prop}
Although the proof of the proposition is just a calculation,
we will give the detail in the rest of this subsection.

\begin{rem}
 \label{rem;06.1.14.3}
Let $h_{\epsilon}'$ be the metric
determined by
$h'_{\epsilon}(u_1,u_1)=L_{\epsilon}$,
$h'_{\epsilon}(u_2,u_2)=L_{\epsilon}^{-1}$
and $h'_{\epsilon}(u_i,u_j)=0$ $(i\neq j)$.
Then, there exist positive constants $C_i$ $(i=1,2)$
such that
$C_1\cdot h_{\epsilon}'\leq h_{\epsilon}
 \leq C_2\cdot h_{\epsilon}'$
for any $0\leq\epsilon\leq 1/2$.
Hence, we have only to consider the norms
for $h_{\epsilon}'$ instead of those for $h_{\epsilon}$.
\hfill\qed
\end{rem}

\subsubsection{Preliminary}
Let $H_{\epsilon}$ be the hermitian matrix valued function
given by $H_{\epsilon}:=H(h_{\epsilon},\vecu)$,
i.e.,
\[
 H_{\epsilon}:=
 \left(
 \begin{array}{cc}
 L_{\epsilon} & -\overline{\lambda}\cdot M_{\epsilon}\\
 -\lambda\cdot M_{\epsilon} &
 (1+|\lambda|^2)L_{\epsilon}^{-1}
 \end{array}
 \right).
\]
Let $N$ be determined by
$\DDlambda \vecu=\vecu\cdot N\cdot dz/z$,
and 
let $N^{\dagger}_{\epsilon}$ denote the adjoint of $N$
with respect to the metric $H_{\epsilon}$,
i.e.,
\[
 N=\left(
 \begin{array}{cc}
 0 & 0 \\
 1 & 0
 \end{array}
\right),
\quad
 N^{\dagger}_{\epsilon}=
\overline{H}_{\epsilon}^{-1}
\cdot \lefttop{t}\overline{N}\cdot
\overline{H}_{\epsilon}
=\frac{1}{1+|\lambda|^2(1-M_{\epsilon}^2)}
 \left(\begin{array}{cc}
 -\overline{\lambda}(1+|\lambda|^2)L^{-1}_{\epsilon}M_{\epsilon}
 & (1+|\lambda|^2)^2L_{\epsilon}^{-2}\\
 -\lambdabar^2M_{\epsilon}^2
 & \lambdabar (1+|\lambda|^2)M_{\epsilon}L_{\epsilon}^{-1}
 \end{array}
 \right).
\]
Recall the calculation given in
Subsection \ref{subsubsection;06.1.14.1}.
Then, $\delbar_{\epsilon}$ and $\theta_{\epsilon}$
can be described with respect to $\vecu$
as follows:
\[
 \delbar_{\epsilon}\vecu=
\vecu\cdot \frac{\lambda}{1+|\lambda|^2}
\Bigl(\lambdabar \cdot\Hbarepsilon^{-1}\delbar\Hbarepsilon
-N^{\dagger}_{\epsilon}\frac{d\zbar}{\zbar}
 \Bigr),
\quad
 \theta_{\epsilon}\vecu=
 \vecu\frac{1}{1+|\lambda|^2}
 \Bigl(
 N\frac{dz}{z}-\lambda\Hbarepsilon^{-1}\del \Hbarepsilon
 \Bigr).
\]
Therefore, $\delbar_{\epsilon}(\theta_{\epsilon})$
is described by the following 
$2\times 2$-matrix valued 2-form
with respect to $\vecu$:
\begin{equation}
 \label{eq;06.1.2.2}
 \frac{1}{1+|\lambda|^2}
 \delbar\Bigl(
 -\lambda\Hbarepsilon^{-1}\del\Hbarepsilon
 \Bigr)
+\frac{\lambda}{(1+|\lambda|^2)^2}
\left(
 \Bigl[
 \lambdabar\cdot \Hbarepsilon^{-1}\delbar\Hbarepsilon,
 N\frac{dz}{z}
 \Bigr]
-\Bigl[
 N_{\epsilon}^{\dagger}\frac{d\zbar}{\zbar},\,
 N\frac{dz}{z}
 \Bigr]
+\Bigl[
 N_{\epsilon}^{\dagger}\frac{d\zbar}{\zbar},\,\,
 \lambda \Hbarepsilon^{-1}\del\Hbarepsilon
 \Bigr]
\right).
\end{equation}
Here we have used 
$\bigl[\Hbarepsilon^{-1}\del \Hbarepsilon,\,\,
 \Hbarepsilon^{-1}\delbar\Hbarepsilon
 \bigr]=0$,
which can be checked easily.
\begin{lem}
 \label{lem;06.1.14.150}
To show Proposition {\rm\ref{prop;06.1.14.2}},
we have only to show the uniform boundedness of
$(1,1)$-entry, $(2,2)$-entry,
$L_{\epsilon}\times$ $(1,2)$-entry
and $L_{\epsilon}^{-1}\times$ $(2,1)$-entry,
in the matrix valued function {\rm(\ref{eq;06.1.2.2})}.
\end{lem}
\pf
It follows from  Remark \ref{rem;06.1.14.3}.
\hfill\qed

\vspace{.1in}

In the following calculation,
we often use the notation $L$ and $M$
instead of $L_{\epsilon}$ and $M_{\epsilon}$,
if there are no risk of confusion.
Let us see $\Hbarepsilon^{-1}\del \Hbarepsilon$.
We have the following equality:
\[
 \overline{H}_{\epsilon}^{-1}
=\frac{1}{1+|\lambda|^2(1-M_{\epsilon}^2)}
 \left(
 \begin{array}{cc}
 (1+|\lambda|^2)\cdot L_{\epsilon}^{-1} &
 \lambda\cdot M_{\epsilon} \\
 \overline{\lambda}\cdot M_{\epsilon}
 & L_{\epsilon}
 \end{array}
 \right),
\quad
 \del \overline{H}_{\epsilon}
=\left(
 \begin{array}{cc}
 \del L_{\epsilon} & -\lambda\cdot \del M_{\epsilon}\\
 -\overline{\lambda}\cdot \del M_{\epsilon} &
 (1+|\lambda|^2)\cdot \del L_{\epsilon}^{-1}
 \end{array}
 \right).
\]
Then, we obtain the following formula
for $\Hbarepsilon^{-1}\del\Hbarepsilon$:
\begin{equation}
 \label{eq;06.1.14.500}
 \overline{H}_{\epsilon}^{-1}
 \del \overline{H}_{\epsilon}
=\frac{1}{1+|\lambda|^2(1-M_{\epsilon}^2)}
 \left(
 \begin{array}{cc}
 (1+|\lambda|^2)L^{-1}\del L
-|\lambda|^2 M\del M &
 \lambda(1+|\lambda|^2)\bigl(
 -L^{-1}\del M+M\del L^{-1}
 \bigr) \\
 \overline{\lambda}(M\del L-L\del M) &
 (1+|\lambda|^2)L\del L^{-1}
-|\lambda|^2M\cdot\del M
 \end{array}
 \right).
\end{equation}
We also have a similar formula for
$\Hbarepsilon^{-1}\delbar \Hbarepsilon$.
We obtain the following formula for
$\delbar \bigl(\Hbarepsilon^{-1}\del\Hbarepsilon\bigr)$:
\begin{multline}
 \label{eq;06.1.2.1}
 \delbar\bigl(
 \overline{H}_{\epsilon}^{-1}
 \del \overline{H}_{\epsilon}
 \bigr)
=\frac{2|\lambda|^2M\delbar M }
 {\bigl(1+|\lambda|^2(1-M^2)\bigr)^2}
 \Hbar_{\epsilon}^{-1}\del \Hbar_{\epsilon} \\
+\frac{1}{1+|\lambda|^2(1-M^2)}
 \left(
 \begin{array}{cc}
 (1+|\lambda|^2)\delbar\del \log L
-2^{-1}|\lambda|^2\delbar\del M^2 &
 \lambda(1+|\lambda|^2)
 (M\delbar\del L^{-1}-L^{-1}\delbar\del M)\\
 \overline{\lambda}
 (M\delbar\del L-L\delbar\del M) &
 (1+|\lambda|^2)\delbar\del\log L^{-1}
-2^{-1}|\lambda|^2\delbar\del M^2
 \end{array}
 \right).
\end{multline}
The commutator of
$\Hbarepsilon^{-1}\delbar\Hbarepsilon$
and $N\cdot dz/z$ is as follows:
\begin{equation}
 \label{eq;06.1.14.200}
 \Bigl[
 \overline{H}_{\epsilon}^{-1}\delbar \overline{H}_{\epsilon},
 N\cdot \frac{dz}{z}
 \Bigr]
=\frac{(1+|\lambda|^2)}{1+|\lambda|^2(1-M^2)}
\left(
 \begin{array}{cc}
  \lambda
 (-L^{-1}\delbar M+M\delbar L^{-1}) & 0 \\
 2L\delbar L^{-1} &
 -\lambda
 (-L^{-1}\delbar M+M\delbar L^{-1}) 
 \end{array}
\right)\frac{dz}{z}.
\end{equation}
Let us see the commutator of
$\Hbarepsilon^{-1}\del \Hbarepsilon$ and $N_{\epsilon}^{\dagger}$.
By direct calculations,
we have the following equality:
\begin{multline}
 \Hbarepsilon^{-1}\del\Hbarepsilon\cdot
 N_{\epsilon}^{\dagger}
=\frac{1}{1+|\lambda|^2(1-M^2)}
 \left(
\begin{array}{cc}
 -\lambdabar (1+|\lambda|^2)L^{-2}M\del L &
(1+|\lambda|^2)^2L^{-3}\del L\\
\lambdabar^2\cdot M\del M  & 
-\lambdabar (1+|\lambda|^2)L^{-1}\del M
\end{array}
 \right) \\
+\frac{1}{\bigl(1+|\lambda|^2(1-M^2)\bigr)^2}
\left(
 \begin{array}{cc}
 2|\lambda|^2\lambdabar(1+|\lambda|^2)M^2L^{-1}\del M
 &
 -2|\lambda|^2(1+|\lambda|^2)^2ML^{-2}\del M\\
 2M^3\del M\lambdabar^2|\lambda|^2 &
 -2\lambdabar|\lambda|^2(1+|\lambda|^2)M^2L^{-1}\del M
 \end{array}
\right).
\end{multline}
We also have the following:
\begin{equation}
 N_{\epsilon}^{\dagger}
\cdot
 \Hbarepsilon^{-1}\del\Hbarepsilon
=\frac{1}{1+|\lambda|^2(1-M^2)}
 \left(
 \begin{array}{cc}
 -\lambdabar (1+|\lambda|^2)L^{-1}\del M &
 (1+|\lambda|^2)^2L^{-1}\del L^{-1}\\
 -\lambdabar^2M\del M &
 \lambdabar (1+|\lambda|^2)M\del L^{-1}
 \end{array}
 \right).
\end{equation}
Therefore, we obtain the following formula:
\begin{multline}
 \label{eq;06.1.14.300}
 \Bigl[
 N_{\epsilon}^{\dagger}\frac{d\zbar}{\zbar},
 \Hbarepsilon^{-1}\del \Hbarepsilon
 \Bigr] \\
=\frac{1}{1+|\lambda|^2(1-M^2)}
 \frac{d\zbar}{\zbar}
 \left(
 \begin{array}{cc}
 -\lambdabar(1+|\lambda|^2)(L^{-1}\del M-L^{-2}M\del L)
 &
 -2(1+|\lambda|^2)^2L^{-3}\del L\\
 -2\lambdabar^2M\del M &
 \lambdabar(1+|\lambda|^2)
 (M\del L^{-1}+L^{-1}\del M)
\end{array}
 \right) \\
- \frac{2|\lambda|^2}{\bigl( 1+|\lambda|^2(1-M^2)\bigr)^2}
 \frac{d\zbar}{\zbar}
 \left(
 \begin{array}{cc}
 \lambdabar(1+|\lambda|^2)M^2L^{-1}\del M &
 -(1+|\lambda|^2)^2ML^{-2}\del M \\
 \lambdabar^2M^3\del M &
 -\lambdabar (1+|\lambda|^2)M^2L^{-1}\del M
 \end{array}
 \right) 
\end{multline}
The commutator of $N$ and $N^{\dagger}_{\epsilon}$ is as follows:
\begin{equation}
 \label{eq;06.1.14.350}
 \bigl[N^{\dagger}_{\epsilon},N
 \bigr]
=\frac{1}{1+|\lambda|^2(1-M^2)}
 \left(
 \begin{array}{cc}
 (1+|\lambda|^2)^2L^{-2} & 0 \\
 2\lambdabar (1+|\lambda|^2)M L^{-1} &
 -(1+|\lambda|^2)^2L^{-2}
 \end{array}
 \right).
\end{equation}

\subsubsection{Estimate}
\label{subsubsection;06.1.14.501}

We have the following:
\begin{equation}
 \label{eq;06.1.2.3}
 \del L_{\epsilon}=-K_{\epsilon}\frac{dz}{z},
\quad
 \del K_{\epsilon}=-\frac{\epsilon^2}{4}L_{\epsilon}
 \frac{dz}{z},
\quad
 \del M_{\epsilon}=
 4\epsilon^2\cdot|z|^{4\epsilon}\cdot
 L_0\cdot \frac{dz}{z}.
\end{equation}
In particular, we have the following estimate:
\[
 M_{\epsilon}\del M_{\epsilon}
=O\Bigl(
\epsilon^2\cdot |z|^{8\epsilon}\cdot 
L_0\cdot \bigl(1+\epsilon L_0\bigr)
 \frac{dz}{z}
 \Bigr).
\]

Let us see the first term in the right hand side of
(\ref{eq;06.1.2.1}):
\begin{equation} 
 \label{eq;06.1.2.4}
\frac{2|\lambda|^2M_{\epsilon}\delbar M_{\epsilon} }
 {\bigl(1+|\lambda|^2(1-M_{\epsilon}^2)\bigr)^2}
 \Hbar_{\epsilon}^{-1}\del \Hbar_{\epsilon}
\end{equation}
For the $(1,1)$-entry and $(2,2)$-entry,
we have the following estimates:
\[
 M_{\epsilon}\delbar M_{\epsilon}\cdot
 L_{\epsilon}^{-1}\del L_{\epsilon}
=O\left(
 \epsilon^2\cdot L_0\cdot |z|^{8\epsilon}(1+\epsilon L_0)
 \frac{K_{\epsilon}}{L_{\epsilon}}
 \right)\frac{d\zbar \cdot dz}{|z|^2}
=O\left(
 |z|^{5\epsilon}(1+\epsilon L_0)
 \frac{L_0}{L_{\epsilon}}
 \right)\cdot \omega_{\epsilon}
\]
\[
 M_{\epsilon}\delbar M_{\epsilon}\cdot
 M_{\epsilon}\del M_{\epsilon}
=O\Bigl(
 \epsilon^4\cdot |z|^{16\epsilon}
 \cdot (1+\epsilon L_0)^2L_0^2
 \Bigr)\frac{dz\cdot d\zbar}{|z|^2}
=O\Bigl(
 |z|^{15\epsilon} \cdot (1+\epsilon L_0)^2(\epsilon L_0)^2
 \Bigr)
\cdot\omega_{\epsilon}.
\]
They are bounded with respect to $\omega_{\epsilon}$
due to Lemma \ref{lem;06.1.14.10}
and Lemma \ref{lem;06.1.14.11}.
Hence, the $(1,1)$-entry and the $(2,2)$-entry 
of (\ref{eq;06.1.2.4}) are bounded
independently of $\epsilon$.
Let us see the $(1,2)$-entry.
Recall Lemma \ref{lem;06.1.14.150}.
Hence, we have only to see the following:
\[
 L_{\epsilon}\times
(M_{\epsilon}\delbar M_{\epsilon})
 \cdot \bigl(
 L_{\epsilon}^{-1}\del M_{\epsilon}
-M_{\epsilon}\del L_{\epsilon}^{-1}
\bigr)
=M_{\epsilon}\delbar M_{\epsilon}\del M_{\epsilon}
+M_{\epsilon}^2\delbar M_{\epsilon}
 L_{\epsilon}^{-1}\del L_{\epsilon}.
\]
Both terms in the right hand side can be
estimated as in the previous case,
by using Lemma \ref{lem;06.1.14.10}
and Lemma \ref{lem;06.1.14.11}:
\[
 M_{\epsilon}\delbar M_{\epsilon}\del M_{\epsilon}
=O\Bigl(
 |z|^{10\epsilon}(1+\epsilon L_0)(\epsilon L_0)^2
 \Bigr)\cdot\omega_{\epsilon}
=O(1)\cdot \omega_{\epsilon}
\]
\[
  M_{\epsilon}^2\delbar M_{\epsilon}
 L_{\epsilon}^{-1}\del L_{\epsilon}
=O\Bigl(
 |z|^{11\epsilon}
 (1+\epsilon L_0)^2\frac{L_0}{L_{\epsilon}}
 \Bigr)
 \cdot\omega_{\epsilon}
=O(1)\cdot
 \omega_{\epsilon}
\]
The $(2,1)$-entry can be estimated similarly:
\[
 L_{\epsilon}^{-1}\times
 (M_{\epsilon}\delbar M_{\epsilon})
 \bigl(
 M_{\epsilon}\del L_{\epsilon}-L_{\epsilon}\del M_{\epsilon}
 \bigr)
=M_{\epsilon}^2L_{\epsilon}^{-1}
 \delbar M_{\epsilon}\del L_{\epsilon}
 -M_{\epsilon}\cdot\delbar M_{\epsilon}\del M_{\epsilon}
=O(1)
 \cdot\omega_{\epsilon}.
\]

Let us see the second term in the right hand side
of (\ref{eq;06.1.2.1}):
\begin{equation}
\label{eq:06.1.2.5}
 \frac{1}{1+|\lambda|^2(1-M^2)}
 \left(
 \begin{array}{cc}
 (1+|\lambda|^2)\delbar\del \log L
-2^{-1}|\lambda|^2\delbar\del M^2 &
 \lambda(1+|\lambda|^2)
 (M\delbar\del L^{-1}-L^{-1}\delbar\del M)\\
 \overline{\lambda}(1+|\lambda|^2)
 (M\delbar\del L-L\delbar\del M) &
 (1+|\lambda|^2)\delbar\del\log L^{-1}
-2^{-1}|\lambda|^2\delbar\del M^2
 \end{array}
 \right).
\end{equation}
It is easy to see the following estimate:
\begin{equation}
\delbar\del M_{\epsilon}^2=
O\Bigl(\epsilon^2\cdot |z|^{6\epsilon}(1+\epsilon L_0)^2\Bigr)
\cdot\omega_{\epsilon}
=O(\epsilon^2)\cdot \omega_{\epsilon}.
\end{equation}
Hence, it is bounded with respect to $\omega_{\epsilon}$
independently of $\epsilon$.
We remark that
$L_{\epsilon}^{-1}M_{\epsilon}\delbar\del L_{\epsilon}$
is also bounded
independently of $\epsilon$:
\[
 L_{\epsilon}^{-1}M_{\epsilon}\cdot\delbar\del L_{\epsilon}=
\frac{\epsilon^2}{4}M_{\epsilon}
 \cdot \frac{d\zbar\cdot dz}{|z|^2}
=O(1)\cdot\omega_{\epsilon}.
\]
Hence, we have the following,
modulo the uniformly bounded term
with respect to $(h_{\epsilon},\omega_{\epsilon})$:
\begin{equation}
 \label{eq;06.1.14.370}
 \delbar\bigl(\Hbarepsilon^{-1}\del\Hbarepsilon\bigr)
\equiv\frac{(1+|\lambda|^2)}{1+|\lambda|^2(1-M_{\epsilon}^2)}
 \left(
 \begin{array}{cc}
 \delbar\del\log L_{\epsilon} & 
 \lambda M_{\epsilon}\delbar\del L_{\epsilon}^{-1}\\
 0 & -\delbar\del \log L_{\epsilon}
 \end{array}
 \right).
\end{equation}

Let us see (\ref{eq;06.1.14.200}).
By the same argument,
we have the following uniform boundedness:
\[
 L_{\epsilon}^{-1}\delbar M_{\epsilon}
\cdot\frac{dz}{z}
=O\left(
 \epsilon^2|z|^{4\epsilon}\frac{L_0}{L_{\epsilon}}
 \right)\cdot\frac{dz\cdot d\zbar}{|z|^2}
=O(1)\cdot\omega_{\epsilon}.
\]
Hence, we have the following,
modulo the uniformly bounded terms
with respect to $(h_{\epsilon},\omega_{\epsilon})$:
\begin{equation}
 \label{eq;06.1.14.380}
\left[\Hbarepsilon^{-1}\delbar \Hbarepsilon,\,\,
 N\cdot \frac{dz}{z}\right]
\equiv
 \frac{(1+|\lambda|^2)}{1+|\lambda|^2(1-M_{\epsilon}^2)}
 \left(
 \begin{array}{cc}
 \lambda M_{\epsilon}\delbar L_{\epsilon}^{-1}
 & 0 \\
 2L_{\epsilon}\delbar L_{\epsilon}^{-1}
 & -\lambda
 M_{\epsilon}\delbar L_{\epsilon}^{-1}
 \end{array}
 \right)\cdot \frac{dz}{z}.
\end{equation}
Let us see (\ref{eq;06.1.14.300}).
We remark the following,
for any $k\geq 1$:
\[
 \frac{d\zbar}{\zbar}
 \frac{M_{\epsilon}^k\del M_{\epsilon}}{L_{\epsilon}}
=O\left(
 \epsilon^2
 |z|^{4(k+1)\epsilon}(1+\epsilon L_0)^{k}
 \frac{L_0}{L_{\epsilon}}
 \right)\cdot \frac{d\zbar\cdot dz}{|z|^2}
=O(1)\cdot \omega_{\epsilon}.
\]
Hence, the terms containing $\del M$ in the right hand side of
(\ref{eq;06.1.14.300}) can be ignored.
Hence, we obtain the following,
modulo the uniformly bounded terms
with respect to $(h_{\epsilon},\omega_{\epsilon})$:
\begin{equation}
 \label{eq;06.1.14.360}
 \left[
 N_{\epsilon}^{\dagger}\frac{d\zbar}{\zbar},\,\,
 \Hbarepsilon^{-1}\del \Hbarepsilon
 \right]
\equiv
 \frac{(1+|\lambda|^2)}{1+|\lambda|^2(1-M_{\epsilon}^2)}
 \frac{d\zbar}{\zbar}
 \left(
 \begin{array}{cc}
 \lambdabar 
 L_{\epsilon}^{-2}M_{\epsilon}\del L_{\epsilon}
 &
 -2(1+|\lambda|^2)L_{\epsilon}^{-3}\del L_{\epsilon}\\
 0 &
 \lambdabar 
 M_{\epsilon}\del L_{\epsilon}^{-1}
 \end{array}
 \right).
\end{equation}
In all,
(\ref{eq;06.1.2.2}) is same as the following,
modulo uniformly bounded terms
due to (\ref{eq;06.1.14.350}),
(\ref{eq;06.1.14.370}),
(\ref{eq;06.1.14.380})
and (\ref{eq;06.1.14.360}):
\begin{multline}
 \label{eq;06.1.14.400}
\frac{1}{1+|\lambda|^2(1-M_{\epsilon}^2)}
 \left(
 \begin{array}{cc}
 -\lambda\delbar\del \log L_{\epsilon} &
 -\lambda^2M_{\epsilon}\cdot\delbar\del L_{\epsilon}^{-1}\\
 0 & \lambda\delbar\del \log L_{\epsilon}
 \end{array}
 \right) \\
+
\frac{1}{1+|\lambda|^2(1-M_{\epsilon}^2)}
\frac{|\lambda|^2}{1+|\lambda|^2}
\frac{d\zbar\cdot dz}{|z|^2}
 \left(
\begin{array}{cc}
 \lambda \cdot M_{\epsilon}\cdot
 K_{\epsilon}\cdot L_{\epsilon}^{-2}
 & 0 \\
 2K_{\epsilon}\cdot L_{\epsilon}^{-1} &
 -\lambda\cdot M_{\epsilon}\cdot K_{\epsilon}\cdot L_{\epsilon}^{-2}
 \end{array}
 \right) \\
+\frac{1}{1+|\lambda|^2(1-M_{\epsilon}^2)}
 \frac{\lambda^2}{1+|\lambda|^2}
 \frac{d\zbar\cdot dz}{|z|^2}
 \left(\begin{array}{cc}
 -\lambdabar\cdot M_{\epsilon}\cdot K_{\epsilon}\cdot L_{\epsilon}^{-2}
 & 2(1+|\lambda|^2)L_{\epsilon}^{-3}\cdot K_{\epsilon}\\
 0 & \lambdabar \cdot M_{\epsilon}\cdot
  K_{\epsilon}\cdot L_{\epsilon}^{-2}
 \end{array}
 \right) \\
-\frac{\lambda}{1+|\lambda|^2(1-M_{\epsilon}^2)}
 \frac{d\zbar \cdot dz}{|z|^2}
 \left(
 \begin{array}{cc}
 L_{\epsilon}^{-2} & 0 \\
 2\lambdabar (1+|\lambda|^2)^{-1}
 M_{\epsilon}\cdot L_{\epsilon}^{-1}
 & -L_{\epsilon}^{-2}
 \end{array}
 \right).
\end{multline}
The summation of the last three term in (\ref{eq;06.1.14.400})
is as follows:
\begin{equation}
\frac{1}{1+|\lambda|^2(1-M_{\epsilon}^2)}
 \frac{d\zbar\cdot dz}{|z|^2}
 \left(
 \begin{array}{cc}
 -\lambda L_{\epsilon}^{-2} 
& 2\lambda^2 L_{\epsilon}^{-3}K_{\epsilon} \\
 2|\lambda|^2(1+|\lambda|^2)^{-1}
(K_{\epsilon}-M_{\epsilon})L_{\epsilon}^{-1}
 & \lambda L_{\epsilon}^{-2}
 \end{array}
 \right).
\end{equation}
By a direct calculation,
we can show the following equalities:
\[
 \delbar\del \log L_{\epsilon}
=-\frac{1}{L_{\epsilon}^2}\frac{d\zbar\cdot dz}{|z|^2},
\quad\quad
 \delbar\del L_{\epsilon}^{-1}
=\frac{2}{L_{\epsilon}^3}\frac{d\zbar\cdot dz}{|z|^2}
-\frac{\epsilon^2}{2}\frac{1}{L_{\epsilon}}
\frac{d\zbar \cdot dz}{|z|^2}.
\]
Therefore, (\ref{eq;06.1.14.400}) can be 
rewritten as follows:
\begin{multline}
 \label{eq;06.1.14.450}
 \frac{1}{1+|\lambda|^2(1-M_{\epsilon}^2)}
  \left( \begin{array}{cc}
 0 & 2\lambda^2L_{\epsilon}^{-3}(K_{\epsilon}-M_{\epsilon}) \\
 2|\lambda|^2(1+|\lambda|^2)^{-1}L_{\epsilon}^{-1}
 (K_{\epsilon}-M_{\epsilon})
 & 0
 \end{array}
 \right)\cdot\frac{d\zbar \cdot dz}{|z|^2} \\
+
 \frac{1}{1+|\lambda|^2(1-M_{\epsilon}^2)}
\left( \begin{array}{cc}
 0 & \lambda^2\epsilon^2M_{\epsilon}(2L_{\epsilon})^{-1} \\
 0 & 0
 \end{array}
\right)\cdot \frac{d\zbar\cdot dz}{|z|^2}.
\end{multline}
Due to $M_{\epsilon}=O\bigl(|z|^{4\epsilon}(1+\epsilon L_0)\bigr)$,
the second term in (\ref{eq;06.1.14.450}) can be ignored.
Due to Lemma \ref{lem;06.1.2.6} and
Lemma \ref{lem;06.1.2.7},
we have the uniform boundedness of
$(M_{\epsilon}-1)\cdot L_{\epsilon}^{-2}
 \cdot dz\cdot d\zbar/ |z|^2$
and 
$(K_{\epsilon}-1)\cdot L_{\epsilon}^{-2}
 \cdot dz\cdot d\zbar/|z|^2$.
Thus, the proof of Proposition \ref{prop;06.1.14.2}
is finished.
\hfill\qed

\subsection{A family of metrics of 
 a parabolic flat bundle on a disc}
\label{subsection;06.1.21.30}

\subsubsection{Simple case}
\label{subsubsection;06.2.7.1}

We put $X:=\Delta=\{z\in\cnum\,|\,|z|<1\}$
and $X^{\ast}:=\Delta-\{O\}$.
Let $V_l$ be a vector space over $\cnum$
with a base $\vece=(e_1,\ldots, e_l)$,
and let $N_l$ be the nilpotent endomorphism of $V_l$
given by $N_l\cdot e_{i+1}=e_{i}$ for $i=1,\ldots,l-1$
and $N_l\cdot e_l=0$.
We put $E_l:=\nbigo_{X}\otimes V_l$.
Then, $e_i$ naturally induce the frame of $E_l$,
which we denote by $\vecv=(v_1,\ldots,v_l)$.
The fiber $E_{|O}$ is naturally identified with $V$,
and we have $\vecv_{|O}=\vece$.
We have the logarithmic $\lambda$-connection
$\DDlambda_l$ of $E_l$ given by
$\DDlambda_lv_i=v_{i+1}\cdot dz/z$ for $i=1,\ldots,l-1$
and $\DDlambda_l v_l=0$.
The residue $\Res(\DDlambda)$ is given by $N_l$.
We have the weight filtration $W$  of $E_{|O}$
with respect to $N_l$.

We have the trivial parabolic structure $F$ of $E_l$.
Take a sufficiently small positive number $\epsilon$.
We consider the $\epsilon$-perturbation $F^{(\epsilon)}$
given by
$F^{(\epsilon)}_{k\epsilon}=W_k$
for $k=-l+1,-l+3\ldots,l-1$
in this case.

\vspace{.1in}

Let us fix a sufficiently small positive number $\epsilon_0$
such that $\rank E\cdot\epsilon_0<\eta/10$.
In the previous subsection,
we have constructed a family of metrics
$h_2^{(\epsilon)}$ $(0\leq \epsilon\leq \epsilon_0)$.
It induces the metric of $\Sym^{l-1}(E_2,\DDlambda_2)
\simeq (E_l,\DD_l)$,
which we denote by $h_l^{(\epsilon)}$.
The following property can be shown
by reducing to the case $l=2$.
\begin{itemize}
\item
 $h_l^{(0)}$ is the harmonic metric.
\item
 $h_{l}^{(\epsilon)}\lrarr h_{l}^{(0)}$
 for $\epsilon\to 0$,
 in the $C^{\infty}$-sense locally on $X^{\ast}$.
\item
 $\bigl|
 \Lambda_{\omega_{\epsilon}}
 G(h_{l}^{(\epsilon)})\bigr|_{h_{l}^{(\epsilon)}}
 <C$.
\item
 $h_l^{(\epsilon)}$ is adapted to the parabolic
 structure $F_l^{(\epsilon)}$.
\item
 Let $t_{\epsilon}:=\det(h_l^{(\epsilon)})\big/\det(h_l^{(0)})$.
 Then, $t_{\epsilon}$ and $t_{\epsilon}^{-1}$ are bounded,
 independently of $\epsilon$.
\end{itemize}

\begin{lem}
\label{lem;06.1.14.600}
Let $H_{\epsilon}= \bigl(h^{(\epsilon)}(v_i,v_j)\bigr)$.
Then, we have the following estimate
on $\{0<|z|<1/2\}$
with respect to $h_l^{(\epsilon)}$:
\[
 \Hbar_{\epsilon}^{-1}\cdot\bigl(\delbar+\lambda\del\bigr)
 \Hbar_{\epsilon}
=O(1)\cdot\frac{dz}{z}+O(1)\cdot\frac{d\zbar}{\zbar}
\]
\end{lem}
\pf
We see only $\Hbar_{\epsilon}^{-1}\del \Hbarepsilon$.
The term $\Hbarepsilon^{-1}\delbar\Hbarepsilon$
can be discussed in the same way.
We have only to check the case $l=2$.
As in Lemma \ref{lem;06.1.14.150},
we have only to see the $(1,1)$-entry, $(2,2)$-entry,
$L_{\epsilon}\times$ $(1,2)$-entry
and $L_{\epsilon}^{-1}\times$ $(2,1)$-entry
in the matrix valued function (\ref{eq;06.1.14.500}).
As is seen in Subsection \ref{subsubsection;06.1.14.501},
the term containing $\del M_{\epsilon}$ is bounded
with respect to $\omega_{\epsilon}$,
and the estimate is uniform for $\epsilon$.
Hence, we can ignore them.
Therefore, we have only to show that
$L_{\epsilon}^{-1}\del L_{\epsilon}
=-L_{\epsilon}\del L_{\epsilon}^{-1}$ 
is $O(1)\cdot dz/z$,
but it can be checked by a direct calculation.
\hfill\qed

\subsubsection{General case}
\label{subsubsection;06.1.21.15}

Let $(E,\vecF,\DDlambda)$ be
a parabolic flat $\lambda$-connection
on $(X,O)$.
Take a positive number $\eta$
such that $10\cdot \eta<\gap(E,\vecF)$.
We will use the metrics:
\begin{equation}
 \label{eq;06.1.15.1}
 \omega_{\epsilon}=
 \epsilon^2|z|^{\epsilon}\frac{dz\cdot d\zbar}{|z|^2}
+|z|^{2\eta}\frac{dz\cdot d\zbar}{|z|^2}.
\end{equation}
Here, $\epsilon$ will be $m^{-1}$ for some $m\in\seisuu_{>0}$
such that $10\rank(E)\cdot\epsilon<\eta$.
We take the $\epsilon$-perturbation $\vecF^{(\epsilon)}$
as in (II) of Subsection \ref{subsubsection;06.1.19.10}.
Let $a(\epsilon)$ be the numbers
which is denoted by $a(\epsilon,i)$
in the explanation there.

We have the endomorphism
$\Res (\DDlambda)$ of $\Gr^F_a(E)$.
It induces the generalized eigen decomposition
$\Gr^{F}_a(E)=
\bigoplus_{\alpha\in\cnum} \Gr^{F,\EE}_{a,\alpha}(E)$.
On $\Gr^{F,\EE}_{u}(E)$,
the endomorphism $\Res(\DDlambda)$ is decomposed as
$\alpha\cdot\id+N_u$,
where $u=(a,\alpha)\in\real\times\cnum$.
Let $W$ be the weight filtration of $N_u$
on $\Gr^{F,\EE}_{u}(E)$.
They induce the filtration $W$ of $\Gr^F_a(E)$.

For $u\in\real\times\cnum$,
we put $V_{u}:=\Gr^{F,\EE}_{u}(E)$
with the induced nilpotent map $N_{u}$.
Then, we can take an isomorphism:
\[
 (V_{u},N_u)\simeq
 \bigoplus_{i=1}^{m(u)} 
 \bigl(V_{l(u,i)},N_{l(u,i)}\bigr).
\]
We put 
$ (E_u,\DDlambda_u):=
\bigoplus
 \bigl(E_{l(u,i)},\DDlambda_{l(u,i)}\bigr)$.
Let $h_u^{\prime\,(\epsilon)}$ denote the metric of $E_u$
induced by $h_{l(u,i)}^{(\epsilon)}$
$(i=1,\ldots,m(u))$.
(See Subsection \ref{subsubsection;06.2.7.1}).

Let $Q(u)$ denote the logarithmic $\lambda$-flat bundle
of rank one for $u=(a,\alpha)$,
which is given by
$ \nbigo_{X}\!\cdot\! e$
with the $\lambda$-connection
$\DDlambda e=e\!\cdot\! \alpha\!\cdot\! dz/z$.
It is equipped with the family of the harmonic metrics
$h^{\prime\prime\,(\epsilon)}_{u,\epsilon}(e,e)
 =|z|^{-2a(\epsilon)}$.
Then, we obtain the vector bundle $E_0$
with the $\lambda$-connection $\DDlambda_0$
and the parabolic structure $F$,
as follows:
\[
 (E_0,\DDlambda_0)
=\bigoplus_{u}
 \bigl(
 E_u,\DDlambda_u \bigr)\otimes Q(u),
\quad\quad
 F_b(E_{0\,|\,O})
=\bigoplus_{a\leq b}
 E_{(a,\alpha)|O} \otimes Q(a,\alpha)_{|O}.
\]
The metrics $h^{\prime\,(\epsilon)}_u$
and $h_{u}^{\prime\prime\,(\epsilon)}$ 
induce the metric $h^{(\epsilon)}_u$
of $E_u\otimes Q(u)$.
Let $h^{(\epsilon)}_0$ denote the direct sum of them.
We can take a holomorphic isomorphism 
$\Psi:E_0\lrarr E$ 
satisfying the following conditions:
\begin{itemize}
\item 
It preserves the filtration $F$.
\item
 $\Gr^F(\Psi)\circ\Gr^F\Res\DDlambda
 =\Gr^F\Res\DDlambda_0\Gr^F(\Psi)$.
\end{itemize}
We identify $E_0$ and $E$ via $\Psi$.
The naturally induced metric of $E$ is denoted
by the same notation $h_0^{(\epsilon)}$.

\begin{lem}
The family $\bigl\{h_0^{(\epsilon)}\,\big|\,0\leq \epsilon\leq\epsilon_0\bigr\}$
of the hermitian metrics has the following properties:
\begin{itemize}
\item
$G(\DDlambda,h_0^{(\epsilon)})$ 
is uniformly bounded with respect to 
$(\omega_{\epsilon},h_0^{(\epsilon)})$.
\item
$\{h_0^{(\epsilon)}\,|\,\epsilon>0\}$ converges to $h_0^{(0)}$ 
in the $C^{\infty}$-sense locally on $X^{\ast}$.
\item
$h_0^{(\epsilon)}$ is adapted to 
the perturbed parabolic structure $F^{(\epsilon)}$.
\item
Let $t_{\epsilon}$ be determined by
$\det (h_0^{(\epsilon)})\big/\det (h_0^{(0)})$.
Then, $t_{\epsilon}$ and $t_{\epsilon}^{-1}$ are bounded,
independently from $\epsilon$.
\end{itemize}
\end{lem}
\pf
We check only the first claim.
The other claims are easy to see.
Let $f$ be determined by
$f\cdot dz/z=\DDlambda-\DDlambda_0$,
and we put
$f_{\epsilon}^{\dagger}:=f^{\dagger}_{h^{(\epsilon)}}$.
We put
$\DDlambdastar_{\epsilon}
:=\DDlambdastar_{h^{(\epsilon)}}$
and $\DDlambdastar_{0,\epsilon}
:=\DDlambdastar_{0,h^{(\epsilon)}}$.
Then, we have the following:
\begin{multline}
 G(\DDlambda,h_0^{(\epsilon)})
=\bigl[\DDlambda,\DDlambdastar_{\epsilon}\bigr]
=\Bigl[
 \DDlambda_0+f \frac{dz}{z},\,\,
 \DDlambdastar_{0,\epsilon}
+f^{\dagger}_{\epsilon}\frac{d\zbar}{\zbar}
\Bigr] \\
=G(\DDlambda_0,h_0^{(\epsilon)})
+\DDlambdastar_{0,\epsilon}(f)\frac{dz}{z}
+\DDlambda_{0}(f^{\dagger}_{\epsilon})
 \frac{d\zbar}{\zbar}
+[f,f^{\dagger}_{\epsilon}]\frac{dz\cdot d\zbar}{|z|^2}.
\end{multline}
We have the decomposition $f=\sum f_{u,u'}$,
where
$ f_{u,u'}\in
 Hom\bigl(E_u\otimes Q(u),\,\,
 E_{u'}\otimes Q(u') \bigr)$.
We have $f_{u,u'\,|\,O}=0$ unless $\alpha=\alpha'$
and $a>a'$.
Hence, there exist positive constants
$C$ and $N$ such that the following holds
for $0<\epsilon<\epsilon_0$:
\[
 |f|_{h_0^{(\epsilon)}}\leq 
 C\cdot |z|^{10\eta} L_{\epsilon}^{N},
\]
Here $N\cdot \epsilon$ is sufficiently smaller than $\eta$.
Hence, we have the following:
\[
 |f|_{h_0^{(\epsilon)}}\leq C\cdot |z|^{9\eta},
\quad\quad
 [f,f^{\dagger}_{\epsilon}]=O\bigl(|z|^{18\eta}\bigr).
\]

We have the induced frames $\vecv_u$
of $E_u\otimes Q(u)$.
They induce the frame $\vecv$ of $E_0$.
Let $B$ and $A_0$ be determined by 
$ f\vecv=\vecv\cdot B\cdot dz/z$
and
$\DDlambda_0\vecv=\vecv A_0\cdot dz/z$.
Then, we have the following:
\[
 \bigl[
 \DDlambda_0,f^{\dagger}
\bigr]\vecv=\vecv
 \left(
 \nbigd B^{\dagger}_{\epsilon}\frac{d\zbar}{\zbar}
+[A_0,B^{\dagger}_{\epsilon}]\frac{dz\cdot d\zbar}{|z|^2}
 \right).
\]
Here we put $\nbigd=\delbar+\lambda\del$
and 
$B^{\dagger}_{\epsilon}=
 \Hbar^{-1}_{\epsilon} \cdot\lefttop{t}\overline{B}\cdot
 \Hbar_{\epsilon}$,
where $H_{\epsilon}=H(h^{(\epsilon)}_0,\vecu)$.
Since $B^{\dagger}_{\epsilon}$ is sufficiently small
with respect to $(\omega_{\epsilon},h^{(\epsilon)}_0)$,
$[A_0,B^{\dagger}_{\epsilon}]$ is also sufficiently small.
Corresponding to the decomposition $f=\sum f_{u,u'}$,
we have $B=\sum B_{u,u'}$.
Then, the following holds:
\[
\bigl(
 B^{\dagger}_{\epsilon}
\bigr)_{u,u'}
=\Hbar^{-1}_{u',\epsilon}
 \lefttop{t}\overline{B}_{u',u}
 \Hbar_{u,\epsilon}.
\]
Here $H_{u,\epsilon}:=H(h^{(\epsilon)}_u,\vecv_u)$.
Hence, we obtain the following:
\[
 \bigl(
 \nbigd B^{\dagger}_{\epsilon}\bigr)
 _{u,u'}\frac{d\zbar}{\zbar}
=\Hbar^{-1}_{u',\epsilon}\cdot
 (\nbigd\lefttop{t}\overline{B}_{u',u})\cdot
 \Hbar_{u,\epsilon}
-\Hbar^{-1}_{u',\epsilon}\nbigd\Hbar_{u',\epsilon}\cdot
 (B^{\dagger}_{\epsilon})_{u,u'}
+(B^{\dagger}_{\epsilon})_{u,u'}\cdot
 \Hbar^{-1}_{u,\epsilon}\nbigd\Hbar_{u,\epsilon}.
\]
Since $B$ is holomorphic,
we have
$ \Hbar^{-1}_{u',\epsilon}\cdot
\bigl(\nbigd\lefttop{t}\overline{B}_{u',u}\bigr)\cdot
 \Hbar_{u,\epsilon}\cdot d\zbar/\zbar=0$.
We put $H'_{u\,\epsilon}:=H(h^{\prime\,(\epsilon)}_u,\vecv_u)$.
Then, we have $H_{u,\epsilon}=|z|^{-2a}H_{u,\epsilon}'$,
and the following holds
with respect to $h_0^{(\epsilon)}$
due to Lemma \ref{lem;06.1.14.600}:
\[
 \Hbar_{u,\epsilon}^{-1}\nbigd \Hbar_{u,\epsilon}
=-a\left(\lambda\frac{dz}{z}+\frac{d\zbar}{\zbar}\right)
+\Hbar^{\prime\,-1}_{u,\epsilon}\nbigd \Hbar'_{u,\epsilon}
=O(1)\frac{dz}{z}+O(1)\frac{d\zbar}{\zbar}.
\]
Since $(B^{\dagger}_{\epsilon})_{u,u'}$ is small
with respect to $(\omega_{\epsilon},h_0^{(\epsilon)})$,
$(B^{\dagger}_{\epsilon})_{u,u'}
\cdot \Hbar^{-1}_{u,\epsilon}
 \del \Hbar_{u,\epsilon}$ is also small.
Therefore,
$ \DDlambda_0 f^{\dagger}\cdot d\zbar/\zbar$
is small with respect to
$(\omega_{\epsilon},h_0^{(\epsilon)})$.
It also follows that
$\DDlambdastar_{0,\epsilon}f\cdot dz/z$
is small.
Thus we are done.
\hfill\qed

\subsection{Proof of Proposition \ref{prop;06.1.21.1}}
\label{subsection;06.1.21.3}

\subsubsection{Construction of a family of initial metrics}

Let $\eta$ be a small positive number such that
$\eta<\gap(E,\vecF)/10$.
Let $\epsilon_0$ be a small positive number
such that $10\rank E\cdot \epsilon_0<\eta$.
For any $0\leq\epsilon<\epsilon_0$,
let us take $\omega_{\epsilon}$ be the Kahler forms of $C-D$
with the following properties:
\begin{itemize}
\item
 Let $(U_P,z)$ be a holomorphic coordinate
 around $P\in D$ such that $z(P)=0$,
 and then $\omega_{\epsilon}$ is given by (\ref{eq;06.1.15.1}).
\item
 $\omega_{\epsilon}\lrarr \omega_0$ for $\epsilon\lrarr 0$
in the $C^{\infty}$-sense locally on $X-D$.
\end{itemize}

\begin{lem}
 \label{lem;06.1.21.20}
We can construct a family of metrics
$h_{0}^{(\epsilon)}$ of $E_{|C-D}$
with the following properties:
\begin{itemize}
\item
 $h^{(\epsilon)}_0$ is adapted to the perturbed 
parabolic structure $\vecF^{(\epsilon)}$.
\item
 $h^{(\epsilon)}_{0}\lrarr h^{(0)}_{0}$
 in the $C^{\infty}$-sense locally on $C-D$.
\item
 $G(h_{0}^{(\epsilon)})$ is uniformly bounded
 with respect to $(\omega_{\epsilon},h_{0}^{(\epsilon)})$.
\item
We put $t_{\epsilon}:=\det (h_0^{(\epsilon)})\big/\det (h_0^{(0)})$.
Then, $t_{\epsilon}$ and $t_{\epsilon}^{-1}$ are bounded
 independently from $\epsilon$.
\end{itemize}
\end{lem}
\pf
We construct a $C^{\infty}$-metric
of $E$ on $\bigcup_{P\in D}(U_P-\{P\})$,
by applying the construction given
in Subsection \ref{subsubsection;06.1.21.15}
to $(E,\vecF,\DDlambda)_{|U_P}$
for each $P\in D$,
and then we prolong it to a $C^{\infty}$-metric
of $E$ on $C-D$.
\hfill\qed

\vspace{.1in}

Let $R(\det h_0^{(0)})$ denote the curvature of 
the metrized holomorphic bundle $\det(E,d'',h_0^{(0)})$,
where $d''$ denote the $(0,1)$-part of $\DDlambda$.
Since $\det h_0^{(0)}$ gives the harmonic metric
around $D$ due to our construction,
$R(\det h_0^{(0)})$ vanishes around $D$.
We also have 
$\int R(\det h_0^{(0)})=-2\pi\sqrt{-1}\cdot \pardeg(E,\vecF)=0$.
Let us take the $C^{\infty}$-function $\chi_0$ on $C$
and satisfies the equality
$\rank(E)\cdot \delbar\del \chi_0+
R\bigl(\det (h_0^{(0)})\bigr)=0$.
We put $h_{in}^{(0)}:=h_0^{(0)}\cdot \exp\bigl(\chi_0\bigr)$.
Then, $R\bigl(\det h_{in}^{(0)}\bigr)=0$,
i.e.,
$\det h_{in}^{(0)}$ is a harmonic metric of $\det (E,\DDlambda)$.
Let $\chi_{\epsilon}$ be the functions determined by 
$\det(h_{in}^{(0)})=
 \det (h_{0}^{(\epsilon)})\cdot \exp\Bigl(\rank(E)\cdot
 \chi_{\epsilon}\Bigr)$.
The following claims immediately follows from
Lemma \ref{lem;06.1.21.20}.
\begin{itemize}
\item
$\chi_{\epsilon}$ and $-\chi_{\epsilon}$ are bounded on $C$,
independently from $\epsilon$.
\item
$\chi_{\epsilon}\lrarr 0$ in the $C^{\infty}$-sense 
locally on $C-D$.
\end{itemize}

We put
$h_{in}^{(\epsilon)}:=h_0^{(\epsilon)}\cdot 
 \exp\bigl(\chi_{\epsilon}\bigr)$,
which is the metric of $E_{|C-D}$.
\begin{lem}
The following claims are easy to check.
\begin{itemize}
\item
 $h_{in}^{(\epsilon)}$ is adapted to 
the parabolic structure $\vecF^{(\epsilon)}$.
\item
 $h^{(\epsilon)}_{in}\lrarr h^{(0)}_{in}$
 in the $C^{\infty}$-sense locally on $C-D$.
\item
 $G(h_{in}^{(\epsilon)})$ is uniformly bounded
 with respect to $(\omega_{\epsilon},h_{in}^{(\epsilon)})$.
\item
 $\det h^{(\epsilon)}_{in}$ is harmonic,
 and we have $\det h^{(\epsilon)}_{in}=\det h^{(0)}_{in}$.
\end{itemize}
In other words, they give initial metrics for
$(E,\vecF^{(\epsilon)},\DDlambda)$
in the sense of Lemma {\rm\ref{lem;06.1.13.200}},
and their pseudo curvature satisfy some uniform finiteness.
\hfill\qed
\end{lem}

\subsubsection{$L_1^2$-finiteness of the sequence}

Due to Proposition \ref{prop;06.1.13.250},
we obtain the harmonic metrics $h^{(\epsilon)}$
for $(E,\vecF^{(\epsilon)},\DDlambda)$
such that $\det h^{(\epsilon)}=\det h^{(0)}_{in}$.
Due to Lemma \ref{lem;06.1.15.70},
we have the following inequalities for any $\epsilon$:
\begin{equation}
 \label{eq;06.1.15.20}
 M_{\omega_{\epsilon}}(h_{in}^{(\epsilon)},h^{(\epsilon)})\leq 0.
\end{equation}
Let $s^{(\epsilon)}$ be determined by
$h^{(\epsilon)}=h^{(\epsilon)}_{in}s^{(\epsilon)}$.
Due to Lemma \ref{lem;06.1.10.6}, (\ref{eq;06.1.15.20})
and $\det s^{(\epsilon)}=1$,
there exists a positive constant $A$
which is independent on $\epsilon$,
with the following property:
\begin{equation}
 \label{eq;06.1.15.100}
 \bigl|
 s^{(\epsilon)}
 \bigr|_{h^{(\epsilon)}_{in}}
\leq A,
\quad
  \bigl|
 s^{(\epsilon)\,-1}
 \bigr|_{h^{(\epsilon)}_{in}}
\leq A.
\end{equation}
Let $\DDlambdastar_{in}$ be the operator
obtained from $\DDlambda$, $\omega_{\epsilon}$ and $h_{in}^{(\epsilon)}$
as in Subsection \ref{subsubsection;06.1.10.7}.
We have the following equalities:
\[
 \Delta_{\omega_{\epsilon}}^{\lambda} \tr s^{(\epsilon)}
=-\sqrt{-1}\tr\bigl(
 s^{(\epsilon)}\Lambda_{\omega_{\epsilon}}
 G(h^{(\epsilon)}_{in})\bigr)
+\sqrt{-1}\tr\Bigl(
 \Lambda_{\omega_{\epsilon}}
 \DDlambda s^{(\epsilon)}
\cdot (s^{(\epsilon)})^{-1} \cdot
 \DDlambdastar_{in} s^{(\epsilon)}
 \Bigr).
\]
See Remark \ref{rem;06.1.26.20}
for $\Delta^{\lambda}_{\omega_{\epsilon}}$.
\begin{lem}
\label{lem;08.2.5.10}
We have 
$\int\Delta_{\omega_{\epsilon}}^{\lambda} 
 \tr s^{(\epsilon)}
 \dvol_{\omega_{\epsilon}}=0$.
\end{lem}
\pf
Let $g$ be a $C^{\infty}$-Kahler metric of $C$.
We have only to show
$\int\Delta_{g}^{\lambda}
 \tr s^{(\epsilon)}
 \dvol_{g}=0$.
We have the following:
\[
 \Delta_{g}^{\lambda} \tr s^{(\epsilon)}
=-\sqrt{-1}\tr\bigl(
 s^{(\epsilon)}\Lambda_{g}
 G(h^{(\epsilon)}_{in})\bigr)
-\bigl|
 \DDlambda s^{(\epsilon)}\cdot
 (s^{(\epsilon)})^{-1/2}
 \bigr|_{h^{(\epsilon)},g}^2
\]
Since 
$\bigl|
 \tr\bigl(
 s^{(\epsilon)}\Lambda_{g}
 G(h^{(\epsilon)}_{in})\bigr)
 \bigr|$ is $O(|z|^{\epsilon-2})$,
we can take a bounded function
$a_{\epsilon}$ such that
$\Delta_ga_{\epsilon}
=\bigl|
 \tr\bigl(
 s^{(\epsilon)}\Lambda_{g}
 G(h^{(\epsilon)}_{in})\bigr)
 \bigr|$.
Hence, we obtain 
$\int_{X-D}\bigl|
 \DDlambda s^{(\epsilon)}\cdot
 (s^{(\epsilon)})^{-1/2}
 \bigr|_{h^{(\epsilon)},g}^2<\infty$,
due to Lemma 2.2 of \cite{s2}.
Since $s^{(\epsilon)}$ is bounded
with respect to $h^{(\epsilon)}$,
we obtain
$\int_{X-D}\bigl|
 \DDlambda s^{(\epsilon)}
 \bigr|_{h^{(\epsilon)},g}^2<\infty$.
Then, it is easy to obtain the vanishing
$\int\Delta_{\omega_{\epsilon}}^{\lambda} 
 \tr s^{(\epsilon)}
 \dvol_{\omega_{\epsilon}}=0$
by Stokes formula
and Lemma 5.2 of \cite{s1}.
\hfill\qed

\vspace{.1in}
Then, there exists a positive constant $A'$
such that the following holds:
\begin{equation}
 \label{eq;06.1.15.101}
 \int |\DDlambda s^{(\epsilon)}\cdot s^{(\epsilon)\,-1/2}|
 _{h_{in}^{(\epsilon)},\omega_{\epsilon}}^2
\dvol_{\omega_{\epsilon}}
\leq A'.
\end{equation}
In particular,
we obtain $\bigl\|\DDlambda s^{(\epsilon)}\bigr\|
 _{L^2,\omega_{\epsilon},h_{in}^{(\epsilon)}}$
 is bounded for $0<\epsilon<\epsilon_0$.

\subsubsection{The end of the proof of Proposition \ref{prop;06.1.21.1}}

Let $Q$ be a point of $C-D$.
Let $(U,z)$ be a holomorphic coordinate around $Q$
such that $z(Q)=0$
and $U\simeq\Delta=\{z\,|\,|z|<1\}$.
We use the standard metric $g=dz\cdot d\zbar$ of $U$.
The harmonic bundle $(E,\DDlambda,h^{(\epsilon)})$
induces the Higgs bundle
$(E,\delbar_{\epsilon},\theta_{\epsilon})$.
We have
$\theta_{\epsilon}=f_{\epsilon}\cdot dz$ on $U$.
On the other hand,
we also obtain
$\delbar_{in,\epsilon}$ and $\theta_{in,\epsilon}$
from $(E,\DDlambda,h^{(\epsilon)}_{in})$,
although $\delbar_{in,\epsilon}\bigl(\theta_{in,\epsilon}\bigr)=0$
is not satisfied, in general.
Let $\delta_{in,\epsilon}'$ be the $(1,0)$-operator obtained from
$h^{(\epsilon)}_{in}$ and $d''$,
as in Subsection \ref{subsubsection;06.1.10.7}.
Then, we have the relation:
\begin{equation}
 \label{eq;06.1.15.125}
 \theta_{\epsilon}=\theta_{in,\epsilon}
-\frac{\lambda}{1+|\lambda|^2}
\bigl(
 s^{(\epsilon)\,-1}\cdot\delta_{in,\epsilon}'s^{(\epsilon)}
\bigr).
\end{equation}
Due to (\ref{eq;06.1.15.100}), (\ref{eq;06.1.15.101})
and (\ref{eq;06.1.15.125}),
there exists a positive constant $C_0$
such that
$ \int_U|f_{\epsilon}|_{h^{(\epsilon)}}^2
 \dvol_g<C_0$
holds for any $0<\epsilon<\epsilon_0$.
Hence, the following inequality holds
for some positive constants $C_i$ $(i=1,2,3)$
and for any $0<\epsilon<\epsilon_0$:
\begin{equation}
 \label{eq;06.1.15.121}
 \int_U\log|f_{\epsilon}|_{h^{(\epsilon)}}^2\dvol_g
\leq C_1+\int_U 
 C_2\cdot |f_{\epsilon}|_{h^{(\epsilon)}}^2
 \dvol_g
\leq C_3.
\end{equation}
Recall the fundamental inequality for the Higgs field
of a harmonic bundle \cite{s2}:
\begin{equation}
 \label{eq;06.1.15.120}
 \Delta_g \log |f_{\epsilon}|_{h^{(\epsilon)}}^2
\leq -\frac{
 \bigl|[f_{\epsilon},f_{\epsilon}^{\dagger}]
 \bigr|_{h^{(\epsilon)}}^2}
 {|f_{\epsilon}|_{h^{(\epsilon)}}^2}\leq 0.
\end{equation}
Due to (\ref{eq;06.1.15.121}) and (\ref{eq;06.1.15.120}),
there exists a positive constant $C_4$
such that the following holds 
for any $Q'\in U(1/2):=\{|z|<1/2\}$:
\begin{equation}
 \bigl|f_{\epsilon}(Q')\bigr|^2_{h^{(\epsilon)}_{in}}
 \leq C_4.
\end{equation}
By using (\ref{eq;06.1.15.125}),
we obtain that 
$\delta'_{in,\epsilon}s^{(\epsilon)}$ is uniformly bounded
with respect to $(\omega_{\epsilon},h^{(\epsilon)}_{in})$
on $U(1/2)$.

Since $\theta^{\dagger}_{\epsilon}$ is the adjoint of
$\theta_{\epsilon}$,
we obtain the uniform boundedness of $\theta^{\dagger}_{\epsilon}$
on $U(1/2)$.
Let $\delta''_{in,\epsilon}$ be the operator obtained from
$h_{in}^{(\epsilon)}$ and $d'$
as in Subsection \ref{subsubsection;06.1.10.7},
where $d'$ denotes the $(1,0)$-part of $\DDlambda$.
Then, we also obtain the uniform boundedness of
$\delta''_{in,\epsilon}s^{(\epsilon)}$ on $U(1/2)$.
Hence, 
$\DDlambdastar_{in,\epsilon} s^{(\epsilon)}$ is uniformly
bounded on $U(1/2)$,
where $\DDlambdastar_{in,\epsilon}=
 \delta'_{in,\epsilon}-\delta''_{in,\epsilon}$.
Since we have
$d''=\lambdabar^{-1}\bigl(
 \delta''_{in,\epsilon}
+(1+|\lambda|^2)\theta^{\dagger}_{in,\epsilon} \bigr)$
and
$d'=\lambda\delta'_{in,\epsilon}
 +(1+|\lambda|^2)\theta_{in,\epsilon}$,
we also obtain $\DDlambda s^{(\epsilon)}$ is uniformly
bounded on $U(1/2)$.
Recall the formula
$ \DDlambda\DDlambdastar_{in} s^{(\epsilon)}
=s^{(\epsilon)}\cdot G(h^{(\epsilon)}_{in})
+\DDlambda s^{(\epsilon)}
\cdot s^{(\epsilon)\,-1}\cdot
 \DDlambdastar_{in}s^{(\epsilon)}$.
Thus $\DDlambda\DDlambdastar_{in} s^{(\epsilon)}$
is also uniformly bounded on $U(1/2)$.
Therefore, $\{s^{(\epsilon)}\}$ is $L_2^p$-bounded
for any $p>1$ and $U(1/2)$.
By taking an appropriate subsequence $(\epsilon_i)$,
$s^{(\epsilon_i)}$ weakly converges
to some $\widetilde{s}$ in $L_2^p$ locally on $C-D$.

It is easy to see that
$h^{(0)}_{in}\cdot\widetilde{s}$ is a harmonic metric.
We have  $\det \widetilde{s}=1$.
We also have the boundedness of $\widetilde{s}$ and $\widetilde{s}^{-1}$
with respect to $h^{(0)}_{in}$.
Thus, we have $h^{(0)}_{in}\cdot\widetilde{s}=h^{(0)}$,
i.e.,
the sequence $\{h^{(\epsilon_i)}\}$ converges to $h^{(0)}$
weakly in $L_2^p$ locally on $C-D$.

Although we take a subsequence in the above argument,
we can conclude that
$h^{(\epsilon)}$ converges to $h^{(0)}$
weakly in $L_2^p$ locally on $C-D$,
due to a general argument.
We can also obtain the $C^{\infty}$-convergence
by a standard bootstrapping argument.
In the above argument,
the convergence of $\{\theta^{(\epsilon)}\}$ is also proved.
\hfill\qed

\begin{rem}
 \label{rem;06.1.21.10}
As for the proof of Proposition {\rm\ref{prop;06.1.18.80}},
we take a $C^{\infty}$-metric $h_{in}$ of $(E,\vecF,\DDlambda)$
such that each restriction $h_{in\,|\,C_t}$ is an
initial metric.
Let $s$ be determined by $h_H=h_{in}\cdot s$.
By applying the same argument,
we obtain the continuity of $s$.
Similarly for $\theta_H$.
\hfill\qed
\end{rem}

%% file: 5.tex
We will prove our main existence theorem
of pluri-harmonic metric for parabolic $\lambda$-flat bundle,
which is adapted to the parabolic structure.
(See Subsection 3.3 of \cite{mochi4}
for the adaptedness.)

\subsection{Preliminary}
\label{subsection;06.1.21.31}

Let $C$ be a smooth projective curve
with a simple effective divisor $D$.
Let $(E,\vecF,\DDlambda)$ be
a stable parabolic $\lambda$-flat bundle
on $(C,D)$ with $\pardeg(E,\vecF)=0$.
For each $P\in D$,
let $(U_P,z)$ be a holomorphic coordinate around $P$
such that $z(P)=0$.
Let $\vecF^{(\epsilon)}$ be an $\epsilon$-perturbation
as in (II) of Subsection \ref{subsubsection;06.1.19.10}
for $\epsilon=m^{-1}$.
We have $h_0^{(\epsilon)}$ be harmonic metrics
for $(E,\vecF^{(\epsilon)},\DDlambda)$.
We assume $\det h_0^{(\epsilon)}=\det h_0^{(0)}$.
As shown in Proposition \ref{prop;06.1.21.1},
$h_0^{(\epsilon)}$ converges to $h_0^{(0)}$
in the $C^{\infty}$-sense locally on $C-D$.
Let $N$ be a large positive number,
for example $N>10$.
In this subsection,
we use Kahler metrics $g_{\epsilon}$ $(\epsilon\geq\,0)$
of $C-D$ which are as follows on $U_P$ for each $P\in D$:
\[
\bigl(
 \epsilon^{N+2}|z|^{2\epsilon}
+|z|^{2}
\bigr)\frac{dz\cdot d\zbar}{|z|^2}.
\]
We assume that $\{g_{\epsilon}\}$ converges
to $g_{0}$ for $\epsilon\lrarr 0$
in the $C^{\infty}$-sense
locally on $C-D$.

\begin{prop}
 \label{prop;06.1.18.15}
Let $h^{(\epsilon)}$ $(\epsilon>0)$ be hermitian metrics of $E_{|C-D}$
with the following properties:
\begin{enumerate}
\item
 Let $s^{(\epsilon)}$ be determined by
 $h^{(\epsilon)}=h_0^{(\epsilon)}\cdot s^{(\epsilon)}$.
 Then, $s^{(\epsilon)}$ is bounded with respect to
 $h^{(\epsilon)}_0$,
and we have $\det s^{(\epsilon)}=1$.
 We also have the finiteness
$ \bigl\|\DDlambda s^{(\epsilon)}\bigr\|_{2,h_0^{(\epsilon)},g_{\epsilon}}<\infty$.
(The estimates may depend on $\epsilon$.)
\item
We have
$  \|G(h^{(\epsilon)})\|_{2,h^{(\epsilon)},g_{\epsilon}}
 <\infty$
and
$ \lim_{\epsilon\to 0} \|G(h^{(\epsilon)})\|_{2,h^{(\epsilon)},g_{\epsilon}}=0$.
\end{enumerate}
Then, the following claims hold.
\begin{itemize}
\item
The sequence $\{s^{(\epsilon)}\}$ is
weakly convergent to the identity
in $L_1^2$ locally on $C-D$.
\item
$\bigl\{
 \sup_{P\in C-D} |s^{(\epsilon)}_{|P}|_{h_0^{(\epsilon)}}
 \,\big|\,\epsilon>0  \bigr\}$ 
and 
$\bigl\{\sup_{P\in C-D}
 |(s^{(\epsilon)})^{-1}_{|P}|_{h^{(\epsilon)}_0}
 \,\big|\,
 \epsilon>0\bigr\}$ are bounded.
\end{itemize}
\end{prop}
\pf
To begin with,
we remark that we have only to show the existence
of a subsequence $\{s^{(\epsilon_i)}\}$
with the desired properties as above.
We put 
$\|s^{(\epsilon)}\|_{\infty,h_0^{(\epsilon)}}:=
\sup_{P\in C-D} \bigl|
s^{(\epsilon)}_{|P}\bigr|_{h^{(\epsilon)}_0}$.
For any point $P\in C-D$,
let 
$ SE(s^{(\epsilon)})(P)$ denote the maximal eigenvalue
of $s^{(\epsilon)}_{|P}$.
There exists a constant $0<C_1<1$
such that
$ C_1\cdot |s^{(\epsilon)}_{|P}|_{h_0^{(\epsilon)}}
\leq SE(s^{(\epsilon)})(P)
\leq |s^{(\epsilon)}_{|P}|_{h^{(\epsilon)}_0}$.
We have $\det s^{(\epsilon)}_{|P}=1$.
Hence, it is easy to see
$\log \tr s^{(\epsilon)}_{|P}\geq \log \rank(E)\geq 0$.
We also have
$SE(s^{(\epsilon)})(P)\geq 1$
for any $P$.

Let us take $b_{\epsilon}>0$
satisfying 
$2\leq b_{\epsilon}\cdot\sup SE(s^{(\epsilon)})(P)
\leq 2+\epsilon$.
We put
 $\widetilde{s}^{(\epsilon)}=b_{\epsilon}s^{(\epsilon)}$
and 
$\widetilde{h}^{(\epsilon)}:=
 h^{(\epsilon)}_0\cdot \widetilde{s}^{(\epsilon)}$.
Then, $\widetilde{s}^{(\epsilon)}$
are uniformly bounded with respect to $h_0^{(\epsilon)}$.
We remark
$G(\widetilde{h}^{(\epsilon)})
=G(h^{(\epsilon)})$.
We also remark that $h^{(\epsilon)}$
and $\widetilde{h}^{(\epsilon)}$ induce
the same metric of $\End(E)$.

\begin{lem}
 \label{lem;06.1.21.50}
After going to an appropriate subsequence,
$\bigl\{\widetilde{s}^{(\epsilon_i)}\bigr\}$ converges
to a positive constant multiplication,
weakly in $L_1^2$ locally on $C-D$.
\end{lem}
\pf
We have the following
(Subsection \ref{subsubsection;06.1.18.3}):
\begin{equation}
 \label{eq;06.1.15.250}
 \Delta^{\lambda}_{g_0,h_0^{(\epsilon)}} \widetilde{s}^{(\epsilon)}
=\widetilde{s}^{(\epsilon)}
 \sqrt{-1}\Lambda_{g_0}
 G(\widetilde{h}^{(\epsilon)})
+\sqrt{-1}\Lambda_{g_0}
 \DDlambda\widetilde{s}^{(\epsilon)}
 \bigl(\widetilde{s}^{(\epsilon)\,-1}\bigr)
 \DDlambdastar_{h_0^{(\epsilon)}} \widetilde{s}^{(\epsilon)}.
\end{equation}
We can show
$\int
 \Delta_{g_0}^{\lambda}
 \tr \widetilde{s}^{(\epsilon)}
 \cdot\dvol_{g_0}=0$,
by the same argument as the proof of
Lemma \ref{lem;08.2.5.10},
we obtain the following inequality
from (\ref{eq;06.1.15.250})
and the uniform boundedness of $\widetilde{s}^{(\epsilon)}$:
\begin{multline}
\int\bigl|
 \DDlambda \widetilde{s}^{(\epsilon)}\cdot
 \widetilde{s}^{(\epsilon)\,-1/2}
 \bigr|_{g_0,h_0^{(\epsilon)}}^2
 \dvol_{g_0}\leq
 A\cdot\int
 \bigl|\tr \Lambda_{g_0} G(\widetilde{h}^{(\epsilon)})\bigr|
 \cdot\dvol_{g_0} \\
=A\cdot\int 
 \bigl|\tr \Lambda_{g_{\epsilon}}G(\widetilde{h}^{(\epsilon)})
 \bigr|\cdot\dvol_{g_{\epsilon}}
\leq
 A'\cdot\bigl\|G(\widetilde{h}^{(\epsilon)})
 \bigr\|_{2,\widetilde{h}^{(\epsilon)},g_{\epsilon}}.
\end{multline}
In particular,
we obtain the uniform estimate
$ \bigl\|\DDlambda \widetilde{s}^{(\epsilon)}
 \bigr\|_{2,g_0,h_0^{(\epsilon)}}^2
\leq
 A''\cdot\bigl\| G(\widetilde{h}^{(\epsilon)})
 \bigr\|_{2,\widetilde{h}^{(\epsilon)},g_{\epsilon}}$.
Therefore, the sequence
$\bigl\{\widetilde{s}^{(\epsilon)}\bigr\}$ is $L^2_1$-bounded
on any compact subset of $C-D$.
By taking an appropriate subsequence,
it is weakly $L^2_1$-convergent
locally on $C-D$.
Let $\widetilde{s}^{(\infty)}$ denote the weak limit.
We obtain $\DDlambda \widetilde{s}^{(\infty)}=0$.
We also know that
$\widetilde{s}^{(\infty)}$ is bounded with respect to
$h^{(0)}_0$.
Therefore,
$\widetilde{s}^{(\infty)}$ gives an automorphism of
$(E,\vecF,\DDlambda)$.
Due to the stability of $(E,\vecF,\DDlambda)$,
$\widetilde{s}^{(\infty)}$ is a constant multiplication.

We would like to show $\widetilde{s}^{(\infty)}\neq 0$.
Let us take any point $Q_{\epsilon}\in C-D$
satisfying the following:
\[
  SE(s^{(\epsilon)})(Q_{\epsilon})
\geq
\frac{9}{10}\cdot
 \sup_{P\in C-D} SE(s^{(\epsilon)})(P).
\]
Then, we have
$ \log \tr \widetilde{s}^{(\epsilon)}(Q_{\epsilon})
\geq \log (9/5)$.
By taking an appropriate subsequence,
we may assume the sequence $\{Q_{\epsilon}\}$ converges
to a point $Q_{\infty}$.
We have two cases
(i) $Q_{\infty}\in D$
(ii) $Q_{\infty}\not\in D$.
We discuss only the case (i).
The other case is similar and easier.

We use the coordinate neighbourhood $(U,z)$
such that $z(Q_{\infty})=0$.
For any point $P\in U$,
we put $\Delta(P,r):=\bigl\{Q\in U\,\big|\,|z(P)-z(Q)|<r\bigr\}$.
When $\epsilon$ is sufficiently small,
$Q_{\epsilon}$ is contained in $\Delta(Q_{\infty},1/2)=\{|z|<1/2\}$.
Let $g=dz\cdot d\zbar$ denote the standard metric of $U$.
We have the following inequality on $U-\{Q_{\infty}\}$
(see Subsection \ref{subsubsection;06.1.18.3}):
\begin{equation}
 \label{eq;08.2.5.15}
 \Delta^{\lambda}_g\log \tr\widetilde{s}^{(\epsilon)}
\leq
 \bigl|\Lambda_g G(\widetilde{h}^{(\epsilon)})
 \bigr|_{\widetilde{h}^{(\epsilon)}}.
\end{equation}
Let $B^{(\epsilon)}$ be the endomorphism of $E$
determined as follows:
\[
  G(\widetilde{h}^{(\epsilon)})
=G(h^{(\epsilon)})
=B^{(\epsilon)}\cdot \frac{dz \cdot d\zbar}{|z|^2}
\]
Then, we have the following estimate
for some constant $A>0$
which is independent of $\epsilon$:
\[
  \int \bigl|B^{(\epsilon)}\bigr|^2_{\widetilde{h}^{(\epsilon)}_0}
 \bigl(\epsilon^{N+1}|z|^{2\epsilon}+|z|^2\bigr)^{-1}
 \frac{\dvol_{g}}{|z|^2}
\leq
 A\int
 \bigl| G(\widetilde{h}^{(\epsilon)})
 \bigr|^2_{\widetilde{h}^{(\epsilon)},g_{\epsilon}}
 \dvol_{g_{\epsilon}}.
\]
Here $A$ denotes a constant independent of $\epsilon$.
Due to Proposition 2.16 in \cite{mochi4},
there exist $v^{(\epsilon)}$ such that
the following inequalities hold
for some positive constant $A'$
which is independent of $\epsilon$:
\[
 \delbar\del v^{(\epsilon)}
=\bigl|B^{(\epsilon)}\bigr|_{\widetilde{h}^{(\epsilon)}}
 \frac{dz \cdot d\zbar}{|z|^2},
\quad\quad
 \bigl|v^{(\epsilon)}(z)\bigr|
\leq
 A'\cdot\bigl(
 \epsilon^{(N-1)/2}|z|^{\epsilon}+|z|^{1/2}
 \bigr)\cdot
 \bigl\| G(\widetilde{h}^{(\epsilon)})
 \bigr\|_{2,\widetilde{h}^{(\epsilon)},g_{\epsilon}}
\]
Then, we have
$ \Delta^{\lambda}_g\bigl(
 \log \tr \widetilde{s}^{(\epsilon)}-v^{(\epsilon)}
 \bigr)\leq 0$ on $U-\{Q_{\infty}\}$.
Since $\log\tr\widetilde{s}^{(\epsilon)}-v^{(\epsilon)}$
is bounded from above,
the inequality holds on $U$.
Therefore, we obtain the following:
\[
 \log \tr\widetilde{s}^{(\epsilon)}(Q_{\epsilon})
-v^{(\epsilon)}(Q_{\epsilon})
\leq
A''\cdot
 \int_{\Delta(Q_{\epsilon},1/2)}
 \Bigl(
 \log \tr\widetilde{s}^{(\epsilon)}-v^{(\epsilon)}
\Bigr)\cdot\dvol_g.
\]
Here $A''$ denotes a positive constant
which is independent of $\epsilon$.
Then, we obtain the following inequalities,
for some positive constants $C_i$ $(i=1,2)$
which are independent of $\epsilon$:
\[
 \log (9/5)\leq
\log \tr\widetilde{s}^{(\epsilon)}(Q_{\epsilon})
\leq
C_1\cdot\int_{\Delta(Q_{\epsilon},1/2)}
 \log\tr\widetilde{s}^{(\epsilon)}
 \cdot\dvol_g
+C_2.
\]
Recall that $\log \tr \widetilde{s}^{(\epsilon)}$ are 
uniformly bounded from above.
Therefore, there exists a positive constant $C_3$
such that the following holds for any sufficiently small
$\epsilon>0$:
\[
 \int_{\Delta(Q_{\epsilon},1/2)}
 -\min (0,\log\tr\widetilde{s}^{(\epsilon)})
 \cdot\dvol_g
\leq C_3.
\]
Due to Fatou's lemma,
we obtain the following:
\[
 \int_{\Delta(Q_{\infty},1/2)}
 -\min\bigl(0,\log \tr \widetilde{s}^{(\infty)}\bigr)
 \cdot\dvol_{g}\leq C_3.
\]
It means $\widetilde{s}^{(\infty)}$
is not constantly $0$ on $\Delta(Q_{\infty},1/2)$.
In all,
we can conclude that
$\widetilde{s}^{(\infty)}$ is a positive constant multiplication.
Thus, the proof of Lemma \ref{lem;06.1.21.50}
is accomplished.
\hfill\qed

\vspace{.1in}

Let $\bigl\{\widetilde{s}^{(\epsilon_i)}\bigr\}$
be a subsequence as in Lemma \ref{lem;06.1.21.50}.
It is almost everywhere convergent to some constant multiplication.
Then, we obtain that the sequence
$\bigl\{\det \widetilde{s}^{(\epsilon_i)}
=b_{\epsilon_i}^{\rank E}\cdot\id_{\det(E)}\bigr\}$
converges to the positive constant.
In particular, $\{b_{\epsilon_i}\}$ is convergent.
Therefore, the sequence $\bigl\{s^{(\epsilon_i)}\bigr\}$
is convergent to the identity.
Thus we are done.
\hfill\qed

\begin{cor}
 \label{cor;06.1.18.16}
\mbox{{}}
\begin{itemize}
\item
 The sequence $\bigl\{ h^{(\epsilon)} \bigr\}$
 is convergent to $h^{(0)}_0$ weakly in $L_1^2$
 locally on $C-D$. 
\item
The sequence 
$\bigl\{\DDlambda s^{(\epsilon)}\bigr\}$ is weakly convergent to $0$
in $L^2$ locally on $C-D$.
\item
The sequence
$\{\theta^{(\epsilon)}\}$ converges to $\theta^{(0)}$
is weakly convergent to $0$ in $L^2$ locally on $C-D$.
\item
In particular,
the sequences are convergent almost everywhere.
\hfill\qed
\end{itemize}
\end{cor}

\subsection{The surface case}
\label{subsection;06.2.12.1}
\subsubsection{Statement}

Let $X$ be a smooth projective surface
with an ample line bundle $L$,
and let $D$ be a simple normal crossing divisor
with the irreducible decomposition $D=\bigcup_{i\in S}D_i$.
We put $X^{\ast}:=X-D$.
Let $\vecc$ be any element of $\real^S$.
Let $(E,\vecF,\DDlambda)$ be 
a $\mu_L$-stable $\vecc$-parabolic $\lambda$-flat bundle
on $(X,D)$ with trivial characteristic numbers
$\pardeg_{L}(E,\vecF)=\int_X\parch_2(E,\vecF)=0$.
Recall that we have already known $\parchern_1(E,\vecF)=0$
due to Bogomolov-Gieseker inequality and Hodge index theorem
(See Corollary 6.2 of \cite{mochi4}.)
Hence, we can take the pluri-harmonic metric
$h_{\det(E)}$ of the determinant bundle $\det(E,\vecF,\DDlambda)$.
The purpose of this subsection is to show the following  existence theorem.
\begin{thm}
 \label{thm;06.1.18.250}
There exists a tame pluri-harmonic metric $h$ of 
$(E,\DDlambda)_{|X^{\ast}}$
with $\det (h)=h_{\det E}$
which is adapted to the parabolic structure.
\end{thm}

The proof will be given in the rest of this subsection.

\subsubsection{The sequence of Hermitian-Einstein metrics
 for the $\epsilon$-perturbations}

Let $\vecF^{(\epsilon)}$ be an $\epsilon$-perturbation
as in (II) of Subsection \ref{subsubsection;06.1.19.10}.
If $\epsilon$ is sufficiently small,
$(E,\vecF^{(\epsilon)},\DDlambda)$ is also $\mu_L$-stable.
We also have
$\parchern_1(E,\vecF^{(\epsilon)})=\parchern_1(E,\vecF)=0$.
Since $(E,\vecF^{(\epsilon)},\DDlambda)$ is graded semisimple
and satisfies (SPW)-condition,
we can apply Proposition \ref{prop;06.1.18.5}.
Let $h^{(\epsilon)}$ be the Hermitian-Einstein metric
for $(E,\vecF^{(\epsilon)},\DDlambda)$
with respect to $\omega_{\epsilon}$,
such that $\det h^{(\epsilon)}=h_{\det(E)}$
and $\Lambda_{\omega_{\epsilon}}G(h^{(\epsilon)})=0$
(Proposition \ref{prop;06.1.18.5}).

Since $h_{\det(E)}$ is pluri-harmonic,
we also have $\tr G(h^{(\epsilon)})=0$.
Therefore, we have the following convergence:
\begin{equation}
 \label{eq;06.1.18.110}
\left(
\frac{\sqrt{-1}}{2\pi}
\right)^2
 \int \bigl|G(h^{(\epsilon)})
 \bigr|_{h^{(\epsilon)},\omega_{\epsilon}}^2
 \dvol_{\omega_{\epsilon}}
=\left(
\frac{\sqrt{-1}}{2\pi}
\right)^2
\int \tr\Bigl(G(h^{(\epsilon)})^2\Bigr)
=2\bigl(1+|\lambda|^2\bigr)^2\cdot\parch_2(E,\vecF^{(\epsilon)})
\lrarr 0.
\end{equation}
We would like to discuss the limit of $h^{(\epsilon)}$
for $\epsilon\to 0$.

\subsubsection{Convergence on almost every curve}
\label{subsubsection;06.2.9.1}

Let $L^m$ be sufficiently ample.
We put $\proj_m:=\proj\bigl(H^0(X,L^m)^{\lor}\bigr)$.
For any $s\in \proj_m$, we put $X_s:=s^{-1}(0)$.
Recall Proposition \ref{prop;06.1.18.6},
and let $\nbigu$ denote the Zariski open subset of $\proj_m$
which consists of the points $s$
with the following properties:
\begin{itemize}
\item
$X_s$ is smooth,
and $X_s\cap D$ is a simple normal crossing divisor.
\item
$(E,\vecF,\DDlambda)_{|X_s}$ is $\mu_L$-stable.
\end{itemize}
If $\epsilon$ is sufficiently small,
we have $\nbigu\neq \emptyset$.

We will use the notation $X_s^{\ast}:=X_s\setminus D$
and $D_s:=X_s\cap D$.
We have the metric $\omega_{\epsilon,s}$ of $X_s^{\ast}$,
induced by $\omega_{\epsilon}$.
The induced volume form is denoted by $\vol_s$.
We put $(E_s,\vecF_s,\DDlambda_s):=(E,\vecF,\DDlambda)_{|X_s}$.
We have the metric $h^{(\epsilon)}_{|X_s^{\ast}}$
of $E_{s\,|\,X_s^{\ast}}$.
Since
$(E_s,\vecF^{(\epsilon)}_s,\DDlambda_s)$ are also stable
for any point $s\in\nbigu$,
we have the harmonic metric 
$h_{s}^{(\epsilon)}$ of $(E_s,\vecF_s^{(\epsilon)},\DDlambda_s)$
with $\det h_s^{(\epsilon)}=h_{\det E\,|\,X_s^{\ast}}$.
Let $u_s^{(\epsilon)}$ be the endomorphism
of $E_{|X_s^{\ast}}$ determined by
$h^{(\epsilon)}_{|X_s^{\ast}}
=h_{s}^{(\epsilon)}\cdot u_s^{(\epsilon)}$.
For a point $x\in X^{\ast}$,
we put $\nbigu_x:=\{s\in \nbigu\,|\,x\in X_s\}$.
We put $Z:=\{x\in X^{\ast}\,\big|\,\nbigu_x=\emptyset\}$.
We remark that $Z$ is a finite set.
Let us fix a sequence $\epsilon_i\lrarr 0$.
We often use the notation ``$\epsilon$''
instead of ``$\epsilon_i$'',
for simplicity of the description.
Let $\DDlambda_s:=\DDlambda_{|X_s^{\ast}}$.

\begin{lem}
 \label{lem;06.1.18.20}
For almost every $s\in \nbigu$,
the following holds:
\begin{itemize}
\item
We have the following convergence when $\epsilon\lrarr 0$:
\begin{equation}
 \label{eq;06.1.18.7}
 \int_{X_s} \bigl|
 G(h^{(\epsilon)}_{|X_s})
 \bigr|^2_{h^{(\epsilon)}_s,\omega_{\epsilon}}
 \dvol_{s}
\lrarr 0.
\end{equation}
\item
 For each $\epsilon$,
 we have the finiteness:
\begin{equation} 
\label{eq;06.1.26.50}
\bigl\|\DDlambda_{s}
 u_s^{(\epsilon)}
 \bigr\|_{L^2,h_s^{(\epsilon)},\omega_{\epsilon}}<\infty.
\end{equation}
\end{itemize}
Let $\widetilde{\nbigu}$ denote the set of $s$
for which both of {\rm(\ref{eq;06.1.18.7})} 
and {\rm(\ref{eq;06.1.26.50})} hold.
\end{lem}
\pf
It can be shown by the same argument
as the proof of Lemma 9.3 of \cite{mochi4}.
($\nbigz_2$ should be corrected to
 $\bigl\{
 (x,s,t)\in X\times U_1\times\nbigb\,\big|\,
 (ts_2+(1-t)s)(x)=0
 \bigr\}$.)
\hfill\qed

\vspace{.1in}

We obtain the following claims from Proposition 
\ref{prop;06.1.18.15} and Corollary \ref{cor;06.1.18.16}.
\begin{cor}
 \label{cor;06.1.18.40}
For any $s\in\widetilde{\nbigu}$,
the sequence
$\{ h^{(\epsilon)}_{|X^{\ast}_s}\}$ converges to $h^{(0)}_s$
weakly in $L_1^2$ locally on $X_s^{\ast}$,
and $\{\theta^{(\epsilon)}_{|X_s^{\ast}}\}$
converges to $\theta^{(0)}_s$
weakly in $L^2$ locally on $X_s^{\ast}$.
In particular, they are almost everywhere convergent.
\end{cor}
\pf
It follows from Lemma \ref{lem;06.1.18.20}
and Proposition \ref{prop;06.1.18.15}
\hfill\qed

\subsubsection{The construction of a metric defined almost everywhere}

Let us take any Kahler form $\omega_{\proj_m}$ of $\proj_m$.
We put
$\nbigz:=\{(s,x)\in \nbigu\times X^{\ast}\,|\, x\in X_s\}$.
Then, we have the induced metric of $\nbigz$.
The induced volume form is denoted by $\dvol_{\nbigz}$.
Let $\nbigt$ denote the set of
$(s,x)\in \widetilde{\nbigu}\times X$
such that $(s,x)\in \nbigz$ and
$\lim_{\epsilon\to 0}h^{(\epsilon)}_{|x}= h^{(0)}_{s|x}$.

\begin{lem}
 \label{lem;06.1.18.50}
The measure of $\nbigt^c:=\nbigz-\nbigt$ is $0$
with respect to $\dvol_{\nbigz}$.
\end{lem}
\pf
Let us consider the naturally defined fibration $\nbigz\lrarr \nbigu$.
Then, the claim follows from Corollary \ref{cor;06.1.18.40}
and Fubini's theorem.
\hfill\qed

\vspace{.1in}

\begin{lem}
For almost every $x\in X^{\ast}$
and almost every $s\in\nbigu_x$,
the sequence
$\{h_{|x}^{(\epsilon)}\}$ converges to $h^{(0)}_{s\,|\,x}$.
\end{lem}
\pf
Let us consider the naturally defined fibration
$\nbigt\lrarr X^{\ast}$.
Then, the claim follows from Lemma \ref{lem;06.1.18.50}
and Fubini's theorem.
\hfill\qed

\vspace{.1in}

Let $\nbigv$ denote the set of $x\in X^{\ast}$
such that the sequence
$\{h_{|x}^{(\epsilon)}\}$ converges to $h^{(0)}_{s\,|\,x}$
for almost $s\in \nbigu_x$.
For any $x\in\nbigv$,
let $\widetilde{\nbigu}_x$ denote the set of $s$
such that 
$\{h_{|x}^{(\epsilon)}\}$ converges to $h^{(0)}_{s\,|\,x}$.

\begin{lem}
 \label{lem;06.1.18.60}
For any $x\in\nbigv$
and for any $s_i\in\widetilde{\nbigu}_x$ $(i=1,2)$,
we have $h^{(0)}_{s_1\,|\,x}=h^{(0)}_{s_2\,|\,x}$.
\end{lem}
\pf
Both of them are same as
the limit $\lim_{\epsilon\to 0} h_x^{(\epsilon)}$.
\hfill\qed

\vspace{.1in}

Let us take any $x\in \nbigv$
and any $s\in\widetilde{\nbigu}_x$.
Then, the metric $h_x$ of $E_{|x}$
is given by $h_x:=h^{(0)}_{s\,|\,x}$.
Due to Lemma \ref{lem;06.1.18.60}, it is well defined.
Thus, we obtain the metric
$h_{\nbigv}:=(h_x\,|\,x\in\nbigv)$ of $E_{|\nbigv}$.

\subsubsection{The $C^1$-property}
\label{subsubsection;06.1.21.100}

We would like to show that $h_{\nbigv}$ is $C^1$
on $X^{\ast}-Z$,
in other words,
we would like to show the existence of
a $C^1$-metric $h$ of $E_{|X^{\ast}-Z}$
such that $h=h_{\nbigv}$ on $\nbigv$.
Let us begin with a preparation.
\begin{lem}
 \label{lem;06.1.18.70}
Let $x\in X^{\ast}-Z$.
Let us take any $s\in\nbigu_x$.
Then, there exists a Lefschetz fibration
$\varphi:\widetilde{X}\lrarr \proj^1$
with the following properties:
\begin{itemize}
\item
 $x$ is not a singular point of $\varphi$.
\item
 $\varphi^{-1}(0)=X_s$.
\item
 Almost every $t\in\proj^1$ belongs to $\widetilde{\nbigu}$.
\end{itemize}
\end{lem}
\pf
Let $\nbigm$ denote the set of 
the lines $\ell$ of $\proj_m$ which contain $s$.
We put as follows:
\[
\widehat{\proj}_m=\bigl\{
 (\ell,s')\in \nbigm\times\proj_m
 \,\big|\, s'\in\ell  \bigr\}
\subset\nbigm\times\proj_m.
\]
It is the blow up of $\proj_m$ at $s$.
We have the projection $\pi_2:\widehat{\proj}_m\lrarr\proj_m$.
We put $\widehat{\nbigu}:=\pi_2^{-1}(\nbigu)$ and
$\widehat{\widetilde{\nbigu}}:=
   \pi_2^{-1}(\widetilde{\nbigu})$.
Since $\nbigu-\widetilde{\nbigu}$ has measure $0$,
the measure of
$\widehat{\proj}_m-\widehat{\widetilde{\nbigu}}$
is also $0$.
Let us consider the projection
$\pi_1:\widehat{\proj}_m\lrarr\nbigm$,
and apply Fubini's theorem.
Then,
for almost every $\ell\in\nbigm$
and for almost every $s_1\in\ell$,
we have
$s_1\in\widehat{\widetilde{\nbigu}}$.
Thus we are done.
\hfill\qed

\vspace{.1in}

Let $x$ be any point of $X^{\ast}-Z$.
Let us take a Lefschetz fibration
$\pi_i:\widetilde{X}_i\lrarr\proj^1$ $(i=1,2)$
with the following properties:
\begin{itemize}
\item
 Both of them satisfy the properties in Lemma 
 \ref{lem;06.1.18.70}.
\item
 Around $x$,
 the fibers of $\pi_1$ and $\pi_2$ are transversal.
 Then, two fibrations give the holomorphic coordinate
 $(z_1,z_2)$  of an appropriate neighbourhood $U_x$  of $x$,
 such that $\{z_i=a\}=\pi_i^{-1}(a)\cap U_x$.
\end{itemize}

For any $t_i\in\proj^1$,
let $X_{t_i}:=\pi_i^{-1}(t_i)$.
If $t_i$ are close to $0$,
$(E,\vecF,\DDlambda)_{|X_{t_i}}$ are stable,
and hence there exist tame harmonic bundles
$h_{t_i}$ for $(E,\vecF,\DDlambda)_{|X_{t_i}}$
such that $\det(h_{t_i})=h_{\det(E)|X_{t_i}}$.
Let $\theta_{t_i}$ denote the operator
obtained from $\DDlambda_{|X_{t_i}}$ and $h_{t_i}$
as in Subsection \ref{subsubsection;06.1.10.7}.

Let us take an appropriate neighbourhoods
$B_i\subset \proj^1$ of $0$.
Recall Proposition \ref{prop;06.1.18.80}.
Then, $\bigl\{h_{t_1}\,\big|\,t_1\in B_1\bigr\}$
are $C^{\infty}$-along $z_2$,
and it is continuous with respect to $(z_1,z_2)$.
The family $\bigl\{\theta_{t_1}\,\big|\,t_1\in B_1\bigr\}$
has a similar property.
Thus, we obtain a continuous metric $h^{(1)}$
and the continuous section
$\theta^{(1)}$ of $\End(E)\otimes\Omega^{1,0}$ around $x$.
Similarly $\bigl\{h_{t_2}\,\big|\,t_2\in B_2\bigr\}$ is $C^{\infty}$ along $z_1$
and it is continuous with respect to $(z_1,z_2)$.
The family $\{\theta_{t_2}\,\big|\,t_2\in B_2\}$ has a similar property.
Thus, we obtain a continuous metric $h^{(2)}$
and the continuous section $\theta^{(2)}$ of $\End(E)\otimes\Omega^{1,0}$
around $x$.

We remark that $h^{(1)}=h_{\nbigv}=h^{(2)}$
on $U_x\cap \nbigv$ due to our construction of $h_{\nbigv}$.
Since $h^{(i)}$ are continuous,
we obtain $h^{(1)}=h^{(2)}$ on $U_x$.
Then, we obtain that $h^{(i)}$ are $C^{1}$ on $U_x$,
due to the continuity of $\theta^{(i)}$.

Therefore, we obtain the $C^1$-metric $h$ of $E$ on $X^{\ast}-Z$
with the following properties:
\begin{itemize}
\item
 $h_{|\nbigv}=h_{\nbigv}$
\item
 For any $s\in\nbigu$,
 we have $h_{|X_s^{\ast}}=h_{s}$
 and $\theta_{h\,|\,X_s^{\ast}}=\theta_{h_s}$.
\end{itemize}

\subsubsection{Pluri-harmonicity}

We would like to show that $h$ is pluri-harmonic.
By the formalism explained in Subsection
\ref{subsubsection;06.1.10.7},
the operators $\delbar_{h}$ and $\theta_{h}$
are given on $X-(D\cup Z)$ from $h$ and $\DDlambda$.
Let us take any $C^{\infty}$ metric $h'$ of $E$ on $X-D$,
and let $s'$ be the endomorphism
determined by $h=h'\cdot s'$.
Then, $s'$ is $C^1$,
and we have the following relation:
\[
 \delbar_h=\delbar_{h'}
  +\frac{\lambda}{1+|\lambda|^2}
 s^{\prime\,-1}\delta_{h'}''s',
\quad
 \theta_h=\theta_{h'}
  -\frac{\lambda}{1+|\lambda|^2}
 s^{\prime\,-1}\delta_{h'}'s'.
\]
Then, we obtain $\delbar_h\theta_h$
as a distribution:
\[
 \delbar_h\theta_h=\delbar_{h'}\theta_{h'}
-\frac{\lambda}{1+|\lambda|^2}
\delbar_{h'}
\bigl(
 s^{\prime\,-1}\delta_{h'}'s'\bigr)
+
\frac{\lambda}{1+|\lambda|^2}
 \bigl[
s^{\prime\,-1}\delta''_{h'}s',\,\theta_{h'}
\bigr] 
-\left(
\frac{\lambda^2}{1+|\lambda|^2}\right)^2
 \bigl[
 s^{\prime\,-1}\delta''_{h'}s',\,\,
 s^{\prime\,-1}\delta_h's'
 \bigr].
\]
Similarly, we obtain $G(h)$ as a distribution.

\begin{lem}
 \label{lem;06.1.21.200}
$\delbar_h\theta_h=0$.
\end{lem}
\pf
For any point $x\in X^{\ast}-D$,
let us take the holomorphic coordinate $(z_1,z_2)$ as before.
We remark that
the curves
$\{z_i=a\}$ $(i=1,2)$,
$\{z_1+z_2=b\}$,
$\{z_1+\sqrt{-1}z_2=c\}$
can be regarded as parts of $X_{s'}$
for some $s'\in \nbigu$.
We have the expression
$\theta=f_1\cdot dz_1+f_2\cdot dz_2$,
where $f_i$ are continuous sections of $\End(E)$.
We have already known 
$\del f_1/\del \zbar_1
=\del f_2/\del \zbar_2=0$.
Thus, we have only to show
$\del f_i/\del \zbar_j=0$ for $i\neq j$.
Let us consider the change of the coordinate
given by $w_1=z_1+z_2$ and $w_2=z_1-z_2$.
Then, we have the following:
\[
 f_1 \cdot dz_1+f_2 \cdot dz_2
=\frac{1}{2}(f_1+f_2)\cdot dw_1
+\frac{1}{2}(f_1-f_2)\cdot dw_2.
\]
Thus, we obtain the following:
\begin{equation}
 \label{eq;06.1.18.150}
 0=\frac{\del}{\del \overline{w}_1} (f_1+f_2)
=\frac{1}{2}\left(
 \frac{\del}{\del \zbar_1}+\frac{\del}{\del \zbar_2}
 \right)(f_1+f_2)
=\frac{1}{2}\left(
 \frac{\del f_2}{\del \zbar_1}+\frac{\del f_1}{\del \zbar_2}
 \right).
\end{equation}
Let us consider the change of the coordinate
given by $u_1=z_1+\sqrt{-1}z_2$
and $u_2=z_1-\sqrt{-1}z_2$.
Then, we have the following:
\[
 f_1\cdot dz_1+f_2\cdot dz_2=
 \frac{1}{2}
\left(f_1+\frac{1}{\sqrt{-1}}f_2\right)
 du_1
+\frac{1}{2}
\left(f_1-\frac{1}{\sqrt{-1}}f_2\right)
 du_2.
\]
Thus, we obtain the following:
\begin{equation}
 \label{eq;06.1.18.151}
 0=\frac{\del}{\del \overline{u}_1}
 \left(f_1+\frac{1}{\sqrt{-1}}f_2\right)
=\frac{1}{2}\left(
 \frac{\del}{\del \zbar_1}
-\frac{1}{\sqrt{-1}}\frac{\del}{\del \zbar_2}\right)
 \left(f_1+\frac{1}{\sqrt{-1}}f_2\right)
=\frac{1}{2}
 \left(
 \frac{1}{\sqrt{-1}}\frac{\del f_2}{\del \zbar_1}
-\frac{1}{\sqrt{-1}}\frac{\del f_1}{\del \zbar_2}
 \right).
\end{equation}
From (\ref{eq;06.1.18.150}) and (\ref{eq;06.1.18.151}),
we obtain $\del f_i/\del\zbar_j=0$ for $i\neq j$.
Thus, we obtain $\delbar_h\theta_h=0$,
and the proof of Lemma \ref{lem;06.1.21.200}
is accomplished.
\hfill\qed

\vspace{.1in}

\begin{lem}
 \label{lem;06.1.18.100}
$h$ is a harmonic metric for $(E,\DDlambda)$
 with respect to $\omega_0$ on $X^{\ast}-Z$.
(Recall
 $Z=\bigl\{x\in X^{\ast}\,\big|\,\nbigu_x=\emptyset\bigr\}$.)
\end{lem}
\pf
Due to Lemma \ref{lem;06.1.21.200},
we have $\Lambda_{\omega}G(h)=
\Lambda_{\omega}(\delbar_h\theta_h)=0$.
Hence, we have only to show that $h$ is $C^{\infty}$.
We obtain the following formula
in the level of distribution,
by the formalism explained
in Subsection \ref{subsubsection;06.1.18.3}:
\[
 \Delta^{\lambda}_{h',\omega}(s')=
 s'\bigl(-\Lambda_{\omega} G(h')\bigr)
+\sqrt{-1}\Lambda_{\omega} \DDlambda s'\cdot  s^{\prime\,-1}
 \cdot \DDlambdastar_{h'} s'.
\]
The right hand side is $C^0$.
Hence, by using the elliptic regularity and 
the standard boot strapping argument,
we obtain that $s'$ is $C^{\infty}$.
Thus, we obtain Lemma \ref{lem;06.1.18.100}.
\hfill\qed

\begin{lem}
\label{lem;06.1.18.200}
$h$ is pluri-harmonic metric of $E_{|X^{\ast}-Z}$.
\end{lem}
\pf
We have already shown $\delbar_h\theta_h=0$
in Lemma \ref{lem;06.1.21.200}.
Because of Corollary \ref{cor;06.1.21.250},
we have only to show $\theta_h^2=0$.
Due to Corollary \ref{cor;06.1.18.40}
and $\theta_{h\,|\,X_s}=\theta_s$,
we know that
the sequence $\{\theta^{(\epsilon)}\}$ converges
to $\theta_h$ almost everywhere.
In particular,
we obtain the almost everywhere convergence
of $\{\theta^{(\epsilon)\,2}\}$ to $\theta_h^2$.
On the other hand,
we know the almost everywhere convergence 
$G(h^{(\epsilon)})\lrarr 0$,
due to (\ref{eq;06.1.18.110}).
We have
$G(h^{(\epsilon)})=
 \delbar^{(\epsilon)\,2}
+\delbar^{(\epsilon)}\theta^{(\epsilon)}
+\theta^{(\epsilon)\,2}$,
which is the decomposition into 
$(2,0)$, $(1,1)$ and $(0,2)$-forms.
Therefore, we obtain $\theta_h^2=0$,
almost everywhere.
Thus, we obtain Lemma \ref{lem;06.1.18.200}.
\hfill\qed

\begin{lem}
\label{lem;06.1.18.310}
$h$ gives a pluri-harmonic metric of $E_{|X^{\ast}}$.
\end{lem}
\pf
We have only to check that
$h$ gives a $C^{\infty}$-metric of $E_{|X^{\ast}}$.
Let $Q$ be a point of $Z$.
Let $(U,z_1,z_2)$ be a holomorphic coordinate
around $Q$ such that $z_1(Q)=z_2(Q)=0$.
The pluri-harmonic metric $h$ of $(E,\DDlambda)_{|U-\{Q\}}$ is given.
We would like to show that $h$ is naturally
extended to the pluri-harmonic metric of $(E,\DDlambda)_{|U}$.

We have
$\theta=f_1 \cdot dz_1+f_2\cdot dz_2$ defined on $U-\{Q\}$.
Let us consider the characteristic polynomials
$\det (t-f_i)$ for $i=1,2$.
The coefficients are holomorphic on $U-\{Q\}$,
and thus on $U$ due to the theorem of Hartogs.
Hence, the eigenvalues of $f_i$ are bounded on $U$.
Let us consider the restriction of
$(E,\DDlambda,h)$ to the discs 
$C(a_j):=\{z_j=a_j\}$  $(a_j\neq 0)$ 
for $j=1,2$.
Then, it can be shown
that the norms
$\bigl|f_{i\,|\,C(a_j)}\bigr|_{h}\leq C$ $(i\neq j)$ 
can be dominated independently from $a_j$.
(See Lemma 2.7 in \cite{s5}, for example.)
Thus, $f_{i}$ are bounded with respect to $h$
on $U-\{Q\}$.
In other words,
$\theta$ is bounded on $U-\{Q\}$.

Let $E':=E_{|U-\{z_1\cdot z_2=0\}}$.
Let us consider the sheaf
$\prolong{E'}$ on $U$
of the sections satisfying the growth condition
$|g|_h=O\bigl(\prod |z_i|^{-\epsilon}\bigr)$
for any $\epsilon>0$
(Subsection \ref{subsubsection;06.1.26.100}).
By using the result of the asymptotic behaviour
of tame harmonic bundle at $\lambda$
(\cite{mochi2}),
$\prolong{E'}$ is locally free on $U$.
Since $\prolong{E'}$ and $E_{|U-\{Q\}}$
are naturally isomorphic on $U-\{Q\}$,
they are isomorphic on $U$.
Let $h'$ be any $C^{\infty}$-metric of $E_{|U}$,
and let $s'$ be the endomorphism determined by
$h=h'\cdot s'$.
Due to the norm estimate given in \cite{mochi2},
the metrics $h$ and $h'$ are mutually bounded.
Hence, $s'$ and $(s^{\prime})^{-1}$ are bounded on $U$.
Let $\delta_{h'}'$ and $\delta_{h'}''$
be obtained from $\DDlambda$ and $h'$
as in Subsection \ref{subsubsection;06.1.10.7}.
Due to the boundedness of $\theta$,
we have the boundedness of
$(s^{\prime})^{-1}\delta'_{h'} s^{\prime}$
on $U-\{Q\}$.
Due to the boundedness of $\theta^{\dagger}$,
we have the boundedness of
$(s^{\prime})^{-1}\delta_{h'}''s^{\prime}$
on $U-\{Q\}$.
Then, we can deduce that
$s^{\prime\,-1}\DDlambda s^{\prime}$
is also bounded on $U-\{Q\}$.
(See Subsection \ref{subsubsection;06.1.18.3}. 
for example.)
Since we have the formula
$ \Delta^{\lambda}_{h',\omega_0} s'
=s'(-\Lambda_{\omega_0} G(h'))
+\Lambda_{\omega_0} \DDlambda_{h'} s'
 \cdot s^{\prime\,-1}\cdot
 \DDlambdastar_{h'} s'$,
we can conclude that $s'$ is $C^{\infty}$
due to the standard bootstrapping argument.
Namely, $h$ is extended to the $C^{\infty}$-metric of $E_{|U}$.
\hfill\qed

\subsubsection{The end of the proof of Theorem \ref{thm;06.1.18.250}}
\label{subsubsection;06.1.26.150}

Now, we have only to show
that $h$ is tame and adapted to the parabolic structure.
Since $h_{|X_s}=h_s$ for any $s\in\nbigu$,
the tameness immediately follows from
the curve test. (See Proposition \ref{prop;08.2.10.1}.)
Then, we obtain the prolongment
$\widetilde{E}:=\prolongg{\vecc}{E}$
with the induced parabolic structure $\vecF$
(Subsection \ref{subsubsection;06.1.26.100}).
We would like to show that
$(E,\vecF,\DDlambda)$
and $(\Etilde,\vecF,\DDlambda)$
are isomorphic.
For that purpose,
we see that the identity
$E_{|X^{\ast}}\lrarr E_{|X^{\ast}}$
can be prolonged to the homomorphism
$\Psi:E\lrarr \Etilde$.
Let $Q$ be any smooth point of $D_i\subset D$.
We take a holomorphic coordinate
$(U_Q,z_1,z_2)$ with the following property:
\begin{itemize}
\item
 The curve  $z_1^{-1}(0)$ is same as $U_Q\cap D$.
\item
 The curves $C(b):=z_2^{-1}(b)$ are parts of
 $X_{s(b)}$ for $s(b)\in \nbigu$.
\end{itemize}
Let $f$ be a holomorphic section of $E_{|U}$.
Since the restriction $h_{|X_{s(b)}}$ is the same as
$h_{s(b)}$,
we have $|f_{|C(b)}|_h=O(|z_1|^{-c_i-\epsilon})$ 
for any $\epsilon>0$.
Then, we obtain
$|f|_h=O(|z_1|^{-c_i-\epsilon})$ for any $\epsilon>0$,
due to the result given in \cite{mochi2}.
Thus, $f$ naturally gives the section of
$\widetilde{E}$ on $U$.
Therefore,
we obtain the morphism
$E\lrarr \widetilde{E}$
on $X-\bigl(\cup_{i\neq j}D_i\cap D_j\bigr)$.
It is naturally extended to the morphism
$E\lrarr \widetilde{E}$.

Recall that the restriction of 
$\widetilde{E}=\prolongg{\vecc}{E}(h)$
to $X_s$ is same as
$\prolongg{\vecc}{(E_{|X_s})}(h_{s})$. (See \cite{mochi2}.)
Therefore, 
the restrictions of $\Psi$ to $X_s$ are isomorphic,
by construction.
Hence, $\Psi$ is isomorphic
on $X-\bigl(\bigcup_{i\neq j}D_i\cap D_j\bigr)$,
and thus on $X$.
By a similar argument,
we can show that the parabolic structures 
are also same.
Thus, the proof of Theorem \ref{thm;06.1.18.250}
is finished.
\hfill\qed

\subsection{Correspondences}

\subsubsection{Kobayashi-Hitchin correspondence
  in the higher dimensional case}

Let $X$ be a smooth projective variety of dimension $n$ $(n\geq 3)$
with an ample line bundle $L$,
and let $D$ be a simple normal crossing divisor
with the irreducible decomposition $D=\bigcup_{i\in S}D_i$.
Let $(\vecE_{\ast},\DDlambda)$ be 
a $\mu_L$-stable regular filtered
$\lambda$-flat bundle on $(X,D)$ in codimension two
with trivial characteristic numbers
$\pardeg_{L}(\vecE_{\ast})=\int_X\parch_{2,L}(\vecE_{\ast})=0$,
and we put $(E,\DDlambda):=(\vecE_{\ast},\DDlambda)_{|X-D}$.
Recall $\parchern_1(\vecE_{\ast})=0$ due to
the Bogomolov-Gieseker inequality and 
the Hodge index theorem.
For each $\vecc\in\real^S$,
we have the determinant line bundle
$\det(\prolongg{\vecc}{E})$ of torsion-free sheaf
$\prolongg{\vecc}{E}$,
on which we have the induced parabolic structure
and the induced flat $\lambda$-connection.
Thus, we obtain the canonically determined
regular filtered $\lambda$-flat bundle
$\bigl(\det\vecE_{\ast},\DDlambda\bigr)$ on $(X,D)$
of rank one.
We also have
$\parchern_1\bigl(\det\vecE_{\ast}\bigr)
=\parchern_1\bigl(\vecE_{\ast}\bigr)=0$.
Therefore, we can take a pluri-harmonic metric 
$h_{\det E}$ of $(\det(E),\DDlambda)$
which is adapted to the parabolic structure
of $\det \vecE_{\ast}$.
By the assumption,
we have a subset $Z\subset D$
with $\codim_X(Z)\geq 3$
such that $(\vecE_{\ast},\DDlambda)_{|X-Z}$
is a regular filtered $\lambda$-flat {\em bundle}.

\begin{thm}
 \label{thm;06.1.23.100}
There exists the unique tame pluri-harmonic metric $h$
of $(E,\DDlambda)$ with the following properties:
\begin{itemize}
\item $\det(h)=h_{\det E}$.
\item It is adapted to the parabolic structure of
 $\vecE_{\ast}$ on $X-Z$.
 Namely, 
 $(\vecE_{\ast}(h),\DDlambda)_{|X-Z}
 \simeq (\vecE_{\ast},\DDlambda)_{|X-Z}$,
where $(\vecE_{\ast}(h),\DDlambda)$ denotes
 the regular filtered $\lambda$-flat bundle on $(X,D)$
 obtained from $(E,\DDlambda,h)$.
 (See Subsection {\rm\ref{subsection;06.2.4.50}}.)
\end{itemize}
\end{thm}
\pf
Due to Mehta-Ramanathan type theorem 
(Proposition \ref{prop;06.1.18.6}),
the uniqueness can be easily reduced
to the $\dim X=1$ case,
by considering the restriction to the generic 
curves $C\subset X$.
We have already known it 
(Proposition \ref{prop;06.1.23.50}).

We will use the induction on the dimension $n$
to show the existence.
The case $n=2$ has already been shown
(Theorem \ref{thm;06.1.18.250}).
Assume that $L^m$ is sufficiently ample.
We put $\proj_m:=\proj(H^0(X,L^m)^{\lor})$.
For any $s\in \proj_m$, we put $X_s:=s^{-1}(0)$.
Recall Proposition \ref{prop;06.1.18.6}.
Let $\nbigu$ be the Zariski open subset of $\proj_m$
which consists of $s\in\proj_m$ with the following properties:
\begin{itemize}
\item
$X_s$ is smooth,
 and $D_s:=X_s\cap D$ is a normal crossing divisor.
\item
 The codimension of $Z\cap X_s$ in $X_s$
 is larger than $3$.
\item
$(\vecE,\DDlambda)_{|X_s}$ is $\mu_L$-stable.
\end{itemize}

We use the existence hypothesis
in the $(n-1)$-dimensional case of the induction.
Then, we may have the tame pluri-harmonic metric
$h_s$ of $(E,\DDlambda)_{|X_s\setminus D}$
with $\det (h_s)=h_{\det E\,|\,X_s\setminus D}$
which is adapted to the parabolic structure on $X_s\setminus W$.
We also use the uniqueness result in the $(n-2)$-dimensional case.
Then, we can show the existence of
a finite subset $Z'\subset X-D$ and a metric $h$ of $E_{|X-D}$
such that $h_{s\,|\,P}=h_{|P}$.
By the arguments given in Subsections
\ref{subsubsection;06.1.21.100}--\ref{subsubsection;06.1.26.150},
we can show that $h$ is the desired metric.
The only different point is 
the argument to show the vanishing of $G(h)=0$.
Due to $\dim X_s\geq 2$,
it can be shown easier.
\hfill\qed

\vspace{.1in}

\begin{thm}
\label{thm;06.2.4.500}
Let $X$, $D$ and $L$ be as above.
Let $(\vecE_{\ast},\DDlambda)$ be 
a saturated $\mu_L$-stable regular filtered $\lambda$-flat sheaf
on $(X,D)$ with the trivial characteristic numbers
$\pardeg_L(\vecE_{\ast})=\int_X\parch_{2,L}(\vecE_{\ast})=0$.
We put $(E,\DDlambda):=(\vecE_{\ast},\DDlambda)_{|X-D}$.
Then, there exists a pluri-harmonic metric $h$
of $(E,\DDlambda)$ such that
the induced regular filtered $\lambda$-flat bundle
$\bigl(\vecE_{\ast}(h),\DDlambda\bigr)$
is isomorphic to $(\vecE_{\ast},\DDlambda)$.
Such a metric is unique
up to positive constant multiplication.
In particular,
$\vecE_{\ast}$ is a filtered bundle.
\end{thm}
\pf
Since a saturated regular filtered $\lambda$-flat sheaf
is a regular filtered $\lambda$-flat bundle in codimension two
(Lemma \ref{lem;06.2.4.3}),
we may apply Theorem \ref{thm;06.1.23.100}.
Then, there exists a pluri-harmonic metric $h$
and a subset $W\subset D$ with $\codim_X(W)\geq 3$
such that the induced regular filtered $\lambda$-flat
bundle $(\vecE_{\ast}(h),\DDlambda)$ is isomorphic 
to $(\vecE_{\ast},\DDlambda)$ on $X-W$.
Since both of 
$(\vecE_{\ast}(h),\DDlambda)$ and 
$(\vecE_{\ast},\DDlambda)$ are saturated,
they are isomorphic on $X$.
\hfill\qed

\subsubsection{The equivalence of the categories}

Let $\Cpoly{\lambda}$ denote the category 
of $\mu_L$-stable regular filtered $\lambda$-flat bundles 
$(\vecE_{\ast},\DDlambda)$ on $(X,D)$
with the trivial characteristic numbers $\pardeg_L(\vecE_{\ast})=
 \int_X\parch_{2,L}(\vecE_{\ast})=0$.
Morphisms $f:(\vecE_{1\,\ast},\DDlambda_1)\lrarr
 (\vecE_{2\,\ast},\DDlambda_2)$ are defined to be
$\nbigo_X$-homomorphism
$f:\vecE_1\lrarr\vecE_2$  satisfying
$\DDlambda_2\circ f=f\circ \DDlambda_1$
and $f\bigl(\prolongg{\vecc}{E}_1\bigr)\subset \prolongg{\vecc}{E}_2$
for any $\vecc$.

\begin{cor}
 \label{cor;06.1.23.101}
Let $\lambda_i$ $(i=1,2)$ be two complex numbers.
We have the natural functor
$\Xi_{\lambda_1,\lambda_2}:
 \Cpoly{\lambda_1}\lrarr \Cpoly{\lambda_2}$,
which is equivalent.
It preserves direct sums, tensor products and duals.
\end{cor}
\pf
Let $(\vecE^{\lambda_1}_{\ast},\DD^{\lambda_1})$
be an object of $\Cpoly{\lambda_1}$.
We put $E^{\lambda_1}:=\vecE^{\lambda_1}_{|D}$.
We have a pluri-harmonic metric $h$
of $(E^{\lambda_1},\DD^{\lambda_1})$,
which is adapted to the parabolic structure.
Then, we obtain the operators
$\delbar_h,\del_h,\theta_h,\theta_h^{\dagger}$,
as in Subsection \ref{subsubsection;06.1.10.7}.
Note that the holomorphic structure of
$E^{\lambda_1}$ is given by
$\delbar_h+\lambda_1\theta^{\dagger}_h$.
The $(0,1)$-operator $\delbar_h+\lambda_2\theta^{\dagger}_h$
also gives a holomorphic structure of
$C^{\infty}$-bundle $E^{\lambda_1}$.
To distinguish them,
we use the notation $E^{\lambda_2}$,
when we consider the holomorphic structure
$\delbar_h+\lambda_2\theta^{\dagger}_h$.
We put 
$\DD^{\lambda_2}:=\delbar_h+\theta_h
+\lambda_2(\del_h+\theta^{\dagger}_h)$,
which gives a flat $\lambda_2$-connection
of $E^{\lambda_2}$.
The metric $h$ is pluri-harmonic for
$(E^{\lambda_2},\DD^{\lambda_2})$.
Since the corresponding Higgs bundle
for $(E^{\lambda_1},\DD^{\lambda_1},h)$ and
$(E^{\lambda_2},\DD^{\lambda_2},h)$
are same,
we obtain the tameness of $(E^{\lambda_2},\DD^{\lambda_2},h)$.
Therefore, we obtain the prolongment
$(\vecE^{\lambda_2},\DDlambda)$,
which are $\mu_L$-polystable regular filtered
$\lambda_2$-flat bundle on $(X,D)$
with trivial characteristic numbers
(Proposition \ref{prop;06.1.23.60}).

We remark that 
$(\vecE^{\lambda_2},\DD^{\lambda_2})$
is independent of a choice of  $h$,
due to the uniqueness in Proposition \ref{prop;06.1.23.50}.
Therefore, we put
$\Xi_{\lambda_1,\lambda_2}(\vecE^{\lambda_1},\DD^{\lambda_1}):=
 (\vecE^{\lambda_2},\DD^{\lambda_2})$.
It is easy to see that $\Xi_{\lambda_1,\lambda_2}$
gives a functor.
It is also easy to see that
$\Xi_{\lambda_2,\lambda_1}\circ\Xi_{\lambda_1,\lambda_2}
 (\vecE^{\lambda_1},\DD^{\lambda_1})$
is naturally isomorphic to
$(\vecE^{\lambda_1},\DD^{\lambda_1})$.
The compatibility with the direct sums,
duals and tensor products are obtained
from the corresponding compatibility statements
of the prolongments for tame harmonic bundles
(\cite{mochi2}).
We also remark that
the categories are semisimple.
Thus, we have only to compare the objects.
\hfill\qed

\begin{rem}
From a $\lambda_1$-connection $\DD^{\lambda_1}=d''+d'$,
a $\lambda_2$-connection is given 
$d''+(\lambda_2/\lambda_1)\cdot d'$.
Hence, we have the obvious functor
$\Obv:\Cpoly{\lambda_1}\lrarr \Cpoly{\lambda_2}$.
This is not same as the above functor $\Xi_{\lambda_1,\lambda_2}$.
\hfill\qed
\end{rem}

%% file: 6.tex
\subsection{Definition}

\subsubsection{Filtered structure}

Let $X$ be a complex manifold,
and let $D$ be a simple normal crossing divisor 
with the irreducible decomposition $D=\bigcup_{i\in S} D_i$.
We will use the notation
$D^{[2]}:=\bigcup_{i\neq j}D_i\cap D_j$
and $D_i^{\circ}:=D_i\setminus \bigcup_{j\neq i}D_j$.
Let $\nbigl$ be a local system on $X-D$.
A filtered structure of $\nbigl$ at $D$
is a tuple of increasing filtrations
$\lefttop{i}\nbigf$ $(i\in S)$
of $\nbigl_{|U_i\setminus D}$
indexed by $\real$,
where $U_i$ denotes an appropriate open neighbourhood
of $D_i$.
Let $U_i'$ be an open neighbourhood of $D_i$
such that $U_i'\subset U_i$,
then we have the induced filtration
$\lefttop{i}\nbigf_{|U_i'}$,
and the filtration $\lefttop{i}\nbigf$ can be reconstructed
from $\lefttop{i}\nbigf_{|U_i'}$.
Hence, we define two filtered structures
$(\lefttop{i}\nbigf,U_i\,|\,i\in S)$
and $(\lefttop{i}\nbigf',U_i'\,|\,i\in S)$ are equivalent,
if there exists an open neighbourhood $U_i''$
of $D_i$ such that $U_i''\subset U_i\cap U_i'$
and $\lefttop{i}\nbigf_{|U_i''}=\lefttop{i}\nbigf'_{|U_i''}$.
A local system $\nbigl$
equipped with an equivalence class of filtered structures
$(\lefttop{i}\nbigf,U_i)$ 
is called a filtered local system,
and it is denoted by $\nbigl_{\ast}$.
We do not have to care
about a choice of open neighbourhoods $U_i$.

Morphisms of filtered local systems
$f:\nbigl_{1\,\ast}\lrarr \nbigl_{2\,\ast}$
are defined to be a morphism
$f:\nbigl_1\lrarr\nbigl_2$ of local systems
preserving the filtered structures
in an obvious sense.
We denote by $\widetilde{\nbigc}(X,D)$ the category
of filtered local systems on $(X,D)$.

\subsubsection{Characteristic numbers}

We put $U_i^{\ast}:=U_i\setminus D$ and 
$\lefttop{i}\Gr_a^{\nbigf}(\nbigl_{|U_i^{\ast}}):=
 \lefttop{i}\nbigf_a(\nbigl_{|U_i^{\ast}}) \big/
 \lefttop{i}\nbigf_{<a}(\nbigl_{|U_i^{\ast}})$.
Since the local monodromy around $D_i$
preserves the filtration $\lefttop{i}\nbigf$,
we obtain the induced endomorphism
of $\lefttop{i}\Gr_a^{\nbigf}(\nbigl_{|U_i^{\ast}})$,
and thus the generalized eigen decomposition:
\[
 \lefttop{i}\Gr_a^{\nbigf}(\nbigl_{|U_i^{\ast}})
=\bigoplus_{\omega}
 \lefttop{i}\Gr_{(a,\omega)}^{\nbigf,\EE}(\nbigl_{|U_i^{\ast}}).
\]
We put as follows:
\[
 \Par\bigl(\nbigl_{\ast},i\bigr):=
 \bigl\{
 a\in\real\,\big|\,
 \lefttop{i}\Gr^{\nbigf}_a
 \bigl(\nbigl_{|U_i^{\ast}}\bigr)\neq 0
 \bigr\},
\quad
  \KMS\bigl(\nbigl_{\ast},i\bigr):=
 \bigl\{
 (a,\omega)\in\real\times\cnum^{\ast}\,\big|\,
 \lefttop{i}\Gr^{\nbigf,\EE}_{(a,\omega)}
 \bigl(\nbigl_{|U_i^{\ast}}\bigr)\neq 0
 \bigr\}.
\]

The parabolic first Chern class is defined 
as follows:
\begin{equation}
 \parchern_1(\nbigl_{\ast}):=
 -\sum_{i\in S}\wt(\nbigl_{\ast},i)\cdot [D_i]\in H^2(X,\real),
\quad\quad
 \wt(\nbigl_{\ast},i):=
 \sum_{a\in\Par(\nbigl_{\ast},i)}
 a\cdot\rank \lefttop{i}\Gr_a^{\nbigf}(\nbigl_{|U_i^{\ast}}).
\end{equation}
Here $[D_i]$ denotes the cohomology class
representing $D_i$.

Let $\Irr(D_i\cap D_j)$ denote the set
of the irreducible components of $D_i\cap D_j$.
For each $P\in \Irr(D_i\cap D_j)$,
let $U_P$ be an appropriate open neighbourhood of $P$
in $X$ such that $U_P\subset U_i\cap U_j$.
We put $U_P^{\ast}:=U_P\setminus D$.
We have the two filtrations
$\lefttop{i}\nbigf$ and $\lefttop{j}\nbigf$
of $\nbigl_{|U_P^{\ast}}$.
The naturally induced graded local system
is denoted as follows:
\[
 \lefttop{P}\Gr^{\nbigf}(\nbigl_{|U_P^{\ast}})
=\bigoplus_{(a_i,a_j)\in\real^2}
 \lefttop{P}\Gr^{\nbigf}_{(a_i,a_j)}(\nbigl_{|U_P^{\ast}}),
\quad
 \lefttop{P}\Gr^{\nbigf}_{(a_i,a_j)}(\nbigl_{|U_P^{\ast}})
:=\frac{\lefttop{i}\nbigf_{a_i}\cap \lefttop{j}\nbigf_{a_j}}
 {\sum_{(b_i,b_j)\lneq (a_i,a_j)}
 \lefttop{i}\nbigf_{b_i}\cap \lefttop{j}\nbigf_{b_j}}.
\]
Here $(b_i,b_j)\lneq (a_i,a_j)$ means
``$b_i\leq a_i$, $b_j\leq a_j$ and $(b_i,b_j)\neq (a_i,a_j)$''.
We have the two endomorphisms 
induced by the local monodromies around
$U_P\cap D_i$ and $U_P\cap D_j$,
which are commutative.
Hence, we obtain the generalized eigen decomposition:
\[
 \lefttop{P}\Gr^{\nbigf}_{\veca}(\nbigl_{|U_P^{\ast}})
=\bigoplus_{\vecomega\in\cnum^{\ast\,2}}
 \lefttop{P}\Gr^{\nbigf,\EE}_{\veca,\vecomega}
 (\nbigl_{|U_P^{\ast}}).
\]
We put  as follows:
\[
 Par(\nbigl_{\ast},P):=
 \bigl\{
 (a_i,a_j)\in\real^2\,\big|\,
 \lefttop{P}\Gr^{\nbigf}_{(a_i,a_j)}(\nbigl_{|U_P^{\ast}})
\neq 0 \bigr\},
\]
\[
  \KMS(\nbigl_{\ast},P):=
 \bigl\{
 (\veca,\vecomega)\in\real^2\times\cnum^{\ast\,2}
 \,\big|\,
 \lefttop{P}\Gr^{\nbigf,\EE}_{(\veca,\vecomega)}
 (\nbigl_{|U_P^{\ast}})\neq 0
 \bigr\}.
\]
The parabolic second Chern character
is defined as follows:
\begin{multline}
 \parch_2(\nbigl_{\ast}):=
 \frac{1}{2}\sum_{i\in S}\sum_{a\in\Par(\nbigl_{\ast},i)}
 a^2\cdot\rank\lefttop{i}\Gr_a^{\nbigf}(\nbigl)\cdot [D_i]^2 \\
+\frac{1}{2}\sum_{i\in S}
 \sum_{j\neq i}\sum_{P\in \Irr(D_i\cap D_j)}
 \sum_{(a_i,a_j)\in\Par(\nbigl_{\ast},P)}
 a_i\cdot a_j\cdot\rank\lefttop{P}\Gr^{\nbigf}_{(a_i,a_j)}
 \bigl(\nbigl_{|U_P^{\ast}}\bigr)\cdot[P].
\end{multline}

When $X$ is a smooth projective variety
with an ample line bundle $L$,
we put as follows:
\[
 \pardeg_L(\nbigl_{\ast}):=\int_X 
\parchern_1(\nbigl_{\ast})\cdot c_1(L)^{\dim X-1},
\quad
 \mu_L(\nbigl_{\ast}):=
 \frac{\pardeg_L(\nbigl_{\ast})}{\rank \nbigl}.
\]
Then, the notion of $\mu_L$-stability,
$\mu_L$-semistability,
and $\mu_L$-polystability for filtered local systems
on $(X,D)$ are defined in the standard manner.
We also put as follows:
\[
  \int_X\parchern^2_{1,L}(\nbigl_{\ast}):=
 \int_X\parchern_1(\nbigl_{\ast})^2\cdot c_1(L)^{\dim X-2},
\quad
 \int_X\parch_{2,L}(\nbigl_{\ast}):=
 \int_X\parch_{2,L}(\nbigl_{\ast})\cdot c_1(L)^{\dim X-2}.
\]

\subsection{Correspondence}

In this subsection,
we give the correspondence of
filtered local systems on $(X,D)$ and 
saturated regular filtered $\lambda$-flat sheaves
($\lambda\neq 0$).
See Subsection \ref{subsubsection;06.2.5.10}
for saturated regular filtered $\lambda$-flat sheaves.
Since we have the obvious equivalence
between flat $\lambda$-connection
and flat $1$-connection,
we only discuss the case $\lambda=1$,
i.e. ordinary flat connections.

Let $\nbigc_{1}^{sat}(X,D)$ denote the category
of saturated regular filtered flat sheaves on $(X,D)$.
Let us see briefly that we have the functor
$\Phi:\widetilde{\nbigc}(X,D)\lrarr \nbigc_1^{sat}(X,D)$
which gives the equivalence.
Since it is given by Simpson in \cite{s2}
essentially in the curve case,
we give only an outline.

\subsubsection{Construction of $\Phi$}
\label{subsubsection;06.2.3.30}

First, we give a construction of $\Phi$.
Let $\nbigl_{\ast}$ be a filtered local system on $(X,D)$.
Let $(E,\nabla)$ be the corresponding flat bundle
on $X-D$.
We have the Deligne extension
$(\widetilde{E},\nabla)$ on $(X,D)$.
We put $\vecE:=\widetilde{E}\otimes\nbigo(\ast D)$.
Thus, we have only to give the way of the construction
of the $\nbigo_X$-coherent submodules
$\prolongg{\veca}{E}\subset \vecE$
such that $\nabla \prolongg{\veca}{E}
 \subset\prolongg{\veca}{E}\otimes\Omega^{1,0}(\log D)$
and $\bigcup_{\veca\in\real^S}\prolongg{\veca}{E}=\vecE$.
Let us consider the case
$X=\Delta^n=\{(z_1,\ldots,z_n)\,|\,|z_i|<1\}$
and $D=\{z_1=0\}$.
Then, the construction is essentially same
as that for the case $\dim X=1$
given by Simpson \cite{s2}.
We briefly recall it.
Let $H(\nbigl)$ denote the space of the multi-valued
flat sections of $\nbigl$.
We have the induced filtration
$\nbigf H(\nbigl)$ and
the generalized eigen decomposition
$H(\nbigl)=\bigoplus_{\omega}\EE_{\omega}(H(\nbigl))$,
which are compatible in the sense
$\nbigf_a=\bigoplus_{\omega}\nbigf_a\cap \EE_{\omega}$.
Let $\vecu=(u_1,\ldots,u_r)$ be
a frame compatible of $H(\nbigl)$,
compatible with $(\nbigf,\EE)$.
Then, for each $u_i$,
the numbers $\omega(u_i)\in\cnum^{\ast}$ and $a(u_i)\in\real$
are determined by
$u_i\in\EE_{\omega(u_i)}$
and $u_i\in \nbigf_{a(u_i)}- \nbigf_{<a(u_i)}$.
The complex number $\alpha(u_i)$ is determined 
by the conditions
$\exp(-2\pi\alpha(u_i))=\omega(u_i)$
and $0\leq \Re \alpha(u_i)<1$.
Let $M^u$ denote the endomorphism of $H(\nbigl)$
or $\nbigl$,
which is the unipotent part of the monodromy around $D$,
and we put $N:=-(2\pi\sqrt{-1})^{-1}\log M^u$.
We regard $u_i$ as a multi-valued $C^{\infty}$-section
of $E$.
Then, it is standard that
$v_i:=\exp\bigl(\log z_1(\alpha(u_i)+N)\bigr)\cdot u_i$
gives a holomorphic section of $E$.
Moreover, $\vecv=(v_1,\ldots,v_r)$ gives
a frame of the Deligne extension $\widetilde{E}$.
Let $b$ be any real number.
Then, we put
$n(b,u_i):=\max\bigl\{
 n\in\seisuu\,\big|\,
 a(u_i)-\Re\alpha(u_i)+n\leq b
 \bigr\}$,
and we put
$v_i(b):=z_1^{-n(b,u_i)}\cdot v_i$.
Let $\prolongg{b}{E}$ denote the $\nbigo_X$-submodule
of $\vecE$ generated by $v_1(b),\ldots,v_r(b)$.
It is easy to check that
$\prolongg{b}{E}$ is locally free and
independent of a choice of $\vecu$.
It is also easy to see 
$\vecE=\bigcup_{b\in\real} \prolongg{b}{E}$.
Thus, we obtain the filtration
in the case $X=\Delta^n$ and $D=\{z_1=0\}$.
It can be checked that the filtration
is independent of a choice of
the coordinate $(z_1,z_2,\ldots,z_n)$ satisfying $D=\{z_1=0\}$.

For any $\vecb\in\real^S$,
we obtain $\prolongg{\vecb}{E}$ on $X-D^{[2]}$
by gluing them.
The subsheaves $\prolongg{\vecb}{E}$
are determined by the condition (\ref{eq;06.2.3.2}).
\begin{lem}
$\prolongg{\vecb}{E}$ is a coherent
$\nbigo_X$-module.
Hence, we obtain the saturated regular filtered
flat sheaf $(\vecE_{\ast},\nabla)$
on $(X,D)$.
\end{lem}
\pf
We may assume that 
$X=\Delta^n$ and 
$D=\bigcup_{i=1}^{\ell}\{z_i=0\}$.
Let $H(\nbigl)$ denote the space of
the multi-valued flat sections of $\nbigl$.
We have the monodromy endomorphisms
$M_i$ $(i=1,\ldots,\ell)$ along the loop
around $D_i$ with counter clockwise direction.
They induce the decomposition
\begin{equation}
 \label{eq;08.2.7.150}
H(\nbigl)=
 \bigoplus_{\vecomega\in\cnum^{\ell}}
 \EE_{\vecomega}H(\nbigl),
\end{equation}
where each $\EE_{\vecomega}H(\nbigl)$
is preserved by $M_i$ $(i=1,\ldots,\ell)$,
and the eigenvalues of $M_i$ on 
$\EE_{\vecomega}H(\nbigl)$ are $\omega_i$.
We also have the filtrations
$\lefttop{i}\nbigf$ $(i=1,\ldots,\ell)$ of $H(\nbigl)$,
corresponding to the divisor $D_i$.
Each $\lefttop{i}\nbigf$ is compatible
with the decomposition (\ref{eq;08.2.7.150}).

Fix $j$ such that $1\leq j\leq \ell$.
We take a frame 
$\vecu=(u_1,\ldots,u_r)$
of $H(\nbigl)$ compatible with 
the filtration $\lefttop{j}\nbigf$ 
and the decomposition (\ref{eq;08.2.7.150}).
For each $u_p$,
the tuple $\vecomega(u_p)\in\cnum^{\ell}$
is determined by $u_p\in \EE_{\omega}$.
Let $\alpha_i(u_p)\in\cnum$ $(i=1,\ldots,\ell)$
be determined by
$\exp\bigl(-2\pi\alpha_i(u_p)\bigr)
=\omega_i(u_p)$
and $0\leq \Re \alpha_i(u_p)<1$.
We also have the numbers
$a(u_p)\in\real$ such that 
$u_p\in
 \lefttop{j}\nbigf_{a(u_p)}
-\lefttop{j}\nbigf_{<a(u_p)}$.
We put
$n(b_j,u_p):=\max\bigl\{
 n\in\seisuu\,\big|\,
 a_j(u_p)-\Re\alpha_j(u_p)+n\leq b_j
 \bigr\}$.
Let $N_i:=-(2\pi\sqrt{-1})^{-1}\log M^u$ 
$(i=1,\ldots,\ell)$,
where $N_i$ denotes the logarithm 
of the unipotent part of $M_i$.
We take a sufficiently large integer $I$.
Then, we put as follows:
\[
 v_p:=
 z_j^{n(b_j,u_p)}\cdot
 \prod_{i\neq j} z_i^{I}
 \prod_{i=1}^{\ell}
 \exp\bigl(
 \log z_i\cdot
 (\alpha_i(u_p)+N_i)
 \bigr)\cdot u_p
\]
If $I$ is sufficiently large,
$v_p$ gives the section of
$\prolongg{\vecb}{E}$ on $X$.
By the correspondence,
we obtain the following morphism,
for $j=1,\ldots,\ell$:
\[
 \Phi_j:
 \bigoplus_{p=1}^r\nbigo_X\cdot v_p
\lrarr
 \prolongg{\vecb}{E}
\]
The morphisms $\Phi_j$ $(j=1,\ldots,\ell)$
induce the morphism
$\Phi:\nbigo^{\oplus \ell\cdot r}
\lrarr \prolongg{\vecb}{E}$.
The image of $\Phi$ is $\nbigo_X$-coherent,
and it is the same as $\prolongg{\vecb}{E}$
on $X-D^{[2]}$.
Then, it is easy to show that
$\prolongg{\vecb}{E}$ is the same as
the double dual of the image of $\Phi$
which is $\nbigo_X$-coherent.
\hfill\qed

\vspace{.1in}

Let $f:\nbigl_{1\,\ast}\lrarr\nbigl_{2\,\ast}$ be a morphism.
Let $(\vecE_{i\,\ast},\nabla_i):=\Phi(\nbigl_i)$.
We have the induced map
$\widetilde{f}:\vecE_1\lrarr\vecE_2$.
It is easy to see that
$\prolongg{\vecc}{E}_{1\,|\,X-D^{[2]}}
\lrarr \prolongg{\vecc}{E}_{2\,|\,X-D^{[2]}}$
is induced.
Due to saturatedness of $(\vecE_{2\,\ast},\nabla)$,
we obtain maps
$\prolongg{\vecc}{E}_{1}
\lrarr \prolongg{\vecc}{E}_{2}$, and thus
$\Phi(f):(\vecE_{1\,\ast},\nabla_1)\lrarr
 (\vecE_{2\,\ast},\nabla_2)$.

\subsubsection{Equivalence}

Let us show that $\Phi$ is equivalent.
To begin with, 
we consider the case $X=\Delta^n$ and $D=\{z_1=0\}$.
Let $\nbigc_1^{vb}(X,D)$ denote the category
of regular filtered flat bundles on $(X,D)$,
which is the subcategory of $\nbigc_1^{sat}(X,D)$.
By the construction,
the image of $\Phi$ is contained in
$\nbigc_1^{vb}(X,D)$.
The following lemma can be shown as in \cite{s2}.
\begin{lem}
 \label{lem;06.2.3.50}
The functor $\Phi$ gives the equivalence of
$\widetilde{\nbigc}_1(X,D)$ and $\nbigc^{vb}_1(X,D)$.
It is also compatible
with direct sums, duals, and tensor products.
\hfill\qed
\end{lem}

\begin{lem}
 \label{lem;06.2.4.100}
In the case $X=\Delta^n$ and $D=\{z_1=0\}$,
we have $\nbigc_1^{vb}(X,D)\simeq \nbigc_1^{sat}(X,D)$
naturally.
In particular,
$\Phi$ gives the equivalence 
$\widetilde{\nbigc}_1(X,D)\simeq \nbigc_1^{sat}(X,D)$.
\end{lem}
\pf
Let $(\vecE_{\ast},\nabla)$ be a saturated regular 
filtered flat sheaf on $(X,D)$.
We put $(E,\nabla):=(\vecE_{\ast},\nabla)_{|X-D}$,
and let $\nbigl$ denote the corresponding local system on $X-D$.
Let $H(\nbigl)$ denote the space of the multi-valued flat
sections of $\nbigl$.

Recall that there exists a subset $W\subset D$
with $\codim_X(W)\geq 3$ such that
$(\vecE_{\ast},\nabla)_{|X-W}$ is regular filtered flat bundle
on $(X-W,D-W)$ (Lemma \ref{lem;06.2.4.3}).
Let $P$ be any point of $D-W$,
and let $(U_P,z_1,\ldots,z_n)$ be a holomorphic coordinate
neighbourhood such that $z_1^{-1}(0)=U_P\cap D$
and $U_P\cap W=\emptyset$.
Due to Lemma \ref{lem;06.2.3.50},
we have the unique filtration $\nbigf$ of 
$H(\nbigl_{|U_P\setminus D})\simeq
 H(\nbigl)$
corresponding to $(\vecE_{\ast},\nabla)_{|U_P}$.
Due to the uniqueness,
it is independent of a choice of $P$ and $U_P$.

Let $\vecu=(u_1,\ldots,u_r)$ be a frame of $H(\nbigl)$
compatible with the filtration $\nbigf$
and the generalized eigen decomposition
with respect to the monodromy around $D$.
For any real number $b\in\real$,
we construct $\vecv(b)=\bigl(v_1(b),\ldots,v_r(b)\bigr)$
as above.
Then, for any $P\in D-W$,
$\vecv(b)$ gives a holomorphic frame of
$\prolongg{b}{E}_{|U_P}$ compatible with the filtration
due to Lemma \ref{lem;06.2.4.100}.
Hence, each $v_i(b)$ gives a section of
$\prolongg{b}{E}_{|X-W}$.
Due to the saturatedness of $(\vecE_{\ast},\nabla)$,
$v_i(b)$ gives a section of $\prolongg{b}{E}$ on $X$.
Now it is easy to see that $\vecv(b)$ gives 
a frame of $\prolongg{b}{E}$,
and in particular, $\prolongg{b}{E}$ is locally free.
Hence, $(\vecE_{\ast},\nabla)$ is a regular filtered 
flat bundle on $(X,D)$.
\hfill\qed

\vspace{.1in}

Now, it is easy to see that $\Phi$ is equivalent
for general $(X,D)$.
Let us see the fully faithfulness of $\Phi$.
The faithfulness is obvious.
Let $f:\Phi(\nbigl_{1\,\ast})\lrarr \Phi(\nbigl_{2\,\ast})$
be a morphism in $\nbigc_{1}^{sat}(X,D)$.
We have the map
$g:\nbigl_1\lrarr\nbigl_2$ corresponding to $f$.
We would like to check that $g$ preserves
 the filtrations $\lefttop{i}\nbigf$.
Let $P$ be any point of $D_i^{\circ}$,
and $(U,z_1,\ldots,z_n)$ be any coordinate neighbourhood
such that $U\cap D=z_1^{-1}(0)$.
Applying Lemma \ref{lem;06.2.4.100},
we obtain that $g$ preserves the filtration $\lefttop{i}\nbigf$
on $U\setminus D_i$.
Thus, we obtain the fully faithfulness.

Let us show the essential surjectivity.
Let $(\vecE_{\ast},\nabla)$ be a saturated filtered flat sheaf
on $(X,D)$.
Let $\nbigl$ denote the local system
corresponding to $(\vecE_{\ast},\nabla)_{|X-D}$.
We have only to construct the appropriate filtrations
$\lefttop{i}\nbigf$ of $\nbigl_{|U_i\setminus D}$
on appropriate neighbourhoods of $D_i$.
Let $P$ be any point of $D_i^{\circ}$,
and $(U_P,z_1,\ldots,z_n)$ denote any coordinate
neighbourhood around $P$ such that
$z_1^{-1}(0)=U_P\cap D$.
Due to Lemma \ref{lem;06.2.3.50},
we obtain the unique filtration
$\lefttop{i}\nbigf$ of $\nbigl_{|U_P\setminus D}$.
We obtain the filtration $\lefttop{i}\nbigf$
on $\bigcup_{P\in D_i^{\circ}}U_P$
by gluing them, due to the uniqueness.
Thus, we obtain that $\Phi$ is essentially surjective,
and hence equivalent.

\subsubsection{The parabolic first Chern class}

We have the $\seisuu$-action on
$\real\times\cnum$ given by
$n\cdot (a,\alpha)=(a+n,\alpha-n)$.
It induces the action of $\seisuu$ on $\KMS(\vecE_{\ast},i)$.
The following lemma is clear from the construction
of $\Phi$.
\begin{lem}
 \label{lem;06.2.5.30}
We have the bijective correspondence
of the sets
$ \KMS(\Phi(\nbigl_{\ast}),i)/\seisuu$
and 
$\KMS(\nbigl_{\ast},i)$,
which is given by
 $(a,\alpha)
\longmapsto
  (b,\omega)=
 \Bigl(a+\Re\alpha,\exp\bigl(-2\pi\sqrt{-1}\alpha\bigr)
 \Bigr)$
for $(a,\alpha)\in\KMS(\Phi(\nbigl_{\ast}),i)$.
Moreover,
$\rank \lefttop{i}\Gr^{F,\EE}_{(a,\alpha)}
=\rank \lefttop{i}\Gr^{\nbigf,\EE}_{(b,\omega)}$.
\hfill\qed
\end{lem}

\begin{cor}
 \label{cor;06.2.3.100}
We have the equality of the parabolic first Chern class
$\parchern_1(\nbigl_{\ast})
=\parchern_1(\Phi(\nbigl_{\ast}))$.
In particular,
when $X$ is a smooth projective variety
with an ample line bundle $L$,
the $\mu_L$-stability of $\nbigl_{\ast}$
and $\mu_L$-stability of $\Phi(\nbigl_{\ast})$
are equivalent.
\end{cor}
\pf
Recall Lemma \ref{lem;08.2.7.105}.
It is shown for the case
where $(\vecE_{\ast},\nabla)$ is graded semisimple
and $\dim X$ is two dimensional.
However, the graded semisimplicity condition
is not necessary as is explained in Remark \ref{rem;06.1.20.50}.
The assumption $\dim X=2$ is also not necessary,
due to the Lefschetz theorem.
Then, the claim of the corollary
follows from Lemma \ref{rem;06.1.20.50}
and the correspondence of the KMS-spectrums
given in Lemma \ref{lem;06.2.5.30}.
\hfill\qed

\subsubsection{The second parabolic Chern character}

\begin{lem}
Let $X=\Delta^n=\{(z_1,\ldots,z_n)\,|\,|z_i|<1\}$,
and $D=D_1\cup D_2$,
where $D_i=\{z_i=0\}$.
Let $(\vecE_{\ast},\nabla)$ be
a saturated regular filtered flat sheaf on $(X,D)$.
\begin{itemize}
\item
$(\vecE_{\ast},\nabla)$ is a regular filtered flat bundle on $(X,D)$.
\item
Let $\vecc$ be any element of $\real^2$,
and let $\prolongg{\vecc}{E}$ denote the $\vecc$-truncation.
Let $\nbigl_{\ast}$ be the corresponding
filtered local system on $(X,D)$.
Then, we have the equality:
\[
 \rank \lefttop{\nibar}
 \Gr^{\nbigf,\EE}_{(\vecb,\vecomega)} (\nbigl)
=\rank\lefttop{\nibar}
 \Gr^{F,\EE}_{(\veca,\vecalpha)}
 (\prolongg{\vecc}{E}).
\]
Here the meaning of the notation is as follows:
\begin{itemize}
 \item
 $\vecb=(b_1,b_2)$ 
 and $\vecomega=(\omega_1,\omega_2)$
 denote elements of
 $\real^2$ and $\cnum^{\ast\,2}$
 respectively.
 \item
 $\veca=(a_1,a_2)$ and $\vecalpha=(\alpha_1,\alpha_2)$
 denote elements of $\real^2$ and $\cnum^2$
 respectively,
 determined by the conditions
 $c_i-1<a_i\leq c_i$,
 $\exp(-2\pi\sqrt{-1}\alpha_i)=\omega_i$
 and $a_i+\Re\alpha_i=b_i$.
\end{itemize}
\end{itemize}
\end{lem}
\pf
Let $\nbigl_{\ast}=(\nbigl,\lefttop{1}\nbigf,\lefttop{2}\nbigf)$
be as above.
Let $\vecu$ be a frame of $H(\nbigl)$
compatible with 
the filtrations $\lefttop{k}\nbigf$ $(k=1,2)$
and the generalized eigen decompositions of
$H(\nbigl)$.
For each $u_j$ and the divisor $D_k$,
the complex number $\alpha_k(u_j)$
and $a_k(u_j)$ are determined as before.
For the monodromies around $D_k$,
we obtain the nilpotent endomorphism $N_k$
as before.
The holomorphic section $v_j$ is given by
$ v_j:=\exp\Bigl(
 \sum \log z_k\bigl(\alpha_k(u_j)+N_k\bigr)
 \Bigr)$.
Let $n_k(u_j)$ be the numbers determined
by the condition
$c_k-1<n_k(u_j)+a_k(u_j)-\Re\alpha_k(u_j)\leq c_k$.
We put 
$\widetilde{v}_j:=\prod z_k^{-n_k(u_j)}\cdot v_j$.
Then, $\widetilde{\vecv}=(\widetilde{v}_1,\ldots,\widetilde{v}_r)$
gives the frame of 
$\prolongg{\vecc}{E}_{|X-(D_1\cap D_2)}$.
Due to the saturatedness,
$\widetilde{\vecv}=(\widetilde{v}_1,\ldots,\widetilde{v}_r)$
gives the frame of $\prolongg{\vecc}{E}$,
and hence $\prolongg{\vecc}{E}$ are locally free.
Thus, the first claim is proved.
The frame $\widetilde{\vecv}$ is
compatible with $\lefttop{i}\EE$ and $\lefttop{i}F$,
and we have
$\lefttop{k}\deg^F(\widetilde{v}_j)
=a_k(u_j)-\Re\alpha_k(u_j)+n_k(u_j)$
and 
$\widetilde{v}_{j\,|\,D_k}
\in \lefttop{k}\EE(\alpha_k(u_j)-n_k(u_j))$.
Thus, the second claim follows.
\hfill\qed

\begin{cor}
 \label{cor;06.2.3.200}
Let $X$ be a projective manifold with an ample line bundle $L$,
and let $D$ be a simple normal crossing divisor.
Let $(\vecE_{\ast},\nabla)$ be a saturated
regular filtered flat sheaf on $(X,D)$,
and let $\nbigl_{\ast}$ denotes the corresponding
filtered local system.
Then, we have the equality of
the parabolic second Chern character numbers
$\int_X\parch_{2,L}(\nbigl_{\ast})
=\int_{X}\parch_{2,L}(\vecE_{\ast})$.
\hfill\qed
\end{cor}

\begin{cor}
 \label{cor;06.2.5.2}
Let $X$ be a smooth projective variety
with an ample line bundle $L$,
and let $D$ be a simple normal crossing divisor.
Let $\nbigl_{\ast}$ be a $\mu_L$-stable
filtered local system on $(X,D)$.
Then, the Bogomolov-Gieseker inequality for $\nbigl_{\ast}$
holds:
\[
 \int_X\parch_{2,L}(\nbigl_{\ast})
\leq
 \frac{\int_X\parchern_{1,L}^2(\nbigl_{\ast})}
 {2\rank \nbigl}.
\]
\end{cor}
\pf
Recall that saturated regular filtered flat shaves
are regular filtered flat bundles in codimension two
(Lemma \ref{lem;06.2.4.3}).
Hence, the claim follows from Corollary \ref{cor;06.2.3.100},
Corollary \ref{cor;06.2.3.200}
and Corollary \ref{cor;06.1.13.150}.
\hfill\qed

\begin{cor}
Let $X$ be a smooth projective variety
with an ample line bundle $L$,
and let $D$ be a simple normal crossing divisor.
Let $\nbigc_1^{poly}$ be the category of
$\mu_L$-polystable regular filtered flat bundle
on $(X,D)$
with trivial characteristic numbers,
and let $\widetilde{\nbigc}_1^{poly}$
be the category of $\mu_L$-polystable filtered
local system on $(X,D)$
with trivial characteristic numbers.
Then, the functor $\Phi$ naturally gives
the equivalence of them.
\end{cor}
\pf
We have only to remark that
saturated $\mu_L$-stable regular filtered flat sheaves
with trivial characteristic numbers
are regular filtered bundles
(Theorem \ref{thm;06.2.4.500}).
\hfill\qed

\begin{rem}
Due to the result in {\rm\cite{mochi2}}
and the existence of a pluri-harmonic metric
for $\Phi(\nbigl_{\ast})$,
the filtrations $\lefttop{i}\nbigf$
for $\mu_L$-stable filtered local systems $\nbigl_{\ast}$
satisfy some compatibility
around the intersection points of $D$.
\hfill\qed
\end{rem}